\documentclass{amsproc}
\usepackage{tikz}
\usepackage{amsmath}
\usepackage{amssymb}
\usepackage{amsthm}
\usepackage{stmaryrd}
\usepackage{cases}
\usepackage{graphicx}
\usepackage{upgreek}

\usepackage{xcolor}
\definecolor{MyGrey}{rgb}{.804,.804,.756}
\usepackage[colorlinks=false,linkbordercolor=MyGrey,citebordercolor=MyGrey,
urlbordercolor=MyGrey, pdfauthor={Carter, Lebed, Yang}, pdftitle={A prismatic classifying space}]{hyperref}

\newtheorem{theorem}{Theorem}

\theoremstyle{definition}

\theoremstyle{remark}

\numberwithin{equation}{section}

\newcommand*\lt{\mathrel{\triangleleft}}

\newcommand{\Z}{\mathbb{Z}}
\newcommand{\R}{\mathbb{R}}

\newcommand{\SL}{{\mbox{\rm SL}}}
\newcommand{\SO}{{\mbox{\rm SO}}}
\newcommand{\SU}{{\mbox{\rm SU}}}

\begin{document}

\vspace*{1.5pc}\large

\title{Amusing Permutation Representations of Group Extensions}

\author{Yongju Bae}\address{Kyungpook National University\\ Daegu, Korea}
\email{ybae@knu.ac.kr}
\author{J. Scott Carter}
\address{University of South Alabama \\ Mobile AL}
\curraddr{Austin, TX 78728}
\email{carter@southalabama.edu}

\author{Byeorhi Kim}
\address{CRT, POSTECH, Pohang, Korea
}
\email{byeorhikim@postech.ac.kr}
\thanks{}
 
\subjclass[2010]{Primary 20F36, 20B30, 20D99; Secondary 57M25 \\ Corresponding Author; J. Scott Carter, Austin, TX 78728, {\sc EMAIL:} carter@southalabama.edu}

\date{}

\begin{abstract}
Semi-direct products of finite groups have permutation representations that are constructed from the permutation representations of their constituents. One can envision these in a metaphoric sense in which a rope is made from a bundle of threads. In this way, subgroups and quotients are easily visualized. The general idea is applied to the finite subgroups of the special unitary group of $(2\times 2)$-matrices.  Amusing diagrams are developed that describe the unit quaternions, the binary tetrahedral, octahedral, and icosahedral group as well as the dicyclic groups. In all cases, the quotients as subgroups of the permutation group are readily apparent. These permutation representations lead to injective homomorphisms into semi-direct products.

\end{abstract}

\maketitle

\section{Introduction}
\label{S:Intro}

The purpose of this paper is to have some fun with the binary extensions of the dihedral and polyhedral groups. 

The study of finite groups is extremely gratifying to a mathematical novice who  has an eye towards symmetry. Even though in the modern setting, we prefer to think of groups in abstracto, our perception is often facilitated when a given group is represented as a specific set of symmetries. For example, the dihedral groups are usually first introduced as the set of symmetries of a polygon before their abstract definition via, say, a group presentation is given. Braid groups (although they are not finite) and finite permutation groups are often conveniently represented via string diagrams. Group theory as a gateway to the study of other more advanced mathematics can both be approached from a visual or an algebraic point of view. Here the emphasis is upon the visual as well as the recreational.

Our use of the word ``visual" is not meant to exclude sight-impaired mathematicians. More generally, ``visual'' can be replaced by ``tactile." Yet within this paper,  diagrams are drawn to represent elements in specific groups. And the group multiplication between two group elements is achieved by the juxtaposition of their diagrams. The authors imagine that the same purposes can be achieved in a tactile realm. 

A  special case of the Krasner-Kaloujine Theorem (\cite{Evens:Coho}, page 47 contains a sketch of a proof), allows us to represent 
the elements of the finite subgroups of $\SU(2)=S^3$ by means of string diagrams that are projections of geometric braids. 
The 
diagrammatic representations
have analogues for arbitrary semi-direct products of finite groups as we describe in Section~\ref{SS:perm} and specifically in Fig.~4, but the scope of the discussion will be limited towards more specific examples.

\begin{center}
\rule{3in}{0.005in}
\end{center}

In the human endeavor of spinning fibers to make thread, long fibers are combined and overlapped until isolated long filaments appear. Of course, each filament is manufactured as a conglomerate, and the conception of this conglomerate as a single entity is only metaphoric. Imagine, for example a sturdy nautical rope. It is thought of as filaments wound around a central core to form a bundle of filaments that are twisted further. In much the same way, we can see the groups that we present here as a bundle of permutations. Unfortunately words such as bundle, grouping,  etc. have technical meanings in mathematics and within this paragraph we mean none of these things. 

Precisely, we partition groups into sets of cosets, and study the group actions on its coset space. In the finite group case, subgroups can be thought of as permutations. These are spun together to make stronger ropes. 

\begin{center}
\rule{3in}{0.005in}
\end{center}

The focus here is upon examples. In particular, a lot details are given in both algebraic  and diagrammatic 
contexts. Here is how material will be presented. While it has naught to do with the rest of the study,  in the next two paragraphs, the initial point of the study is 
described to give these ideas some 
more context. The material summarized there is the starting point of the paper \cite{BCK1}. Section~\ref{S:wreath} gives an overview
of the techniques. We discuss semi-direct products in which the  groups involved are subgroups of permutation 
groups. In this way, a broad class of examples is possible. A (finite) group action upon a set of cosets
allows a variety of  string diagram depictions. Section~\ref{S:Pr2} gives  an elementary  proof of  the known result Theorem~\ref{two}. 
Then one-by-one, the quaternions (Section~\ref{S:Quat}), dicyclic groups (Section~\ref{S:Dic}), binary tetrahedral 
(Section~\ref{S:BinTet}),  and octahedral groups (Section~\ref{S:BinOct}) are given. Generators for the binary icosahedral group are described via the same diagrammatic method in Section~\ref{S:BinIco}. The depictions closely resembles the techniques 
that we learned from \cite{Kauff:OnKnots} in which polyhedra or polygons are suspended from a ceiling with a band 
that supports twisting.  But here the underlying
polygonal object is replaced by an abstract point set. Also, many other groups are easily visualized via string-like diagrams.

While these diagrammatic techniques were being developed, a fairly large number of our colleagues have become intrigued. From several points of view, the paper may be considered elementary. In particular, not much material beyond a course in group theory is assumed of the reader. However, the techniques are useful, and in many instances we have provided much more information than what one typically finds in a mathematical paper. For example, more than one set of diagrams for the elements in both the binary tetrahedral and the binary octahedral groups are presented. And rather than presenting only the generators, we either depict half or all of the elements in one format or the other. When only half are illustrated the other half can be obtained by juxtaposing the diagram that represents $(-1)$.  The paper is not written with an expert, but rather with a novice, in mind.  In this sense, the paper performs a service. Our hope is that others will find the methods and diagrams useful.

\begin{center}
\rule{3in}{0.005in}
\end{center}

Even though the next two paragraphs are the only place in which quandles are mentioned, let us tell you where
 the main idea originated. It is a small but beautiful outgrowth of a study about the non-trivial cocycle extension 
 of the four element tetrahedral quandle which is known as $QS_4$ or $Q(4,1)$.{\footnote{A quandle is a set $Q$ 
 that has a binary operation $\lt$ that satisfies (i) $(\forall a \in Q)$ $(a\lt a=a)$; (ii) $(\forall a,b \in Q)$ $(\exists! 
 c\in Q)(c\lt a=b)$; (iii) $(\forall a,b,c \in Q)$ $((a\lt b)\lt c= (a\lt c)\lt(b \lt c)).$}} The most concise description uses 
 the  field $\Z/2[t]/(t^2+t+1)$ as the underlying set with $a\lt b = t a + (1+t)b$ as the quandle operation. The  
 elements can be labeled $0$, $1$, $2\leftrightarrow t$, and $3 \leftrightarrow (t+1)$. In this way, they are 
 represented by the coefficients (00,01,10,11)  of the affine expressions and these pairs of coefficients corresponds 
 to an integer  written in binary notation. The function that takes values in $\Z/2$ $$f(a,b)= \left\{ \begin{array}{lr} 0 & {\mbox{\rm if}}   \ a=b, 
 a=3, \ {\mbox{\rm or}} \  b=3, \\ 1 & {\mbox{\rm otherwise}} \end{array} \right.$$ satisfies the quandle cocycle 
 conditions $$f(a,b)+f(a\lt b,c)-f(a,c)-f(a\lt c, b\lt c)=0, \quad f(a,a)=0.$$ So there is a quandle extension $Q(8,1)$ 
 that maps surjectively onto $Q(4,1)$. The inner automorphism group ${\mbox{\rm Inn}}(Q(4,1))$ is generated by
 the $3$-cycles $(1,2,3)$, $(0,3,2)$, $(0,1,3)$, and $(0,2,1)$ since these correspond respectively to the right actions
 of the elements $0$,$1$, $2$, and $3$ ($y:  x\lt y \mapsfrom x$). Meanwhile, the inner automorphism group, $
 {\mbox{\rm Inn}}(Q(8,1))$, is the binary tetrahedral group. It is generated by the following four elements each of
 which is a product of disjoint $3$-cycles $(1,6,3)(2,7,5)$, $(2,4,3)(0,7,6)$, $(0,5,3)(1,7,4)$, $(0,2,1)(4,6,5)$   in the 
  group of symmetries of the set $\{0,1,2,3,4,5,6,7\}$. 

The details of the preceding paragraph are given more explicitly in \cite{BCK1}. Here we hope that an interested reader will take the time out to verify these claims as an exercise in quandle theory. The rest of the paper does not depend upon quandle theory in any respect. 

\begin{center}
\rule{3in}{0.005in}
\end{center}

Finite subgroups of $\SU(2)=S^3$ will be considered.  The $3$-sphere is represented as $$S^3=\{x+y {\boldsymbol i}+z {\boldsymbol j}+w {\boldsymbol k}: x^2+y^2+z^2+w^2 =1 \}.$$The subgroups consist of the {\it unit quaternions} 
$$Q_8=\{\pm1,\pm {\boldsymbol i}, \pm {\boldsymbol j}, \pm {\boldsymbol k}\}, $$ the dicyclic groups  $${\mbox{\rm Dic}}_{n}= \langle \rho,x: \rho^{2n}=1, \ x^2=\rho^n, \  \rho x= x\rho^{-1} \rangle,$$
 the {\it binary tetrahedral group}
$$\widetilde{A_4}
=\langle a,b:  (ab)^2=a^3=b^3 \rangle,$$
the {\it binary octahedral group}
$$\widetilde{\Sigma_4}= \langle a,f: a^3=f^4=(af)^2 \rangle,$$
and the {\it binary icosahedral group}
$$\widetilde{A_5}= \langle a,t: (at)^2=a^3=t^5 \rangle.$$

We have been informed (anonymously and without sources) that the  following theorem is known.

\begin{theorem}\label{main}
There  are injective homomorphisms of the finite subgroups of $\SU(2)$ into wreath products as follows.
\begin{enumerate}
\item $Q_8\subset (\Z/2)^4 \rtimes K_4$ where $K_4$ denotes the Klein $4$-group ($\approx \Z/2 \times \Z/2)$;
\item  $Q_8\subset (\Z/4)^2 \rtimes \Z/2;$
\item ${\mbox{\rm Dic}}_n \subset (Z/4)^n \rtimes \Sigma_n;$
\item ${\mbox{\rm Dic}}_n \subset (Z/n)^2 \rtimes \Sigma_2;$
\item $\widetilde{A_4} \subset (\Z/2)^4 \rtimes A_4;$
\item  $\widetilde{A_4} \approx Q_8 \rtimes \Z/3;$
\item  $\widetilde{A_4} \subset (\Z/6)^4 \rtimes A_4; $
\item ${\mbox{\rm GL}}_2(\Z/3) \subset (\Z/2)^4 \rtimes \Sigma_4;$
\item  $\widetilde{\Sigma_4} \subset (Q_8)^3 \rtimes \Sigma_3;$
\item  $\widetilde{\Sigma_4} \subset \left({\mbox{\rm Dic}}_4\right)^3\rtimes \Sigma_3;$
\item  $\widetilde{\Sigma_4} \subset \left({\mbox{\rm Dic}}_3\right)^4\rtimes \Sigma_4;$
\item $\widetilde{A_5} \subset (\widetilde{A_4})^5 \rtimes A_5.$
\item $\widetilde{A_5} \subset (\Z/10)^{12} \rtimes \Sigma_{12}.$
\end{enumerate}
\end{theorem}
In each case, the group on the right side of the semi-direct product acts on the groups on the left by permuting coordinates. Here is a sketch of the proof.

Section~\ref{S:wreath} demonstrates how to represent elements in the semi-direct products of subgroups of permutation groups via string diagrams. The broader technique  is exemplified for permutation groups $\Sigma_3$ and $\Sigma_4$, the Klein $4$ group $K_4$, the cyclic group $\Z/4$, and the group of symmetries of the square (that we call $H_2$). Then for each of the finite subgroups of $S^3$, permutation representations are developed. In some cases, all the elements are listed via diagrams. In others, only the generators are given. Nonetheless, once a permutation representation is given, a pointer to the group embedding articulated in the theorem will be given.

Let us encapsulate this result in a more concise way. The following theorem is a special case of the Krasner-Kaloujine Theorem.{\footnote{Thanks to Roger Alperin, David Benson, and Greg Kuperberg.}}  To our knowledge its use in the current context is novel. A proof of  Theorem~\ref{two} appears in Section~\ref{S:Pr2}.

\begin{theorem}\label{two}
  Let $G$ denote a finite group of order $nk$. Let $H$ denote a subgroup of order $k$. Then there is an inclusion $G \subset H^n \rtimes \Sigma_n$ where the second factor permutes  the coordinates of  $H^n$. \end{theorem}

\section{Overview of the technique}
\label{S:wreath}

\begin{center}
\includegraphics[width=0.6\paperwidth]{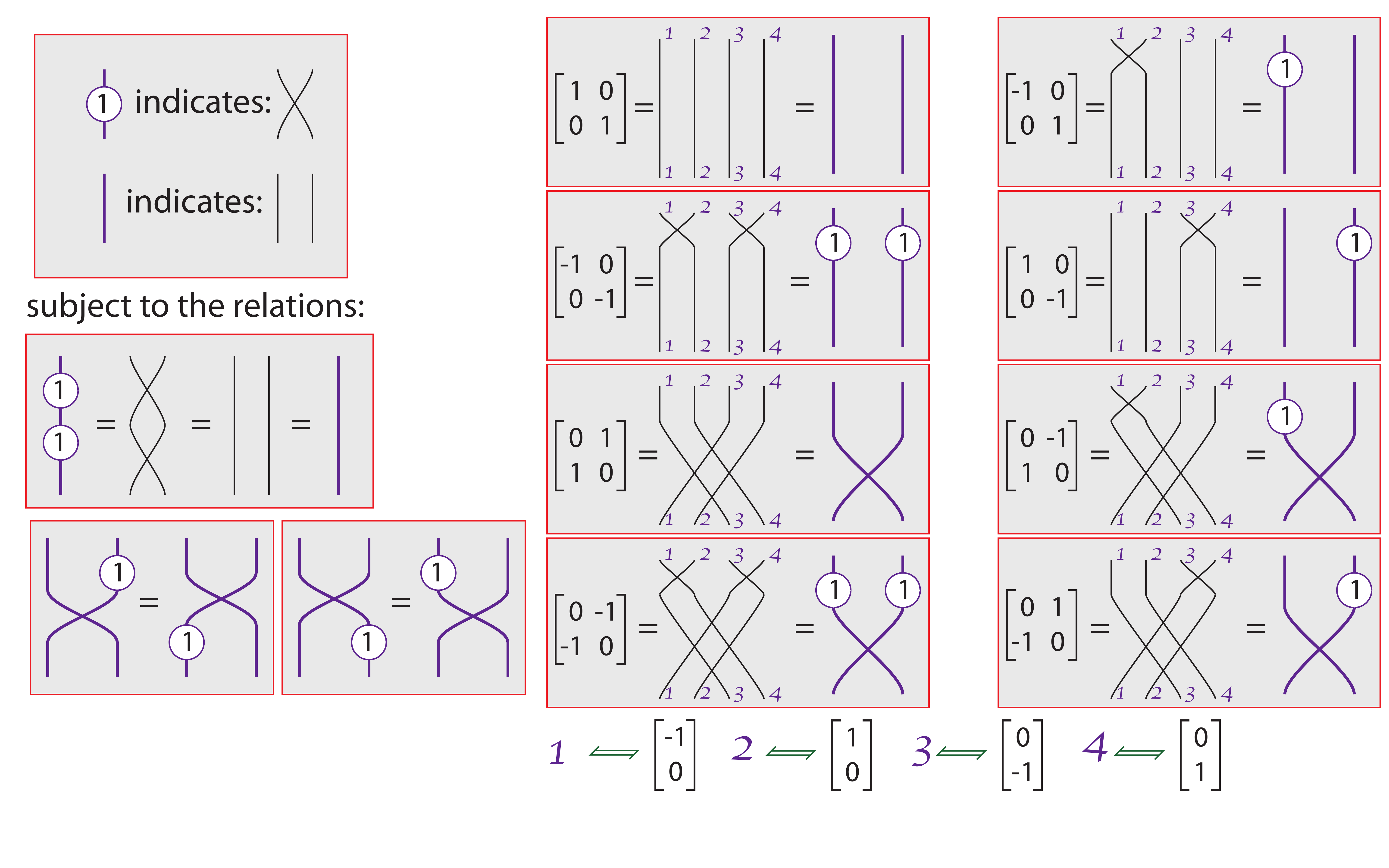}

Figure 1. The permutation representation of the $(2 \times 2)$ signed permutations
\end{center}

An example is presented to initiate the discussion. The group of $(2\times 2)$ signed permutation is the set of eight elements
$$H_2=\left\{ \left[ \begin{array}{lr} \pm 1 & 0 \\ 0 & \pm 1 \end{array} \right], \left[ \begin{array}{lr} 0 & \pm 1  \\ \pm 
1 & 0  \end{array} \right] \right\}$$ in which any combination of signs is allowed. It acts upon the end points 
$\{ (\pm 1, 0)^t,(0,\pm 1)^t\}$ of the coordinate arcs $\{(x,y) \in \R^2: -1 \le x,y \le 1;  \ \ xy=0\}${\footnote{As is 
customary, we will often drop the transpose when discussing these vectors.}} via left matrix
multiplication on the associated column vector. This group coincides with the dihedral group of order $8$. Label $1 \leftrightarrow (-1, 0),$ $2\leftrightarrow (1,0),$ $3\leftrightarrow (0,-1)$, and $4\leftrightarrow (0,1)$. 
Then since \begin{eqnarray*} \left [   \begin{array}{lr} \epsilon_1 & 0 \\ 0 & \epsilon_2 \end{array} \right]  \cdot  \left [   \begin{array}{lr} \delta_1 & 0 \\ 0 & \delta_2 \end{array} \right]  &= & 
\left [   \begin{array}{lr} \epsilon_1 \delta_1 & 0 \\ 0 & \epsilon_2 \delta_2 \end{array} \right], {\mbox{\rm and}}   \\
\left [   \begin{array}{lr} 0 & \epsilon_1  \\  \epsilon_2 & 0 \end{array} \right]  \cdot  \left [   \begin{array}{lr}  \delta_1 & 0 \\ 0 & \delta_2 \end{array} \right]  &= & 
\left [   \begin{array}{lr} 0&  \epsilon_1 \delta_2  \\  \epsilon_2 \delta_1 & 0 \end{array} \right]  \end{eqnarray*}
for $\epsilon_i, \delta_j \in \{ \pm 1\}, (i,j=1,2)$,  the permutation representation of $H_2$ are as follows.

$$\begin{array}{lcllcl}
\left [   \begin{array}{lr}  1 & 0  \\   0 & 1 \end{array} \right] & \leftrightarrow& (1)(2)(3)(4); \ \ \ \ &
\left [   \begin{array}{lr}  -1 & 0  \\   0 & 1 \end{array} \right] & \leftrightarrow &(12); \\
 \left [   \begin{array}{lr}  -1 & 0  \\   0 & -1 \end{array} \right] & \leftrightarrow &(12)(34); \ \ \ \ &
\left [   \begin{array}{lr}  1 & 0  \\   0 & -1 \end{array} \right] & \leftrightarrow& (34); \\
\left [   \begin{array}{lr}   0 & 1 \\ 1 & 0     \end{array} \right] & \leftrightarrow& (13)(24); \ \ \ \ &
\left [   \begin{array}{lr}   0 & -1 \\ 1 & 0 \end{array} \right] & \leftrightarrow& (1324); \\
\left [   \begin{array}{lr}   0 & -1 \\ -1 & 0 \end{array} \right] & \leftrightarrow& (14)(23)\ \ \ \ &
\left [   \begin{array}{lr}   0 & 1 \\ -1 & 0 \end{array} \right] & \leftrightarrow &(1423); 
\end{array}$$

\begin{center}
\includegraphics[width=0.6\paperwidth]{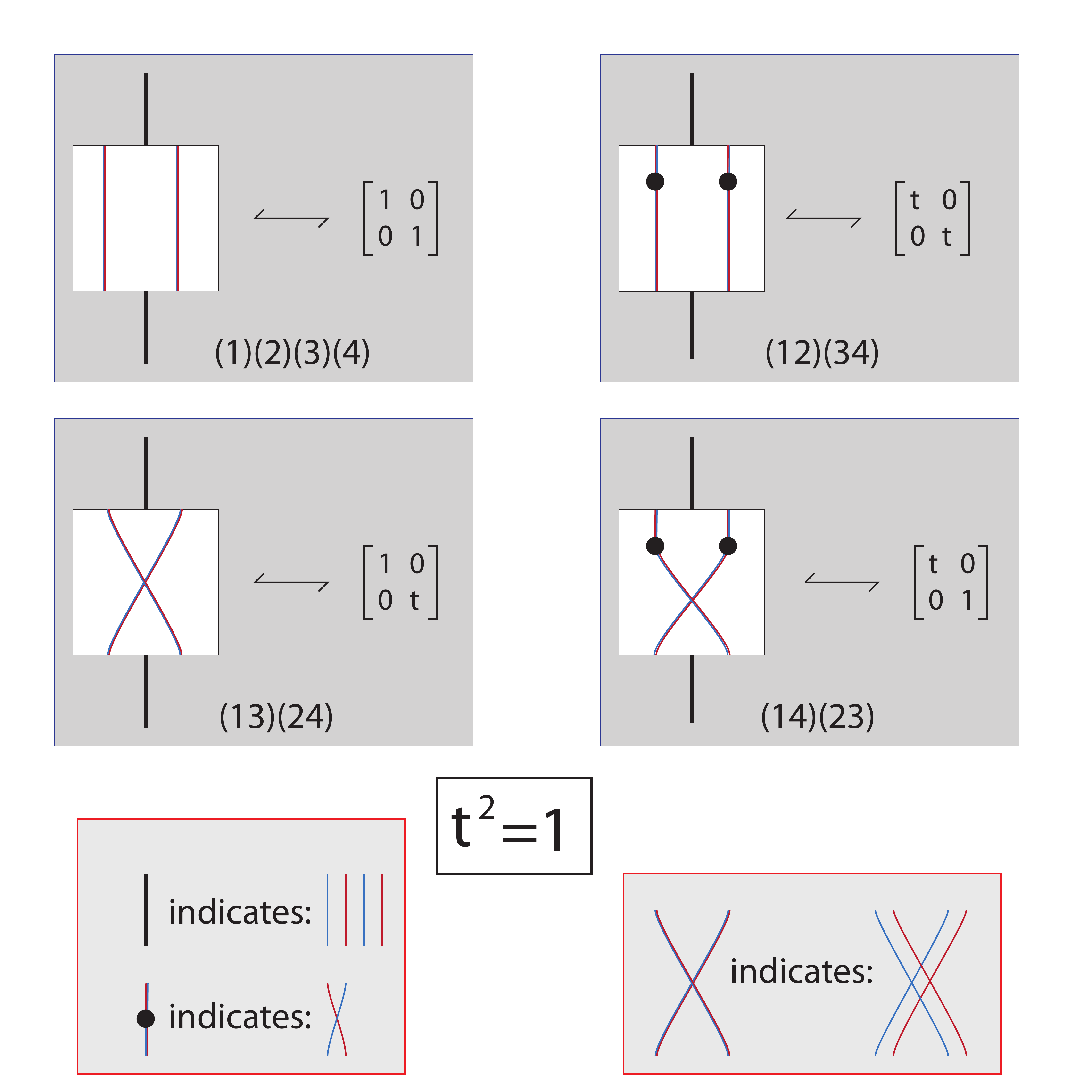}

Figure 2. A diagrammatic and a matrix representation of the Klein $4$ group
\end{center}

\newpage

\begin{sloppypar}
In Fig.~1, these permutations are illustrated as string diagrams. The illustration demonstrates that the permutation respects the grouping $\{1,2\}\{3,4\}$, and thus  the
associated permutation is imagined as a a pair of  ``strings-with-beads" The beads may slide through the crossings of the strings, and two beads upon a string cancel.  The domain of the permutations is written at the bottom of the diagram. So the $4$-cycle $(1324)$ indicates that the first point on the left bottom migrates up the string to the third position from the left at the top. 

\vspace{1cm} 
The matrices 

\vspace{1cm}

\[ \left [   \begin{array}{lr}  1 & 0  \\   0 & 1 \end{array} \right], \quad \left [   \begin{array}{lr}  -1 & 0  \\   0 & -1 \end{array} \right], \quad \left [   \begin{array}{lr}   0 & 1 \\ 1 & 0     \end{array} \right], \quad \left [   \begin{array}{lr}   0 & -1 \\ -1 & 0 \end{array} \right] \]  

\vspace{1cm}

\noindent
appear in the left column, and these constitute a subgroup that corresponds to the Klein $4$ group, $K_4$, which also is a normal subgroup of the alternating group $A_4$ on $4$-strings. The permutation representation consists of  $K_4= \{ (1), (12)(34), (13)(24),(14)(23) \}$.\end{sloppypar}

\vspace{1cm}
In Fig.~2, the string-with-beads illustration is distilled further, and an alternative $(2 \times 2)$-matrix representation is given in which off-diagonal entries are $0$, and the non-zero diagonal entries are elements of the cyclic group $\Z/2 =\langle t: t^2 =1 \rangle$. Since this is an index $2$ subgroup of the group of signed permutation matrices $H_2$, the larger group is written as a pair of cosets $H_2 = K \cup (-e_1,e_2)K$ where $e_j$ indicates the standard column vector basis element for $\R^2$, $j=1,2$.

\begin{sloppypar}
Here and below subgroups are ordered as sets and an ordering is induced upon the cosets.
So the subgroup $K$ is written as $\widehat{K}= \left( (1), (12)(34), (13)(24),(14)(23) \right)$, and its coset (in the permutation representation) is written as $\widehat{(12)K}= \left( (12), (34), (1324),(1423) \right).$ The actions of the elements of $H_2$ upon these pair of ordered sets is computed as follows:\end{sloppypar}
\newpage

\begin{eqnarray*}
(1)\widehat{K} &=&\left( (1), (12)(34), (13)(24),(14)(23) \right)  \\
(1)\widehat{(12)K}&=& \left( (12), (34), (1324),(1423) \right)\\
(12)(34)\widehat{K}= &=&  \left( (12)(34),(1), (14)(23), (13)(24) \right) \\
(12)(34)\widehat{(12)K}&=&\left( (34),(12),  (1423),(1324) \right) \\
(13)(24)\widehat{K} &=&\left((13)(24), (14)(23), (1), (12)(34)  \right)\\
(13)(24)\widehat{(12)K}&=&\left( (1423), (1324), (34), (12) \right)\\
(14)(23)\widehat{K} &=&\left((14)(23), (13)(24), (12)(34), (1) \right)\\
(14)(23)\widehat{(12)K}&=&\left( (1324), (1423), (12),(34) \right)\\
(12)\widehat{K}&=&  \left( (12), (34), (1324),(1423) \right) \\
(12)\widehat{(12)K}&=& \left( (1), (12)(34), (13)(24),(14)(23) \right)\\
(34)\widehat{K}&=& \left(  (34),(12),(1423), (1324) \right)  \\
(34)\widehat{(12)K}&=& \left( (12)(34), (1),(14)(23), (13)(24) \right)\\
(1324)\widehat{K}&=& \left( (1324),(1423),(12), (34) \right)  \\
(1324)\widehat{(12)K}&=& \left((14)(23),(13)(24),(12)(34), (1)\right)\\
(1423)\widehat{K}&=&   \left( (1423), (1324), (34),(12)\right) \\
(1423)\widehat{(12)K}&=&\left((13)(24), (14)(23),(1), (12)(34)\right) 
\end{eqnarray*}

Figure~3 indicates the diagrammatic representation that is induced by these actions. The remaining matrices have $1$ or $t$ on the off-diagonals and $0$s along the diagonals. These string representations are not any more convenient that the previous ones, but they easily generalize to Figs.~6 and~7 that indicate an interesting representation of the group $\Sigma_4$ of symmetries of the set $\{1,2,3,4\}$. In subsequent sections, the ``hat" notation will be dropped even though subgroups and cosets will be considered to be ordered.

\begin{center}
\includegraphics[width=0.65\paperwidth]{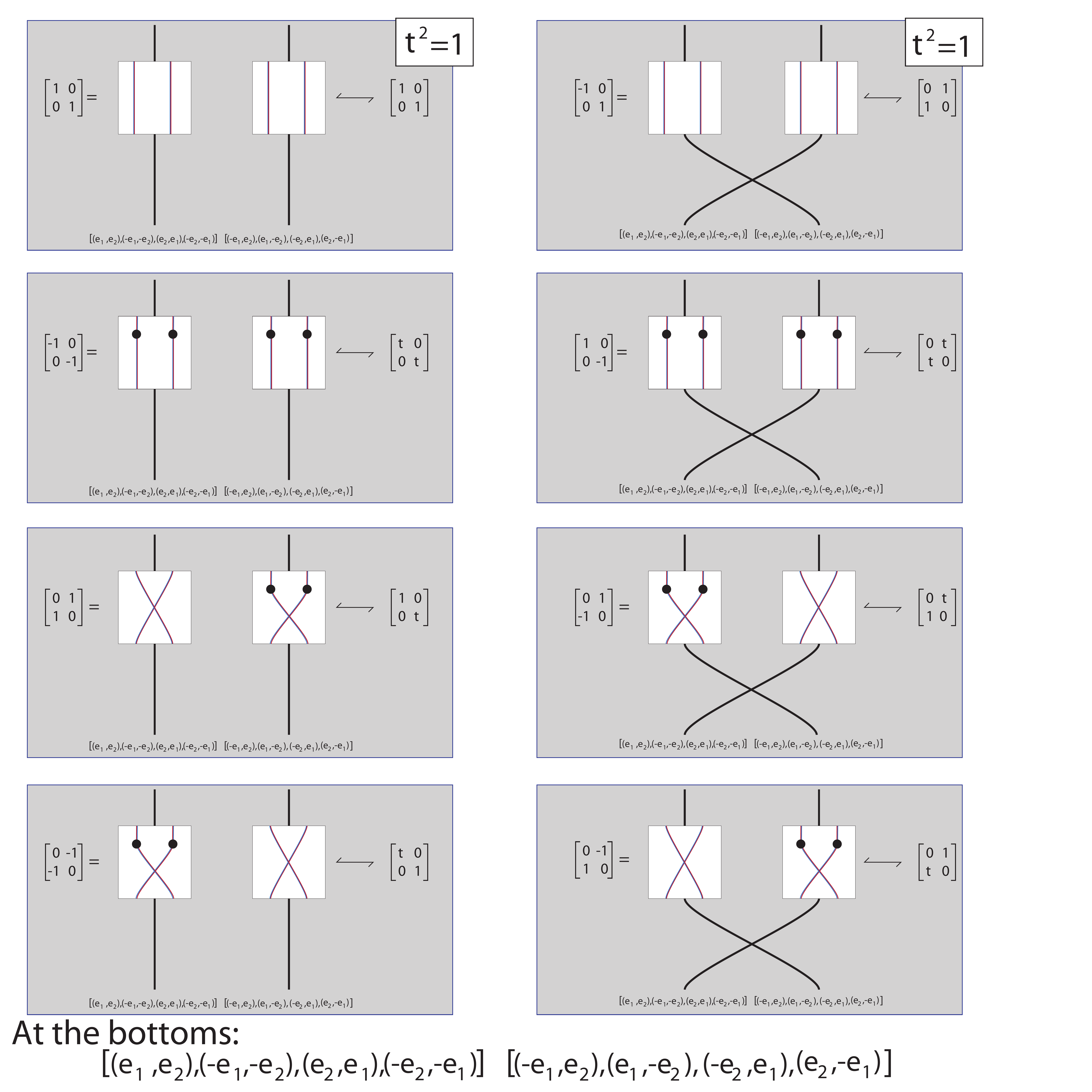}

Figure 3: An alternative representation of $H_2$
\end{center}

\begin{center}
\rule{3in}{0.005in}
\end{center}
The graphical expressions\[
\begin{tikzpicture}
\node[scale = 1] at (1,0) {$(1)$};
\node[scale = 1] at (3,0) {$(12)(34)$};
\node[scale = 1] at (5,0) {$(13)(24)$};
\node[scale = 1] at (7,0) {$(14)(23)$};
\node[scale = .5] at (1,1){
\begin{tikzpicture}
\draw[line width = 1] (1,1) rectangle (4,4);
\draw[line width = 1.5,purple] (1.5,1) to [in=-90, out = 90]  (1.5,4);
\draw[line width = 1.5, purple] (3.5,1) to [in=-90, out = 90]  (3.5,4);
\end{tikzpicture}};

\node[scale = .5] at (3,1){
\begin{tikzpicture}
\draw[line width = 1] (1,1) rectangle (4,4);
\draw[line width = 1.5,purple] (1.5,1) to [in=-90, out = 90]  (1.5,4);
\draw[line width = 1.5, purple] (3.5,1) to [in=-90, out = 90]  (3.5,4);
\draw[line width = 0, fill] (1.5,3) circle [radius=.1];
\draw[line width = 0, fill] (3.5,3) circle [radius=.1];
\end{tikzpicture}};

\node[scale = .5] at (5,1){
\begin{tikzpicture}
\draw[line width = 1] (1,1) rectangle (4,4);
\draw[line width = 1.5,purple] (1.5,1) to [in=-90, out = 90]  (3.5,4);
\draw[line width = 1.5, purple] (3.5,1) to [in=-90, out = 90] (1.5,4) ;
\end{tikzpicture}};

\node[scale = .5] at (7,1){
\begin{tikzpicture}
\draw[line width = 1] (1,1) rectangle (4,4);
\draw[line width = 1.5,purple] (1.5,1) to [in=-90, out = 90]  (3.5,3);
\draw[line width = 1.5, purple] (3.5,1) to [in=-90, out = 90] (1.5,3) ;
\draw[line width = 1.5,purple] (1.5,3) to [in=-90, out = 90]  (1.5,4);
\draw[line width = 1.5, purple] (3.5,3) to [in=-90, out = 90] (3.5,4) ;
\draw[line width = 0, fill] (1.5,3) circle [radius=.1];
\draw[line width = 0, fill] (3.5,3) circle [radius=.1];
\end{tikzpicture}};
\end{tikzpicture}\]
that indicate the Klein $4$ group, which is isomorphic to $\Z/2 \times \Z/2$, might make one ponder the nature of similar expressions that represent the elements of $\Z/4$. After all, both groups fit into the short exact sequence of groups
\[ 0 \rightarrow \Z/2 \rightarrow G \rightarrow \Z/2 \rightarrow 0 \] in which the middle group $G$ is either $\Z/4$ or $\Z/2 \times \Z/2$. By maintaining the convention that two beads upon a string cancel, i.e.
\[
\begin{tikzpicture}\node[scale=.5] {
\begin{tikzpicture}
\draw[line width = 1.5,purple] (1.5,1) to [in=-90, out = 90]  (1.5,4);
\draw[line width = 0, fill] (1.5,2) circle [radius=.1];
\draw[line width = 0, fill] (1.5,3) circle [radius=.1];
\node[scale = 1] at (2.5,2.5) {$=$};
\draw[line width = 1.5,purple] (3.5,1) to [in=-90, out = 90]  (3.5,4);
\end{tikzpicture}};\end{tikzpicture}\]
the elements in $\Z/4$ are represented by the glyphs indicated below.

\[
\begin{tikzpicture}
\node[scale =1] at (0,2.5) {The elements of $\Z/4$:};
\node[scale = 1] at (1,0) {$0$};
\node[scale = 1] at (3,0) {$1$};
\node[scale = 1] at (5,0) {$2$};
\node[scale = 1] at (7,0) {$3$};
\node[scale = .5] at (1,1)
{\begin{tikzpicture}
\draw[line width = 1] (1,1) rectangle (4,4);
\draw[line width = 1.5,purple] (1.5,1) to [in=-90, out = 90]  (1.5,4);
\draw[line width = 1.5, purple] (3.5,1) to [in=-90, out = 90] (3.5,4) ;
\end{tikzpicture}};

\node[scale = .5] at (3,1){
\begin{tikzpicture}
\draw[line width = 1] (1,1) rectangle (4,4);
\draw[line width = 1.5,purple] (1.5,1) to [in=-90, out = 90]  (3.5,3);
\draw[line width = 1.5, purple] (3.5,1) to [in=-90, out = 90] (1.5,3) ;
\draw[line width = 1.5,purple] (1.5,3) to [in=-90, out = 90]  (1.5,4);
\draw[line width = 1.5, purple] (3.5,3) to [in=-90, out = 90] (3.5,4) ;
\draw[line width = 0, fill] (3.5,3) circle [radius=.1];
\end{tikzpicture}};

\node[scale = .5] at (5,1){
\begin{tikzpicture}
\draw[line width = 1] (1,1) rectangle (4,4);
\draw[line width = 1.5,purple] (1.5,1) to [in=-90, out = 90]  (1.5,3);
\draw[line width = 1.5, purple] (3.5,1) to [in=-90, out = 90] (3.5,3) ;
\draw[line width = 1.5,purple] (1.5,3) to [in=-90, out = 90]  (1.5,4);
\draw[line width = 1.5, purple] (3.5,3) to [in=-90, out = 90] (3.5,4) ;
\draw[line width = 0, fill] (1.5,3) circle [radius=.1];
\draw[line width = 0, fill] (3.5,3) circle [radius=.1];
\end{tikzpicture}};

\node[scale = .5] at (7,1){
\begin{tikzpicture}
\draw[line width = 1] (1,1) rectangle (4,4);
\draw[line width = 1.5,purple] (1.5,1) to [in=-90, out = 90]  (3.5,3);
\draw[line width = 1.5, purple] (3.5,1) to [in=-90, out = 90] (1.5,3) ;
\draw[line width = 1.5,purple] (1.5,3) to [in=-90, out = 90]  (1.5,4);
\draw[line width = 1.5, purple] (3.5,3) to [in=-90, out = 90] (3.5,4) ;
\draw[line width = 0, fill] (1.5,3) circle [radius=.1];
\end{tikzpicture}};
\end{tikzpicture}\]
Moreover, these glyphs can be found among the original glyphs depicted in Fig.~1  that represent the matrices $\pm (e_1,e_2)$  and $\pm(-e_2,e_1).$

It is an amusing exercise to generalize the graphical representations of $\Z/4$ and $\Z/2 \times \Z/2$ towards representations of $\Z/(2n)$ and $\Z/n \times Z/n$. 

In this section, we generalize the diagrammatic representations that have been given towards other groups.

\begin{center}
\includegraphics[width=0.5\paperwidth]{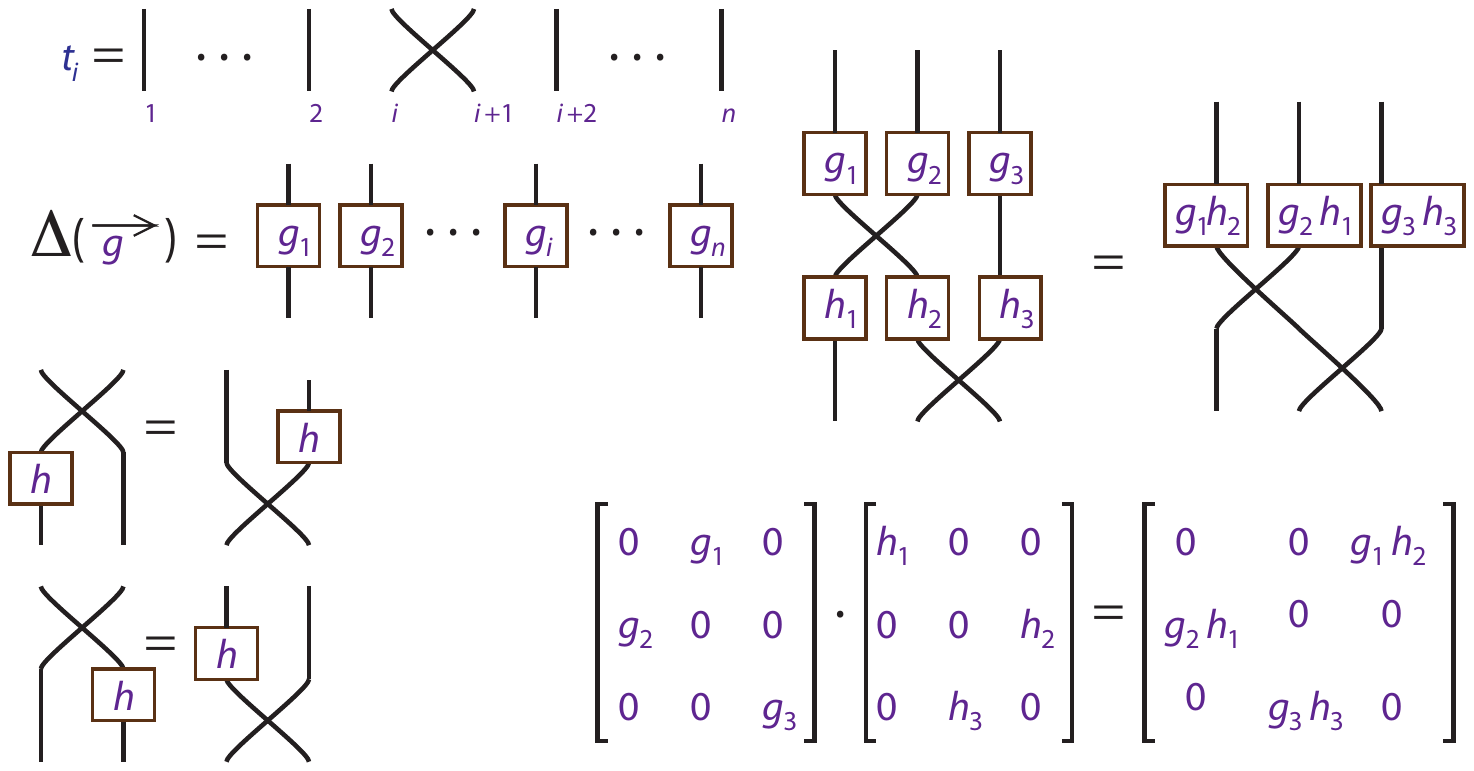}

Figure 4. Permutation representation in semi-direct products
\end{center}

\subsection{Generalities}\label{Generalities}
Let $G$  and $H$ denote  finite groups where both are given as subgroups of permutation groups:
Say $G\subset \Sigma_m$ and $H \subset \Sigma_n$. For example, if $h\in H$, then the map $f_h: k \mapsto hk$
represents $H$ as a subgroup of the set of permutations of the underlying set of $H$. Consider the semi-direct product $W=G^n \rtimes H$ in which $H$ acts on the factors of $G^n$ by permuting the coordinates. The group $W$ fits into a split short
exact sequence 
$$ 1 \rightarrow G^n \stackrel{\ i \ }{\longrightarrow} W  \stackrel{\ p \ }{\longrightarrow} H \rightarrow 1.$$
For convenience of computation, we express  the permutation group $\Sigma_n$ of which $H$ is a subgroup as
the set of permutation matrices. Let $$e_j =(0,\ldots, 0, \underbrace{1}_{j\/{\mbox{\rm th}}}, 0 \ldots, 0)^t$$ denote the $j\/$th standard unit vector in $\R^n$. Then a permutation $\sigma \in \Sigma_n$ is given as the matrix
$(e_{\sigma 1}, e_{\sigma 2}, \ldots, e_{\sigma n})$. In our experience, there is always some miscommunication in
this formula. So for example, the $3$-cycle $(1,2,3) \in \Sigma_3$ corresponds to the matrix
$$(1,2,3) \leftrightarrow \left[ \begin{array}{ccc} 0 & 0 & 1 \\ 1 & 0 & 0 \\ 0 & 1 & 0 \end{array} \right]. $$
Let $g_1, g_2, \ldots, g_n \in G$ denote elements in $G$. Let $1$ denote the identity element of $G$. 
Then an element of $ (\vec{g}, \sigma) \in W$ is expressed as the matrix $(g_1 e_{\sigma 1},g_2 e_{\sigma 2},
\ldots, g_n e_{\sigma n})$. In this way, the entry $g_i$ appears in the $i\/$th row. For example, 
$$ \left[ \begin{array}{ccc} g_1 & 0 & 0  \\ 0 & g_2 & 0  \\  0 & 0 & g_3 \end{array} \right] \cdot \left[ \begin{array}{ccc} 0 & 0 & 1 \\ 1 & 0 & 0 \\ 0 & 1 & 0 \end{array} \right] = \left[ \begin{array}{ccc} 0 & 0 &g_1 \\ g_2 & 0 & 0 \\ 0 & g_3 & 0 \end{array} \right].$$
Thus if $\Delta(\vec{g})$ denotes the  matrix in which $g_1$ through $g_n$ appear along the diagonal and $0\/$s
appear elsewhere, then the element $(\vec{g}, \sigma)\in W$ is the matrix product $\Delta(\vec{g})\cdot
(e_{\sigma 1}, e_{\sigma 2}, \ldots, e_{\sigma n}).$ When two such matrix representatives are multiplied, only one 
non-zero entry appears in any row and column as a product of pairs of elements in $G$. In this way, computations 
in the {\it 
semi-direct   product of $G$ and $H$} are computed as formal matrix products. 
In general this process can be translated into string diagrams as follows. Recall that the symmetric group has the
presentation
$$\Sigma_n= \left\langle t_1, \ldots, t_{n-1}:  \begin{array}{c} t_it_j=t_jt_i \ {\mbox{\rm for}} \ |i-j|>1, \\ 
t_it_jt_i=t_jt_it_j \ {\mbox{\rm for}} \ |i-j|=1,  \\ {\mbox{\rm and}}  \  t_i^2=1 \  {\mbox{\rm for}} \ i=1,\ldots, n-1 
\end{array} \right\rangle.$$ The generator $t_i$ corresponds to the transposition $t_i = (i, i+1)$. We can draw this
as a string diagram as
in Fig.~4. In this illustration the element $\Delta(\vec{g})$ is also depicted, and a sample 
product is shown as a string product and in matrix form. 
Elements in the semi-direct product are drawn as beads labeled by elements in $G$ on the top of string 
permutation diagrams.

Please note that if $G$ also has a permutation representation (as will often be the case here), then the
permutation representation of the semi-direct product is obtained by bundling strings together and implementing 
the representatives of the elements $g_1$, $\ldots$,  $g_n$ at the top of the diagrams. This string replication is exactly 
what happened in our first examples. We apply similar ideas in the next subsection.

\subsection{The permutation groups}
\label{SS:perm}

\begin{center}
\includegraphics[width=0.65\paperwidth]{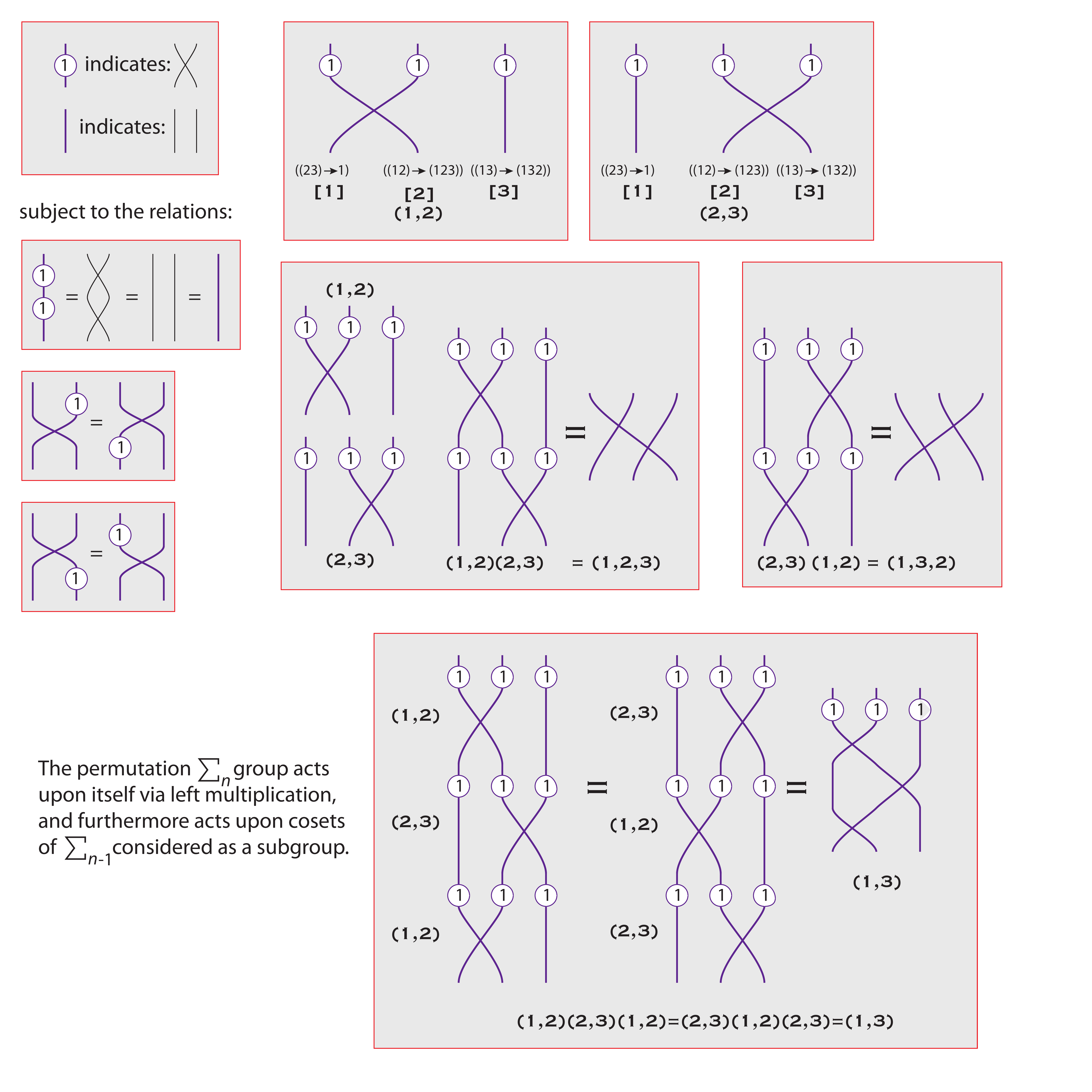}

Figure 5. The group $\Sigma_3$ acting upon itself and its ribbon presentation
\end{center}


The symmetric group, $\Sigma_n$ on $n$ elements has order $n!$. Its presentation in terms of the {\it standard
generators,} $t_i$, is given above. The group acts upon itself by left multiplication. This action represents the group
as a set of permutations on an excessively large set. 

Rather than attempting to gain control of these actions, we present pictures that describe $\Sigma_3$ and $\Sigma_4$.

We start from the even easier case $n=2$.  An oriented edge is considered as an 
arrow $(t_1 \rightarrow 1)$. The action of $t_1$ upon this arrow is  $t_1(t_1 \rightarrow 1)= (1 \leftarrow t_1)$. By 
convention, arrows point from odd permutations towards even permutations. 

In $\Sigma_3$, consider the
subgroup $\{t_2,1\}$, and label and orient the cosets as follows $[1] =  (t_2 \rightarrow 1), [2]=(t_1 \rightarrow 
t_1t_2),$ and $[3] = (t_2t_1t_2 \rightarrow t_2t_1)$.  The idea is to consider the first element to be the subgroup, act
 by the transposition $t_1=(1,2)$, label the result $[2]$, and subsequently act upon this by $t_2=(2,3)$. Each transposition reverses an arrow. Fig.~5 indicates, the resulting permutation representation. It has a slight advantage over the standard string picture because the existence or non-existence of the three beads at the top of the strings redundantly indicates the parity of the permutation represented; a permutation in $\Sigma_3$ is odd if and only if there are three beads at the tops of the strings.  The notation is redundant because the parity is also obtained by counting the number of crossings between pairs of strings.

The result is a faithful representation of $\Sigma_3 \rightarrow (\Z/2)^3 \rtimes \Sigma_3$. The target  is  a
subgroup of the $(3\times 3)$ dimensional signed permutation matrices with determinant $1$.  Specifically, $(1,2) 
\leftrightarrow (-e_2,-e_1, -e_3)$ while $(2,3) \leftrightarrow (-e_1,-e_3,-
e_2)$. Other elements are computed as matrix products. We can envision arrows$ ((-1,1,1) \stackrel{ \ [1] \ }
{\longrightarrow} (1,-1,-1)),$ $((1,-1,1) \stackrel{ \ 
[2] \ }{\longrightarrow} (-1,1,-1)),$  and $((1,1,-1) \stackrel{ \ [3] \ }{\longrightarrow} (-1,-1,1))$ between antipodal points in 
the cube $[-1,1]^3$ that pass through the origin. The signed permutation matrices act upon this configuration, and 
the beaded depiction of  Fig.~5 encodes this   action.

\begin{center}
\rule{3in}{0.005in}
\end{center}

Consider the sequence of  normal subgroups
\[ 0 \lt \Z/2 \lt K_4 \lt A_4 \lt \Sigma_4 \]
with quotient groups 
\[ \Sigma_4/ A_4 \cong \Z/2, \quad A_4/ K_4 \cong \Z/3, \quad K_4/(\Z/2) \cong \Z/2. \]
The {\it alternating group $A_4$} consists of the permutations that can be written as a product of an even number of transpositions taken from the generating set $\{t_1=(12),t_2=(23),t_3=(34)\}$. The depiction 
\[
\begin{tikzpicture}
\node[scale = 1] at (1,0) {$(1)$};
\node[scale = 1] at (3,0) {$(12)(34)$};
\node[scale = 1] at (5,0) {$(13)(24)$};
\node[scale = 1] at (7,0) {$(14)(23)$};
\node[scale = .5] at (1,1){
\begin{tikzpicture}
\draw[line width = 1] (1,1) rectangle (4,4);
\draw[line width = 1.5,purple] (1.5,1) to [in=-90, out = 90]  (1.5,4);
\draw[line width = 1.5, purple] (3.5,1) to [in=-90, out = 90]  (3.5,4);
\end{tikzpicture}};

\node[scale = .5] at (3,1){
\begin{tikzpicture}
\draw[line width = 1] (1,1) rectangle (4,4);
\draw[line width = 1.5,purple] (1.5,1) to [in=-90, out = 90]  (1.5,4);
\draw[line width = 1.5, purple] (3.5,1) to [in=-90, out = 90]  (3.5,4);
\draw[line width = 0, fill] (1.5,3) circle [radius=.1];
\draw[line width = 0, fill] (3.5,3) circle [radius=.1];
\end{tikzpicture}};

\node[scale = .5] at (5,1){
\begin{tikzpicture}
\draw[line width = 1] (1,1) rectangle (4,4);
\draw[line width = 1.5,purple] (1.5,1) to [in=-90, out = 90]  (3.5,4);
\draw[line width = 1.5, purple] (3.5,1) to [in=-90, out = 90] (1.5,4) ;
\end{tikzpicture}};

\node[scale = .5] at (7,1){
\begin{tikzpicture}
\draw[line width = 1] (1,1) rectangle (4,4);
\draw[line width = 1.5,purple] (1.5,1) to [in=-90, out = 90]  (3.5,3);
\draw[line width = 1.5, purple] (3.5,1) to [in=-90, out = 90] (1.5,3) ;
\draw[line width = 1.5,purple] (1.5,3) to [in=-90, out = 90]  (1.5,4);
\draw[line width = 1.5, purple] (3.5,3) to [in=-90, out = 90] (3.5,4) ;
\draw[line width = 0, fill] (1.5,3) circle [radius=.1];
\draw[line width = 0, fill] (3.5,3) circle [radius=.1];
\end{tikzpicture}};
\end{tikzpicture}\] of the elements in $K_4$
is dependent upon considering the ordering  $T=((1),(12)(34))$ of the subgroup $\Z/2$ and the ordered coset  $(13)(24)T=((13)(24),(14)(23))$ and then subsequently considering the actions of the elements in $K_4$ upon these ordered cosets. The ordered cosets of the ordered subgroup
\[K_4 = ( (1), (12)(34), (13)(24),(14)(23)) \] are listed below (both as cosets in the alternating group $A_4$ and as ordered cosets in the ordered coset $(12)A_4$. In this listing,  a convenient naming convention is presented that uses typographically economical representatives.

\[\begin{array}{lcr}
\left.  \begin{array}{lcl} [1] &=& [(1),(12)(34),(13)(24),(14)(23)]
\\   

[123] &=& [(123),(134),(243),(142)]  
\\ 
  
[132] &=&  [(132),(234),(124),(143)] \end{array} \right\} & =& A_4;  \\ && \\
 \left. \begin{array}{lcl} [12] &=& [(12),(34),(1324),(1423)] 
 \\
 
[13] &=& [(13),(1234),(24),(1432)]  
\\

[14] &=&  [(23),(1342),(1243),(14)] \end{array} \right\} & =& (12)A_4.  \end{array} \]
A systematic, tedious and thorough calculation results in depicting the $24$ elements in the symmetric group $\Sigma_4$ as indicated in the Figs.~6 and~7. 

These depictions may only serve as curiosities. However, each diagrammatic depiction of a given group element implicitly contains the corresponding column of the multiplication table for that element. Moreover, the composition series
\[ 0 \lt \Z/2 \lt K_4 \lt A_4 \lt \Sigma_4 \]
with quotient groups 
\[ \Sigma_4/ A_4 \cong \Z/2, \quad A_4/ K_4 \cong \Z/3, \quad K_4/(\Z/2) \cong \Z/2\]
is visually apparent. 

A possible alternative notation for the diagrams that are illustrated in Figs.~6 and~7 is proposed here. It depends upon the hierarchy of the quotient groups, and an observation about the six elements of $K_4$ that appear in each diagram. Relabel the elements of $K_4$ as follows.  

\[
\begin{tikzpicture}
\node[scale = 1] at (1,0) {$[0]$};
\node[scale = 1] at (3,0) {$[1]$};
\node[scale = 1] at (5,0) {$[2]$};
\node[scale = 1] at (7,0) {$[3]$};
\node[scale = .5] at (1,1){
\begin{tikzpicture}
\draw[line width = 1] (1,1) rectangle (4,4);
\draw[line width = 1.5,purple] (1.5,1) to [in=-90, out = 90]  (1.5,4);
\draw[line width = 1.5, purple] (3.5,1) to [in=-90, out = 90]  (3.5,4);
\end{tikzpicture}};

\node[scale = .5] at (3,1){
\begin{tikzpicture}
\draw[line width = 1] (1,1) rectangle (4,4);
\draw[line width = 1.5,purple] (1.5,1) to [in=-90, out = 90]  (1.5,4);
\draw[line width = 1.5, purple] (3.5,1) to [in=-90, out = 90]  (3.5,4);
\draw[line width = 0, fill] (1.5,3) circle [radius=.1];
\draw[line width = 0, fill] (3.5,3) circle [radius=.1];
\end{tikzpicture}};

\node[scale = .5] at (5,1){
\begin{tikzpicture}
\draw[line width = 1] (1,1) rectangle (4,4);
\draw[line width = 1.5,purple] (1.5,1) to [in=-90, out = 90]  (3.5,4);
\draw[line width = 1.5, purple] (3.5,1) to [in=-90, out = 90] (1.5,4) ;
\end{tikzpicture}};

\node[scale = .5] at (7,1){
\begin{tikzpicture}
\draw[line width = 1] (1,1) rectangle (4,4);
\draw[line width = 1.5,purple] (1.5,1) to [in=-90, out = 90]  (3.5,3);
\draw[line width = 1.5, purple] (3.5,1) to [in=-90, out = 90] (1.5,3) ;
\draw[line width = 1.5,purple] (1.5,3) to [in=-90, out = 90]  (1.5,4);
\draw[line width = 1.5, purple] (3.5,3) to [in=-90, out = 90] (3.5,4) ;
\draw[line width = 0, fill] (1.5,3) circle [radius=.1];
\draw[line width = 0, fill] (3.5,3) circle [radius=.1];
\end{tikzpicture}};
\end{tikzpicture}\]
Then observe that within the cosets $A_4$ and $(12)A_4$ the three boxed quantities are ordered and of the form $[000]$, $[123]$, $[231]$, and $[312]$.  Moreover, the two pairs in a given element are always of the form $[000][000]$, $[123][123]$, $[231][312]$, or $[312][231]$ reading this triple from the point of view of the strings at the bottom from left to right.{\footnote{We regret the ambiguity in this sentence. The accompanying figures should eliminate the ambiguity.}} 
\newpage

\begin{center}
\includegraphics[width=0.6\paperwidth]{AltCosets.pdf}

Figure 6. An alternative graphical depiction of the elements in $A_4$
\end{center}

\begin{center}
\includegraphics[width=0.6\paperwidth]{NonAltCosets.pdf}

Figure 7.An alternative graphical depiction of the elements in $(12)A_4$
\end{center}

We write:
\[
\begin{array}{|lccc|} \hline
K_4: & (1) & \leftrightharpoons & \left( {\sf I}, 0, [000] \right) \\
& (12)(34) & \leftrightharpoons & \left( {\sf I}, 0, [123] \right) \\
& (13)(24) & \leftrightharpoons & \left( {\sf I}, 0, [231] \right) \\
& (13)(24) & \leftrightharpoons & \left( {\sf I}, 0, [312] \right) \\ \hline
\end{array} \quad 
\begin{array} {|lccc|} \hline
(12)K_4: & (12) & \leftrightharpoons & \left( {\sf X}, 0, [000] \right) \\
&(34) & \leftrightharpoons & \left( {\sf X}, 0, [123] \right) \\
&(1324) & \leftrightharpoons & \left( {\sf X}, 0, [231] \right) \\
& (1432) & \leftrightharpoons & \left( {\sf X}, 0, [312] \right) \\ \hline
\end{array}\]
\[
\begin{array} {|lccc|} \hline
(123)K_4: & (123) & \leftrightharpoons & \left( {\sf I}, 1, [000] \right) \\
&(134) & \leftrightharpoons & \left( {\sf I}, 1, [312] \right) \\
&(243) & \leftrightharpoons & \left( {\sf I}, 1, [123] \right) \\
&(142) & \leftrightharpoons & \left( {\sf I}, 1,[231]  \right]) \\ \hline
\end{array} \quad
\begin{array} {|lccc|} \hline
(23)K_4:& (23) & \leftrightharpoons & \left( {\sf X}, 1, [000] \right) \\
&(1342) & \leftrightharpoons & \left( {\sf X}, 1, [312]  \right) \\
&(1243) & \leftrightharpoons & \left( {\sf X}, 1, [123] \right) \\
&(14) & \leftrightharpoons & \left( {\sf X}, 1, [231]\right]) \\ \hline
\end{array} \]
\[
\begin{array} {|lccc|} \hline
(132)K_4:&(132) & \leftrightharpoons & \left( {\sf I}, 2, [000] \right) \\
&(234) & \leftrightharpoons & \left( {\sf I}, 2,[231]  \right) \\
&(124) & \leftrightharpoons & \left( {\sf I}, 2, [312] \right) \\
&(143) & \leftrightharpoons & \left( {\sf I}, 2,[123]  \right]) \\ \hline
\end{array} \quad
\begin{array} {|lccc|} \hline
(13)K_4:&(13) & \leftrightharpoons & \left( {\sf X}, 2, [000] \right) \\
&(1234) & \leftrightharpoons & \left( {\sf X}, 2,[231]  \right) \\
&(24) & \leftrightharpoons & \left( {\sf X}, 2, [312] \right) \\
&(1432) & \leftrightharpoons & \left( {\sf X}, 2,[123]  \right]) \\ \hline
\end{array}\]

The first character {\sf I} or {\sf X} indicates whether the element is in $A_4$ or in $(12)A_4$, respectively. The second character $0,1,$ or $2$ indicates which of the cosets in $\Z/3 \cong A_4/K_4$ the element represents. And the triple $[000]$, $[123]$, $[231]$, $[312]$ indicates the triple of elements in $K_4$ on the left bottom of each element as indicated in the figures. It is possible, but perhaps unhelpful, to develop arithmetic rules for combining these elements to reflect the multiplication in $\Sigma_4$. The standard $4$-string picture is probably more intuitive and less prone to algorithmic errors. 

\section{Proof of Theorem~\ref{two}}
\label{S:Pr2}
 The  proof  will be a consequence of studying the permutation action of $G$ upon itself by left multiplication: for $g'\in G$, let  $g:g' \mapsto gg'$.
Denote the elements of $H$ by $h_1,h_2, \ldots, h_k$. In fact, let us observe that this is an arbitrary ordering of $H$. We write $(H)=(h_1,h_2,\ldots, h_k)$ to indicate $H$ as an ordered set. Choose cosets $a_1H$, $\ldots,$ $a_nH$, and observe that this, too, indicates an ordering so that as an ordered set $(G)=(a_1H, \ldots, a_{n}H)$. The action of $G$ upon itself leads to a permutation action of $g$ upon the set of cosets.
So for $g\in G$, there is a permutation $\eta_g=\eta \in \Sigma_n$ such that $g(a_1H, \ldots, a_{n}H)=(a_{\eta1}H,a_{\eta2}H, \ldots, a_{\eta n}H).$ Let us write $ga_iH=a_jH$ to indicate that $\eta i=j.$ Furthermore, there is a permutation $\sigma_{i,j}=\sigma \in \Sigma_k$ such that  $(ga_i h_1, \ldots, ga_i h_k)= (a_j h_{\sigma1}, \ldots, a_jh_{\sigma k})$. Thus far, we have established a group embedding
$G\subset (\Sigma_k)^n \rtimes \Sigma_n$. To complete the proof, we establish that the permutation representation of $g$ acting upon any coset corresponds to the action of $H$ upon itself.

Since $ga_iH=a_jH$, we have that there is an $h_{i,j}\in H$ with $h_{i,j}=a_j^{-1}ga_i$. 
So $h_{i,j}h_\ell= a_j^{-1} g a_i h_\ell= a_j^{-1}a_j h_{\sigma \ell}= h_{\sigma_{i,j} \ell}$. 
In other words, the permutation action between $ga_iH$ and $a_jH$ coincides with multiplication by some element $h_{i,j}\in H$. Therefore, $G\subset (H)^n \rtimes \Sigma_n$ as desired. This completes the proof.

\section{The quaternions}
\label{S:Quat}

Two diagrammatic representations of the quaternions are developed herein. One or the other will be used in the subsequent descriptions of the dicyclic groups and the binary polyhedral groups. Let $Q_8$ denote the group with eight elements $Q_8=\{\pm 1, \pm {\boldsymbol i}, \pm {\boldsymbol j}, \pm {\boldsymbol k}\}$ where the multiplication is defined as:

$${\boldsymbol i}\cdot {\boldsymbol j}={\boldsymbol k }=-{\boldsymbol j}\cdot {\boldsymbol i},$$ $${\boldsymbol j}\cdot {\boldsymbol k}={\boldsymbol i} =-{\boldsymbol k}\cdot {\boldsymbol j},$$ 
$${\boldsymbol k}\cdot {\boldsymbol i}={\boldsymbol j} =-{\boldsymbol i}\cdot {\boldsymbol k},$$ and 
$${\boldsymbol i}\cdot {\boldsymbol i}={\boldsymbol j}\cdot {\boldsymbol j}={\boldsymbol k}\cdot {\boldsymbol k}=-1.$$

All of our own calculations among ${\boldsymbol i}, {\boldsymbol j},$ and ${\boldsymbol k}$ are facilitated by means of using the standard diagram: 
\[
\begin{tikzpicture}
\node[scale=.6] at (0,0) {\begin{tikzpicture}
\draw[line width = 1.5] (3,3) circle [radius=3];
\draw[line width = 0, fill, white] (3,6) circle [radius=.3];
\draw[line width = 0, fill, white] (5.56,1.5) circle [radius=.3];
\draw[line width = 0, fill, white] (6-5.56,1.5) circle [radius=.3];
\node[scale =1.5] at (3,6) {${\boldsymbol i}$};
\node[scale =1.5] at (6-5.56,1.5) {${\boldsymbol j}$};
\node[scale =1.5] at (5.56,1.5) {${\boldsymbol k}$};
\draw[line width = 1.5] (0,2)--(.25,1.8)--(.3,2.1);
\draw[line width = 1.5] (5.1,1.05)--(5.39,1.2)--(5.4,.9);
\draw[line width = 1.5] (3.5,5.8)--(3.25,6)--(3.5,6.1);
\end{tikzpicture}};
\end{tikzpicture}\]
Moreover, at several junctures, the following table (and variations upon it that depend upon signs) was helpful.
\begin{center}
\begin{tabular}{||c||c|c|c|c||}
\hline \hline
$\cdot$ & $1$ & ${\boldsymbol i}$ & ${\boldsymbol j}$ &  ${\boldsymbol k}$ \\ \hline \hline 
$1$ & $1$ & ${\boldsymbol i}$ & ${\boldsymbol j}$ &  ${\boldsymbol k}$ \\ \hline
${\boldsymbol i}$ & ${\boldsymbol i}$ & $-1$ & ${\boldsymbol k}$ &  $-{\boldsymbol j}$ \\ \hline
${\boldsymbol j}$ & ${\boldsymbol j}$ &  $-{\boldsymbol k}$ & $-1$ &  ${\boldsymbol i}$ \\ \hline
${\boldsymbol k}$ & ${\boldsymbol k}$ &  ${\boldsymbol j}$  &  $-{\boldsymbol i}$ & $-1$ \\ \hline \hline
\end{tabular}\end{center}

\begin{center}
\includegraphics[width=0.6\paperwidth]{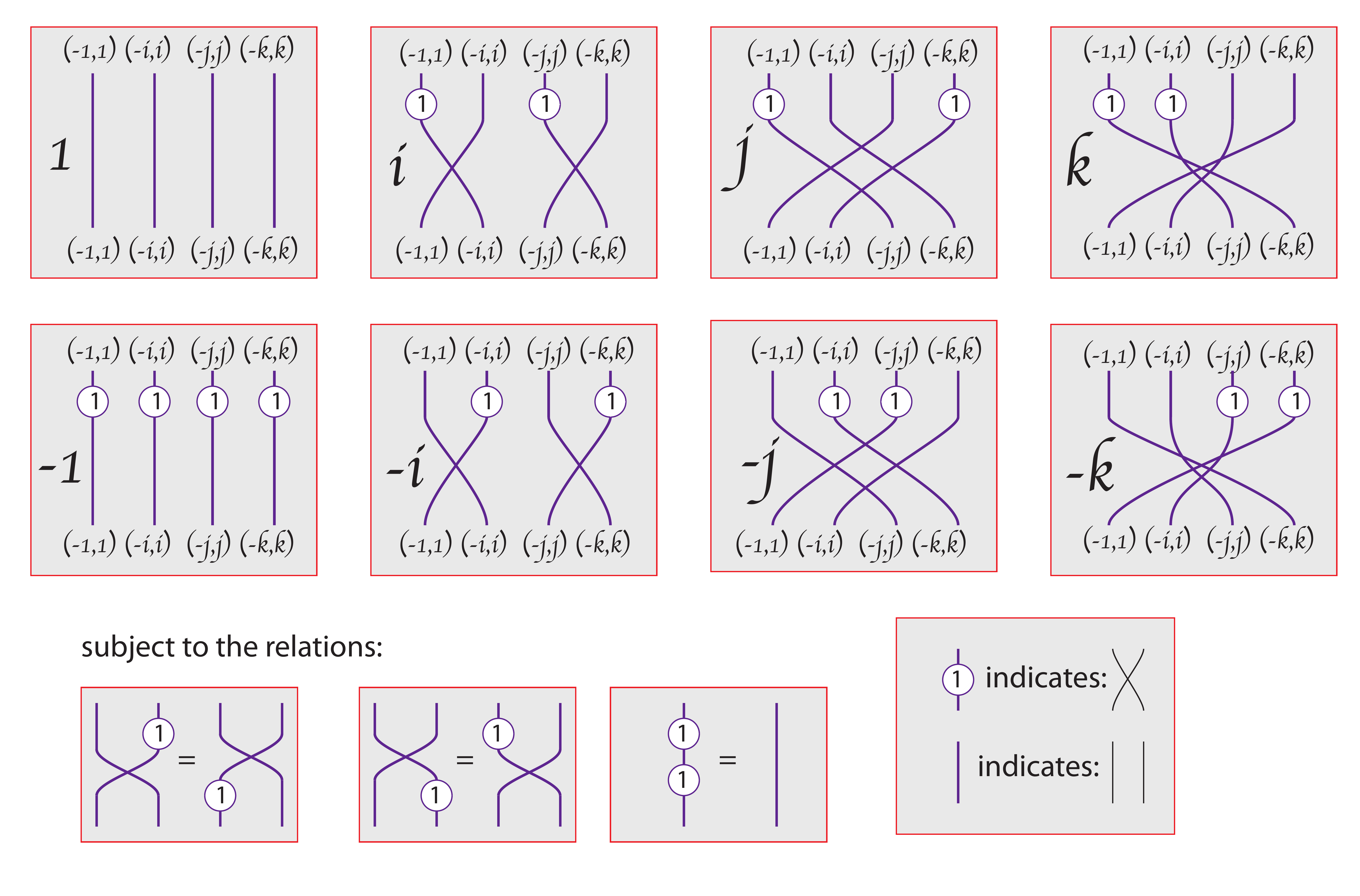}

Figure 8. A ribbon presentation of the group $Q_8$

\end{center}

Permutation representations of a group depend upon the ordering of the group and the choice of ordered subgroup upon which the group acts. Consider the short exact sequence
$$0 \rightarrow \Z/2 \rightarrow Q_8 \stackrel{ \ p \ }{\rightarrow} \Z/2 \times \Z/2 \rightarrow 0$$ 
for which the inclusion identifies $\Z/2$ with $\{\pm 1\}$ and the quotient group is identified with the Klein $4$-group.  The subgroup $\{ -1, 1\}$ is replaced by an ordered set $(-1, 1)$, and the cosets are also ordered as
$((-1, 1),$ $(-{\boldsymbol i} , {\boldsymbol i}),$ $(- {\boldsymbol j}, {\boldsymbol j}),$ 
$(-{\boldsymbol k} , {\boldsymbol k})).$ 
We have
\begin{eqnarray*}
1((-1, 1),(-{\boldsymbol i} , {\boldsymbol i}),(- {\boldsymbol j}, {\boldsymbol j}),(-{\boldsymbol k} , {\boldsymbol k})) & =& 
((1, 1),
(-{\boldsymbol i} , {\boldsymbol i}),
(- {\boldsymbol j}, {\boldsymbol j}),
(-{\boldsymbol k} , {\boldsymbol k})); \\
{\boldsymbol i}((-1, 1),(-{\boldsymbol i} , {\boldsymbol i}),(- {\boldsymbol j}, {\boldsymbol j}),(-{\boldsymbol k} , {\boldsymbol k})) & =& 
((-{\boldsymbol i}, {\boldsymbol i}),
(1 , -1),
( -{\boldsymbol k}, {\boldsymbol k}),
({\boldsymbol j} , -{\boldsymbol j})); \\
{\boldsymbol j}((-1, 1),(-{\boldsymbol i} , {\boldsymbol i}),(- {\boldsymbol j}, {\boldsymbol j}),(-{\boldsymbol k} , {\boldsymbol k})) & =& 
((-{\boldsymbol j}, {\boldsymbol j}),
({\boldsymbol k} , -{\boldsymbol k}),
(1, -1),
(-{\boldsymbol i} , {\boldsymbol i})); \\
{\boldsymbol k}((-1, 1),(-{\boldsymbol i} , {\boldsymbol i}),(- {\boldsymbol j}, {\boldsymbol j}),(-{\boldsymbol k} , {\boldsymbol k})) & =& 
((-{\boldsymbol k}, {\boldsymbol k}),
(-{\boldsymbol j} , {\boldsymbol j})
( {\boldsymbol i}, -{\boldsymbol i}),
(1 , -1)).
\end{eqnarray*}

Fig.~8 indicates a beaded string depiction of the quaternions.  The thick strings indicate a pair of parallel strings, and a bead indicates a crossing between them.  Beads are allowed to pass through crossings, and a pair of beads on a single string cancel. Under this identification, the group of unit quaternions, $Q_8$, is a subgroup of $(\Z/2)^4 \rtimes K_4,$   where $K_4$ indicates the Klein $4$ group. Specifically, the projection $p$ maps elements as follows: $p(\pm 1)=(1),$ $p(\pm {\boldsymbol i})=(12)(34),$  $p(\pm {\boldsymbol j})=(13)(24),$  and $p(\pm {\boldsymbol k})= (14)(23)$. 

This depiction does not work well with the previous depiction of the Klein $4$ group. In that representation two thinner strings were conglomerated into one, and here each of these thinner strings represents two strings that are thinner still. So the beads in the quaternions are viewed as ``half beads" in the two string view of $K_4$. A calculus could be developed in that context which would involve different types of string markers, but it is much less intuitive that the one in Fig~8. 

It is perhaps instructive to ponder the difference between the dihedral group  $H_2$ and the unit quaternions $Q_8$. These are the two non-abelian groups of order $8$. The former has $K_4$ as a subgroup, and the latter has $K_4$ as a quotient. Their distinct diagrammatic depictions reflect their distinct nature. 

This completes the proof of item (1)  of Theorem~\ref{main}.

\begin{center}
\rule{3in}{0.005in}
\end{center}

For another representation, consider the short exact sequence
$$0 \rightarrow \Z/4 \rightarrow Q_8 \stackrel{ \ p \ }{\rightarrow} \Z/2  \rightarrow 0$$
for which the kernel is the subgroup  $(-1,-{\boldsymbol i},1, {\boldsymbol i})$ --- which is ordered as indicated. More precisely, we are considering these in a cyclic order since they represent cardinal points on the great circle in the $3$-sphere $$S^3=\{w+x {\boldsymbol i} + y {\boldsymbol j} + z {\boldsymbol k}: w^2+x^2+y^2+ z^2 =1 \}.$$  See also Fig.~17. Then there are two cosets $((-1,-{\boldsymbol i},1, {\boldsymbol i}),$ $(-{\boldsymbol j},{\boldsymbol k},{\boldsymbol j}, -{\boldsymbol k})).$ The action of $Q_8$ on the ordered subgroup is computed as follows. 
Multiplication by $(-1)$ puts a half-twist in both cosets. Meanwhile,
\begin{eqnarray*}
{\boldsymbol i}((-1,-{\boldsymbol i},1, {\boldsymbol i}),(-{\boldsymbol j},{\boldsymbol k},{\boldsymbol j}, -{\boldsymbol k}))
& = & 
((-{\boldsymbol i},1,{\boldsymbol i}, -1),(-{\boldsymbol k}, -{\boldsymbol j},{\boldsymbol k},{\boldsymbol j})); \\
{\boldsymbol j}((-1,-{\boldsymbol i},1, {\boldsymbol i}),(-{\boldsymbol j},{\boldsymbol k},{\boldsymbol j}, -{\boldsymbol k}))
& = &
((-{\boldsymbol j},{\boldsymbol k}, {\boldsymbol j},-{\boldsymbol k}),(1,-{\boldsymbol i},-1,-{\boldsymbol i})); \\
{\boldsymbol k}((-1,-{\boldsymbol i},1, {\boldsymbol i}),(-{\boldsymbol j},{\boldsymbol k},{\boldsymbol j}, -{\boldsymbol k}))
&=&
((-{\boldsymbol j},{\boldsymbol k},{\boldsymbol j},-{\boldsymbol k}),(1,{\boldsymbol i},-1, -{\boldsymbol i})).\end{eqnarray*}

\begin{center}
\includegraphics[width=0.6\paperwidth]{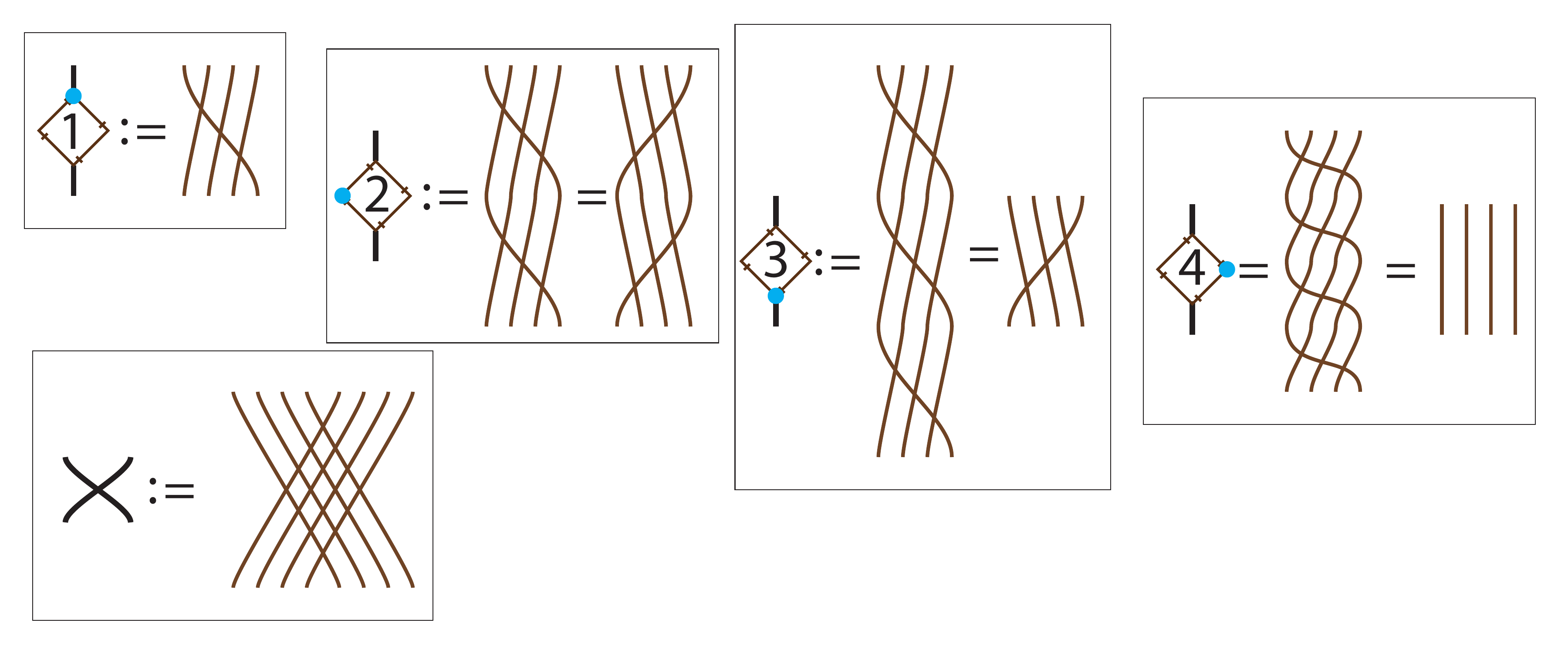}

Figure 10. Bundling strings into groups of four
\end{center}

\begin{center}
\includegraphics[width=0.6\paperwidth]{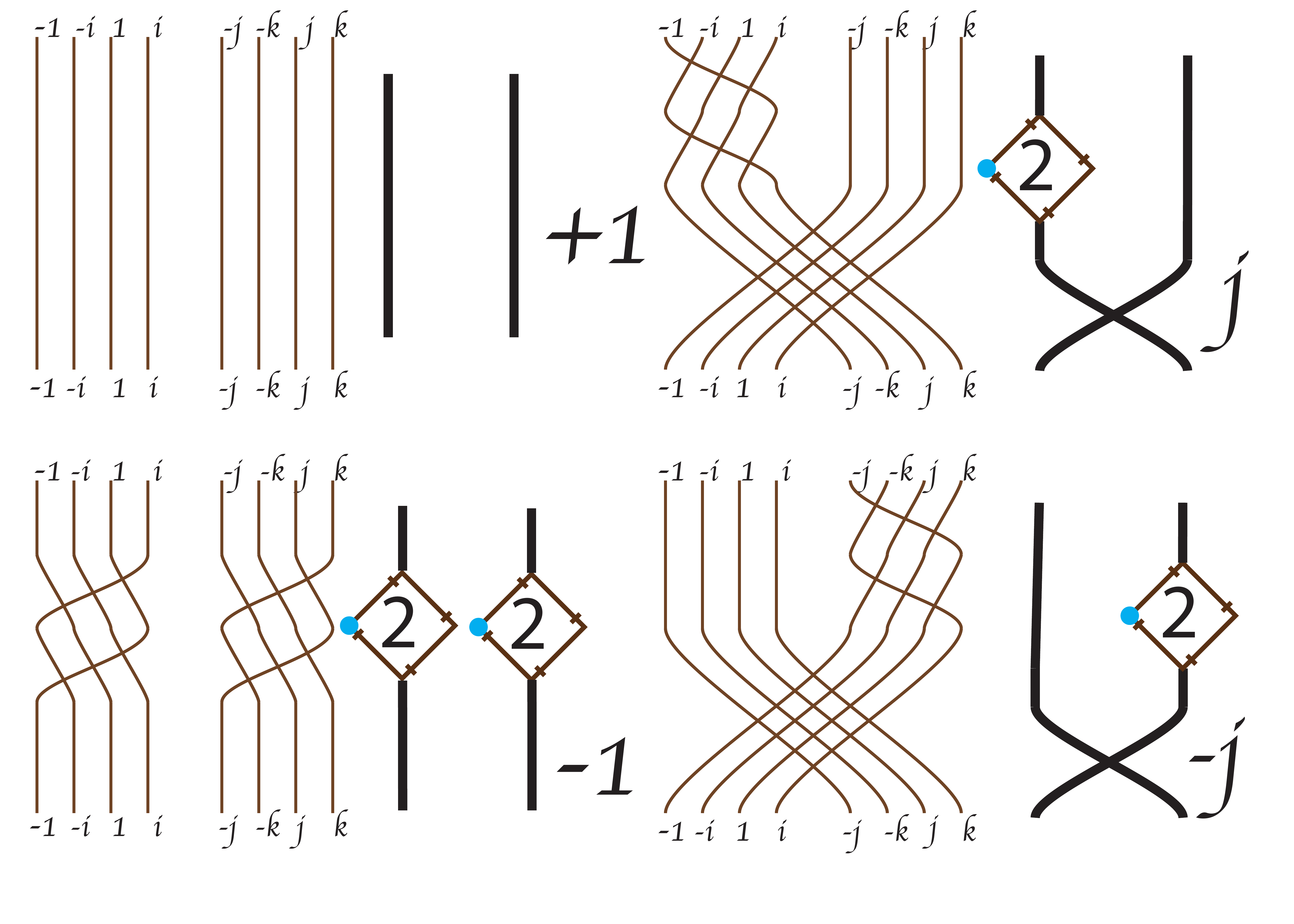}

Figure 11. The string representation of the  quaternionic group $Q_8$, part 1
\end{center}

\begin{center}
\includegraphics[width=0.6\paperwidth]{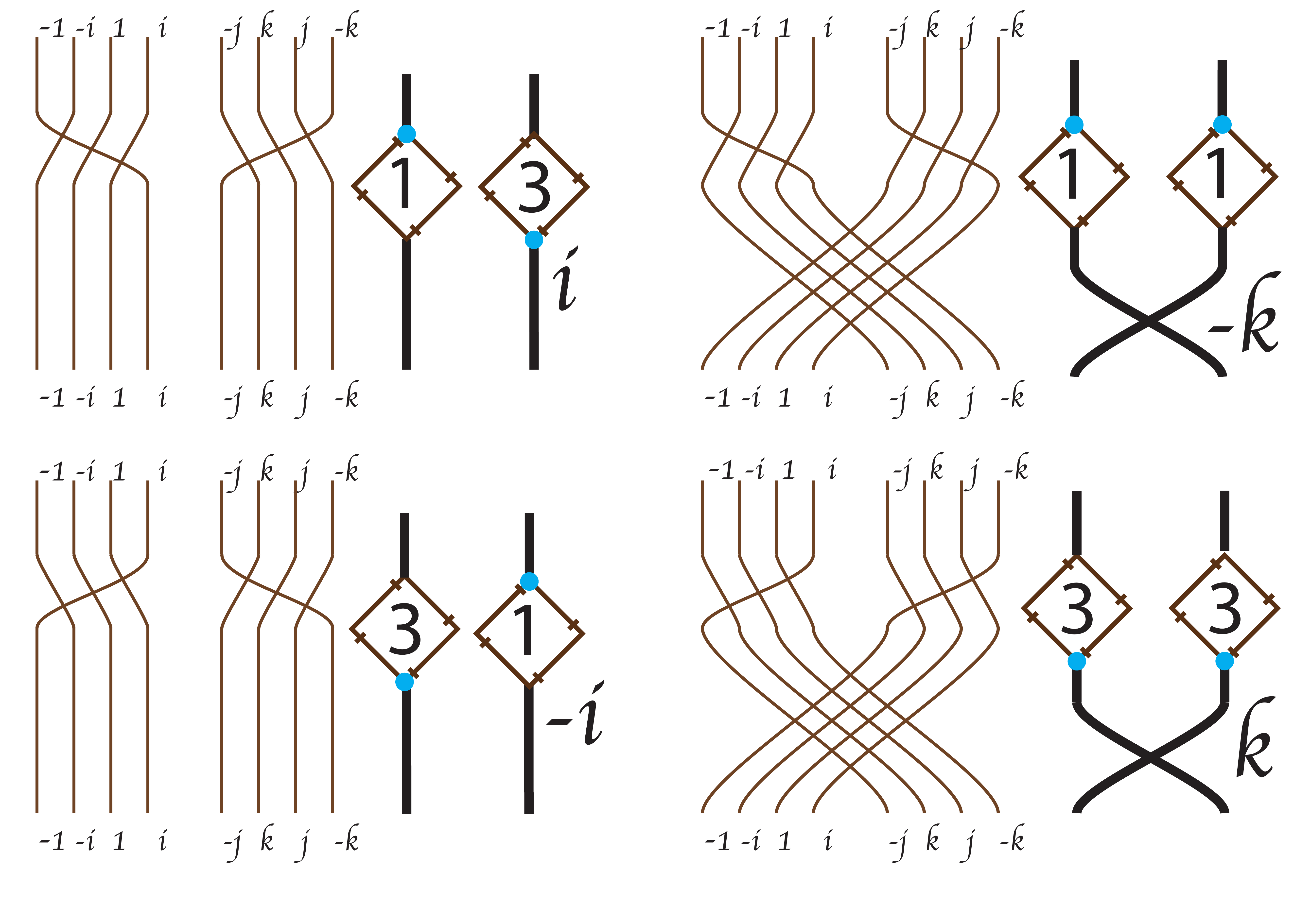}

Figure 12. The string representation of the  quaternionic group $Q_8$, part 2
\end{center}

\begin{center}
\includegraphics[width=0.6\paperwidth]{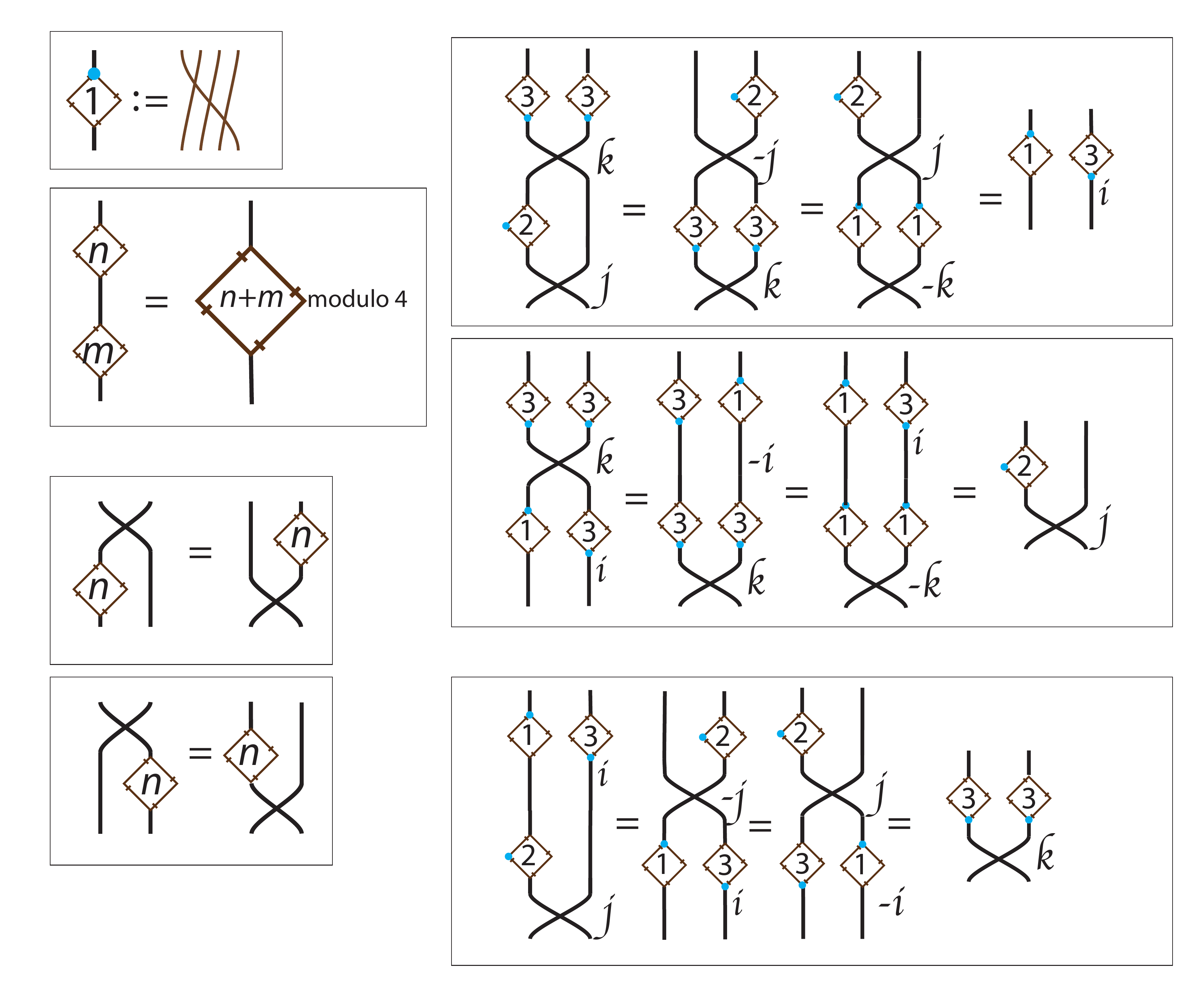}

Figure 13. Some products in $Q_8$
\end{center}

Fig.~9 indicates the method for bundling the eight strings into two sets of four. The four string bundles acquire $0$ through $3$ positive  $(1/2)$-twists, and  a $(3/2)$-twist is equivalent to  a $(-1/2)$-twist. Fig.~10 and ~11 indicate the representations of elements $\pm 1,$ $\pm{\boldsymbol i},$ $\pm{\boldsymbol j},$ and $\pm{\boldsymbol k}.$ There is a faithful permutation representation  $Q_8 \rightarrow (\Z/4)^2 \rtimes \Sigma_2$ into the semi-direct product. The multiplication rules are indicated in Fig.~13.
Some mathematicians will point out that group $Q_8$ is not a semi-direct product, but  there are faithful representations into semi-direct products. This completes the representation for item (2) in Theorem~\ref{main}.

\begin{center}
\rule{3in}{0.005in}
\end{center}

\subsection{Matrix representation.}\label{MatRep}
In Section~\ref{Generalities}, we observed that elements in a semi-direct product with a permutation group could be represented as ``permutation matrices with group entries." In the case of the quaternions, the group elements are elements in a cyclic group of order $4$ which may be represented as being generated by any of ${\boldsymbol i}$, ${\boldsymbol j}$, or ${\boldsymbol k}.$ So the string representation in Figs.~10 and~11, corresponds to the representation via Pauli matrices
\[{\boldsymbol i} \leftrightharpoons\left[ \begin{array}{cc} {\boldsymbol i} & 0 \\ 0 &  -{\boldsymbol i}\end{array} \right], \quad
{\boldsymbol j} \leftrightharpoons   \left[ \begin{array}{cc} 0 & -1 \\ 1 & 0 \end{array} \right], \quad {\boldsymbol k} \leftrightharpoons \left[ \begin{array}{cc} 0& -{\boldsymbol i}\\    -{\boldsymbol i} & 0 \end{array} \right]\] as trace $0$ matrices. Each corresponds to a $90^\circ$ rotation of a complex plane. By dividing the entries by $2$, one can obtain the standard basis of 
the lie algebra ${\it su}(2)$, and the Lie bracket $[A,B]=[AB-BA]$ induces the cross product multiplication that is, perhaps, familiar from vector calculus. 

A direct computation shows that
\[ w \left[ \begin{array}{cc} 1 & 0 \\ 0 &  1\end{array} \right] 
+ x\left[ \begin{array}{cc} {\boldsymbol i} & 0 \\ 0 &  -{\boldsymbol i}\end{array} \right] 
+ y \left[ \begin{array}{cc} 0 & -1 \\ 1 & 0 \end{array} \right] 
+ z \left[ \begin{array}{cc} 0& -{\boldsymbol i}\\    -{\boldsymbol i} & 0 \end{array} \right]\] \[
= \left[ \begin{array}{cc} w+x{\boldsymbol i}& -y-z{\boldsymbol i}\\    y-z{\boldsymbol i}  &  w-x{\boldsymbol i} \end{array} \right] \] is a unitary matrix   when $w^2+x^2+y^2+z^2=1$. And this expression 
rewrites an element in the $3$-sphere $S^3$  as an element in the matrix group $\SU(2)$.

\section{The dicyclic groups}
\label{S:Dic}

The second representation of the group $Q_8$ is part and parcel of the next diagrammatic descriptions of the dicyclic groups. The dicyclic group  of order $4n$ is given by the presentation
$${\mbox{\rm Dic}}_{n}= \langle \rho,x: \rho^{2n}=1, \ x^2=\rho^n, \  \rho x= x\rho^{-1} \rangle.$$
It maps $2$-to-$1$ to the dihedral group of order $2n$. Here the inclusion 
${\mbox{\rm Dic}}_{n} \subset S^3$  is given by $\rho\mapsto \cos{(\pi/n)} + \sin{(\pi/n)} {\boldsymbol i}$ and $x \mapsto {\boldsymbol j}$. In particular,  in this representation ${\mbox{\rm Dic}}_2 \approx Q_8$ with $\rho \mapsto {\boldsymbol i},$ but because the diagrams for ${\boldsymbol i}$ and ${\boldsymbol j}$ are interchanged, the representations of $\pm{\boldsymbol k}$ will also be switched. The background material here has been compiled from a number of wikipedia sources, specifically \cite{wiki:SO(3),wiki:Dic}.

It is a good time to introduce a formula for the projection from $\SU(2)=S^3$ to the group,  ${\SO}(3)$, of $(3\times 3)$ orthogonal matrices that are of determinant $1$.
 The continuous $2$-to-$1$ covering $ S^3 \stackrel{\ p \ }{\longrightarrow} \SO(3) $ is given by the formula
$$p(w+x {\boldsymbol i} + y {\boldsymbol j} + z {\boldsymbol k}) = 
\left[ \begin{array}{ccc} 1 - 2(y^2+ z^2) & 2(xy-zw) & 2(xz+yw) \\
2(xy+zw) & 1-2(x^2+z^2) & 2(yz-xw) \\
2(xz-yw) & 2(yz+xw) & 1-2(x^2+y^2) \end{array} \right].$$ It is a straightforward, yet tedious calculation to show that the columns of the matrix form an orthogonal basis for $\R^3$. It is easy to see that antipodal points 
have the same image: $p(-w-x {\boldsymbol i} - y {\boldsymbol j} - z {\boldsymbol k})=p(w+x {\boldsymbol i} + y {\boldsymbol j} + z {\boldsymbol k})$ because each entry is a homogeneous degree $2$ polynomial.	

\begin{center}
\includegraphics[width=0.6\paperwidth]{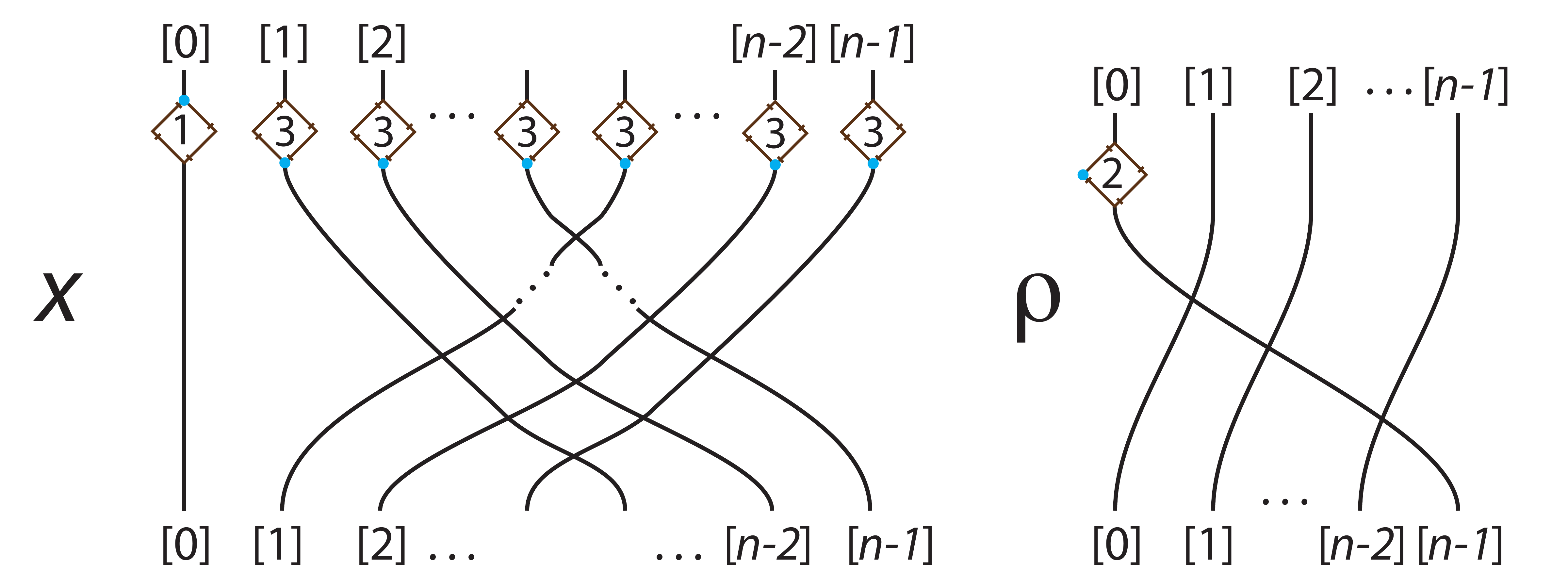}

Figure 14. Generators of the dicyclic groups in permutation form
\end{center}

\begin{center}
\includegraphics[width=0.6\paperwidth]{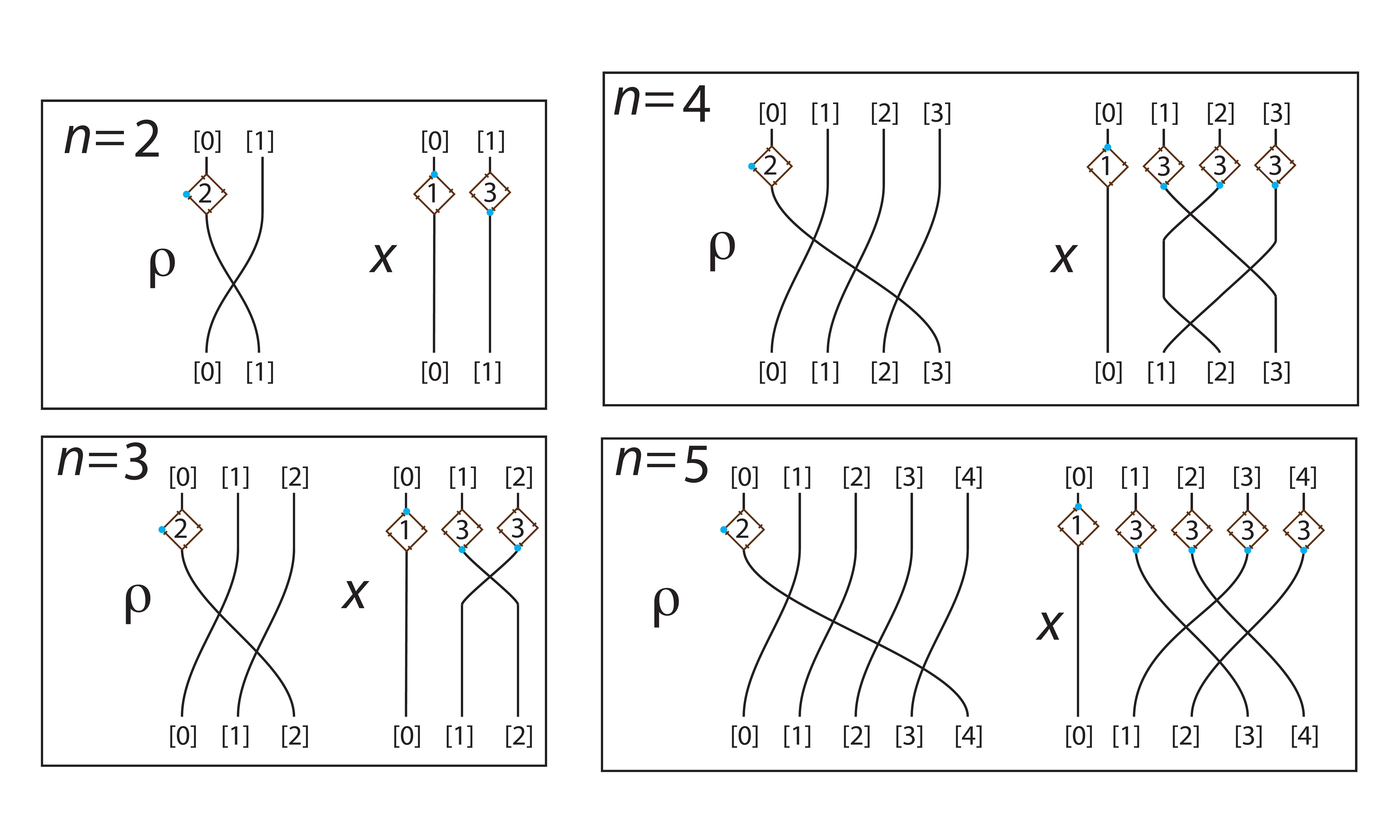}

Figure 15. Generators of the dicyclic groups  for $n=2,3,4,5$
\end{center}

Let ${\boldsymbol u}_n(\ell)	= \cos{\left( \ell\pi/n\right)}	+  \sin {\left( \ell\pi/n\right)} {\boldsymbol i}$, and let
${\boldsymbol v}_n(\ell)	= \cos{\left( \ell\pi/n\right)}{\boldsymbol j}	+  \sin {\left( \ell\pi/n\right)} {\boldsymbol k}$. A direct computation shows  that
$$p({\boldsymbol u}_n(\ell)) =\left[ \begin{array}{ccc} 1 & 0 &0 \\ 0 & \cos{\left( 2\ell\pi/n\right)} & -\sin{\left( 2\ell\pi/n\right)} \\
0 & \sin{\left( 2\ell\pi/n\right)} & \cos{\left( 2\ell\pi/n\right)}\end{array} \right],$$ 
and $$
p({\boldsymbol v}_n(\ell)) =\left[ \begin{array}{ccc} -1 & 0 &0 \\ 0 & \cos{\left( 2\ell\pi/n\right)} & \sin{\left( 2\ell\pi/n\right)} \\
0 & \sin{\left( 2\ell\pi/n\right)} & -\cos{\left( 2\ell\pi/n\right)}\end{array} \right].$$
Thus   a power of the generator $\rho^\ell \mapsto {\boldsymbol u}_n(\ell)$  projects to a rotation in the $(y,z)$ plane in $\R^3$, and $p(\pm {\boldsymbol j})$ is a reflection. Furthermore, $\rho^\ell x \mapsto {\boldsymbol v}_n(-\ell) \in S^3$. 

In a manner analogous to that of the quaternions, the orbits of the subgroup $J=\{-1,-{\boldsymbol j}, 1, {\boldsymbol j}\}$ which is ordered as $(-1,-{\boldsymbol j}, 1, {\boldsymbol j})$ are considered. For a more concise notation,  write the ordered cosets as $[\ell]= \rho^\ell ( \rho^n, \rho^nx, 1, x)= (\rho^{\ell+n}, \rho^{\ell+n} x, \rho^\ell, \rho^\ell x),$ for $\ell =0, \ldots, n-1.$ Note that
$[0]=  (-1, -x, 1, x).$
We compute $$x[\ell]=(x\rho^{\ell+n}, x\rho^{\ell+n} x, x\rho^\ell, x\rho^\ell x)=(\rho^{n-\ell}x,-\rho^{n-\ell }, \rho^{-\ell}x, - \rho^{-\ell}),$$
since $\rho^{2n}=1$, 
and $$\rho[n-1]=(\rho^{n-1+n+1}, \rho^{n-1+n+1} x, \rho^{n}, \rho^n x)=(1,  x, -1, - x).$$
In particular, $x[0]=(-x,1, x, - 1).$ In Fig.~14, the actions of $x$ and $\rho$ are depicted using the $4$-string convention from the presentation of the quaternions that were given above. The $n=2,3,4$ and $5$ cases are shown in Fig.~15. Note that $x$ corresponds to ${\boldsymbol j}$ and $\rho$ corresponds to ${\boldsymbol i}$, as expected. In Fig.~16, we demonstrate the relations in the group for $n=3.$

This completes the embedding of the dicyclic group in Theorem~\ref{main}, item (3).

\begin{center}
\includegraphics[width=0.57\paperwidth]{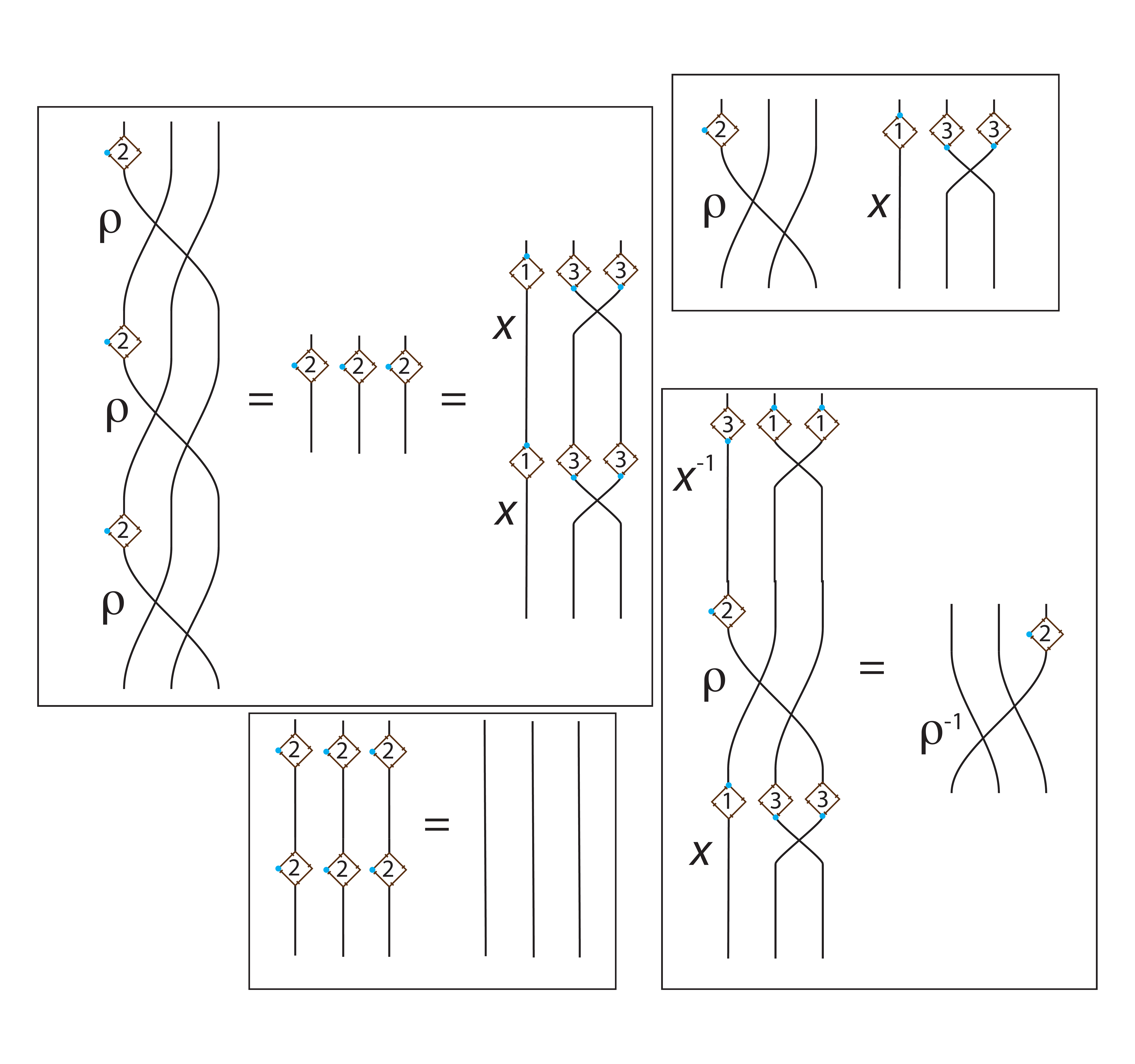}

Figure 16. Geometric manifestations of the relations in ${\mbox{\rm Dic}}_3$
\end{center}

\subsection{Two string representation and the corresponding matrices.}
Consider the cyclic group $\langle \zeta \rangle$ where $\zeta = {\boldsymbol u}_n(1)= \cos{\left( \pi/n\right)}	+  \sin {\left( \pi/n\right)} {\boldsymbol i}$.  Notice that 
${\boldsymbol j} \zeta^\ell = \cos{\left( \ell\pi/n\right)}{\boldsymbol j}	- \sin {\left( \ell\pi/n\right)} {\boldsymbol k}= {\boldsymbol v}_n(-\ell)$. So ${\mbox{\rm Dic}}_{n}$ is the union of the two ordered cosets 
\begin{eqnarray*} {\mbox{\rm Dic}}_{n} & =& ( \zeta, \zeta^2, \ldots , \zeta^{2n-1}, 1 ) \cup {\boldsymbol j} ( \zeta, \zeta^2, \ldots , \zeta^{2n-1}, 1 )\\ & =&  ( {\boldsymbol u}_n(1), \ldots, {\boldsymbol u}_n(2n-1), 1 ) \cup ( {\boldsymbol v}_n(-1), \ldots, {\boldsymbol v}_n(1-2n), {\boldsymbol j}).\end{eqnarray*} Multiplication by $\zeta$ rotates $\langle \zeta \rangle$ through one angle in a counterclockwise direction, and rotates $\boldsymbol j \langle \zeta \rangle$ clockwise. Multiplication by ${\boldsymbol j}$ interchanges the two cosets, but causes a $180^\circ$ rotation when the right coset moves back to the left. 
These two cosets lie upon a Hopf link in the $3$-sphere $S^3$. In the cases, $n=3$ and $n=4$, the situation is illustrated in Fig.~17. These two dicyclic groups are isomorphic to subgroups of the binary octahedral group $\widetilde{\Sigma_3}$. It may be helpful to try to envision these groups as they sit inside $S^3$. 

\begin{center}
\includegraphics[width=0.4\paperwidth]{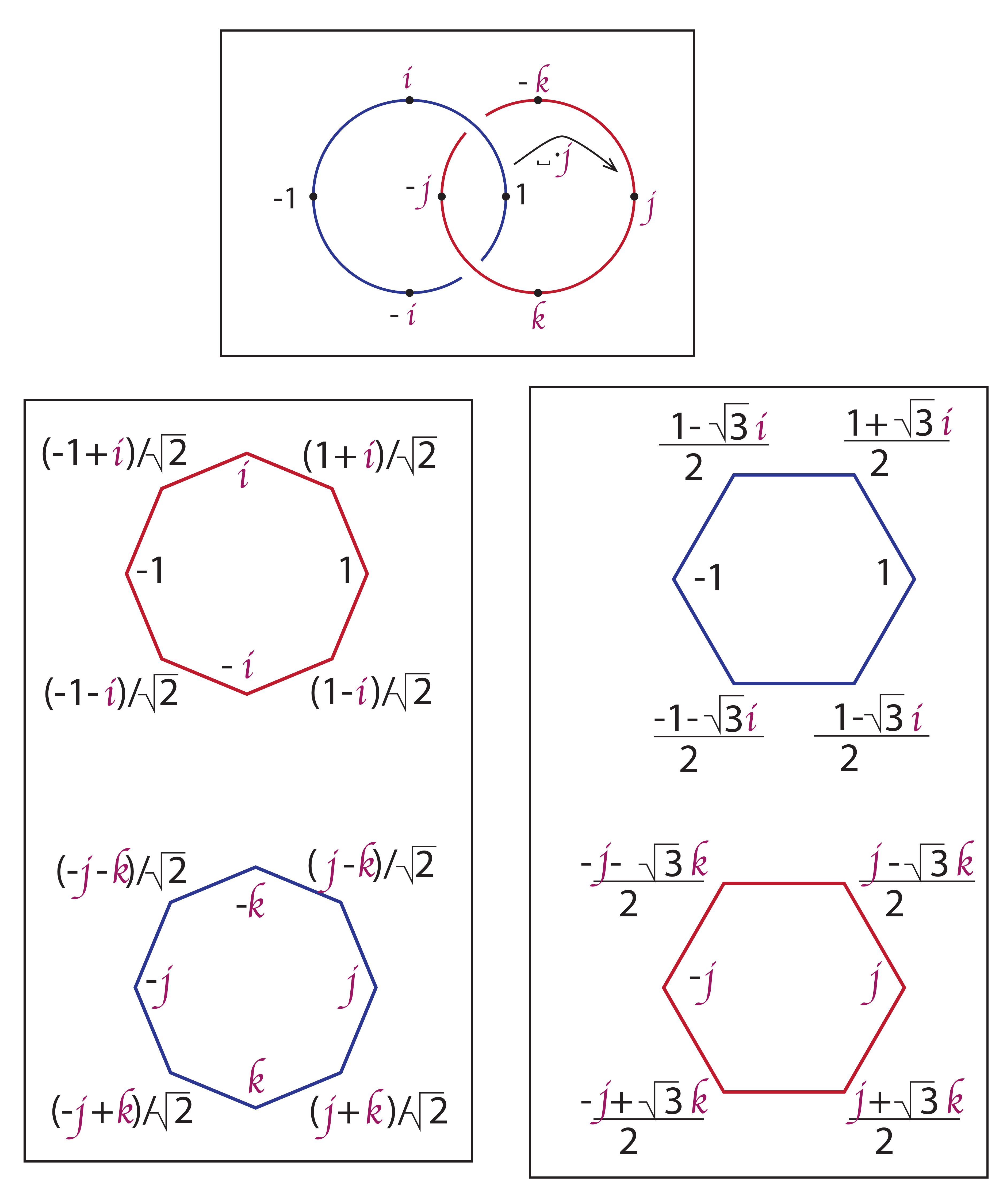}

Figure 17. The dicyclic groups as they lie upon the $(1,{\boldsymbol i}) \cup ({\boldsymbol j},{\boldsymbol k})$ Hopf link.
\end{center}

The two-string representation is depicted in Fig.~18. In the case of ${\mbox{\rm Dic}}_4$, the two string representatives are depicted encircling the  corresponding octagons in Fig.~19.

Since $\zeta = e^{(\pi {\boldsymbol i})/n}$, there is a corresponding matrix representation
\[{\rho} \leftrightharpoons \left[ \begin{array}{cc} \zeta & 0 \\ 0 & \zeta^{-1}
\end{array} \right], \quad
{\boldsymbol j} \leftrightharpoons \left[ \begin{array}{cc} 0& -1 \\ 1 & 0 \end{array} \right].\]

\begin{center}
\includegraphics[width=0.3\paperwidth]{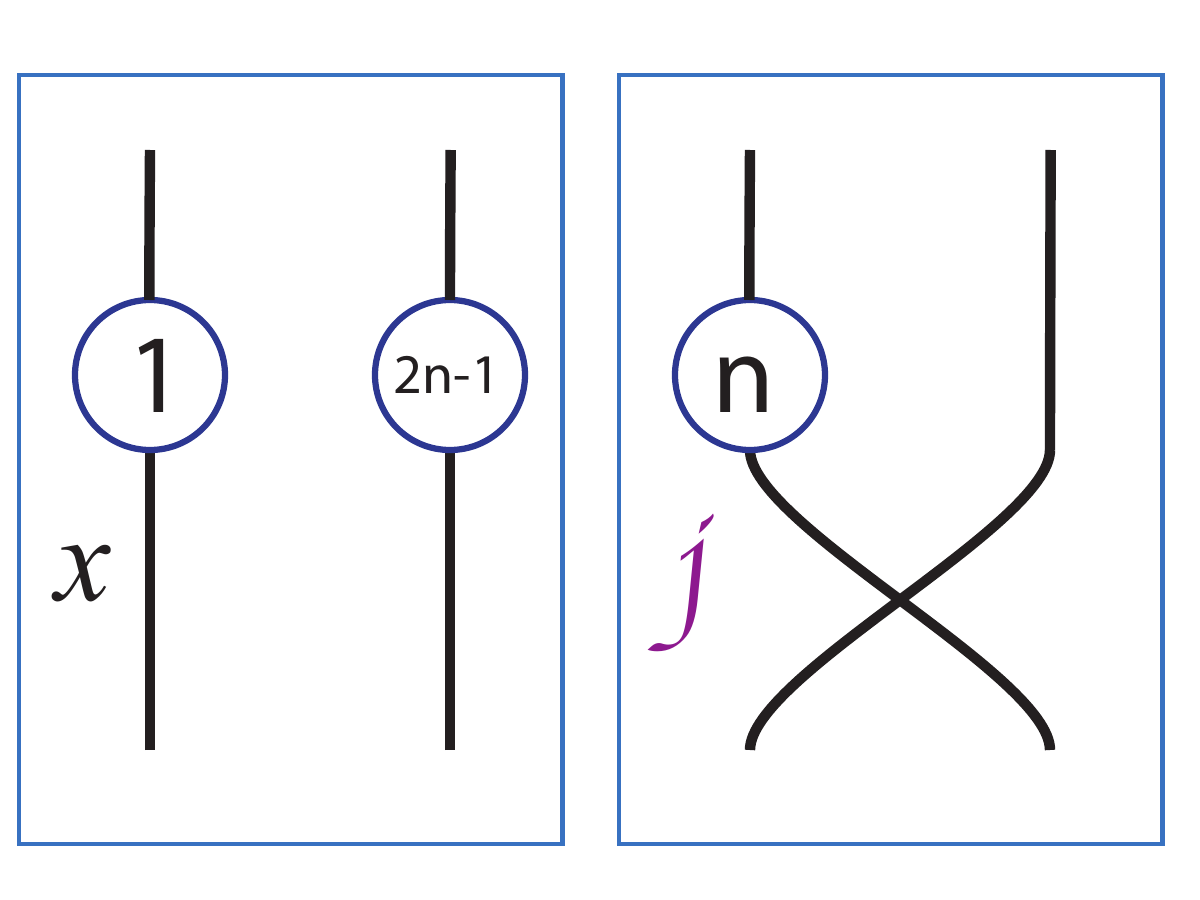}

Figure 18. Two string presentation of the dicyclic groups

\end{center}

\begin{center}
\includegraphics[width=0.6\paperwidth]{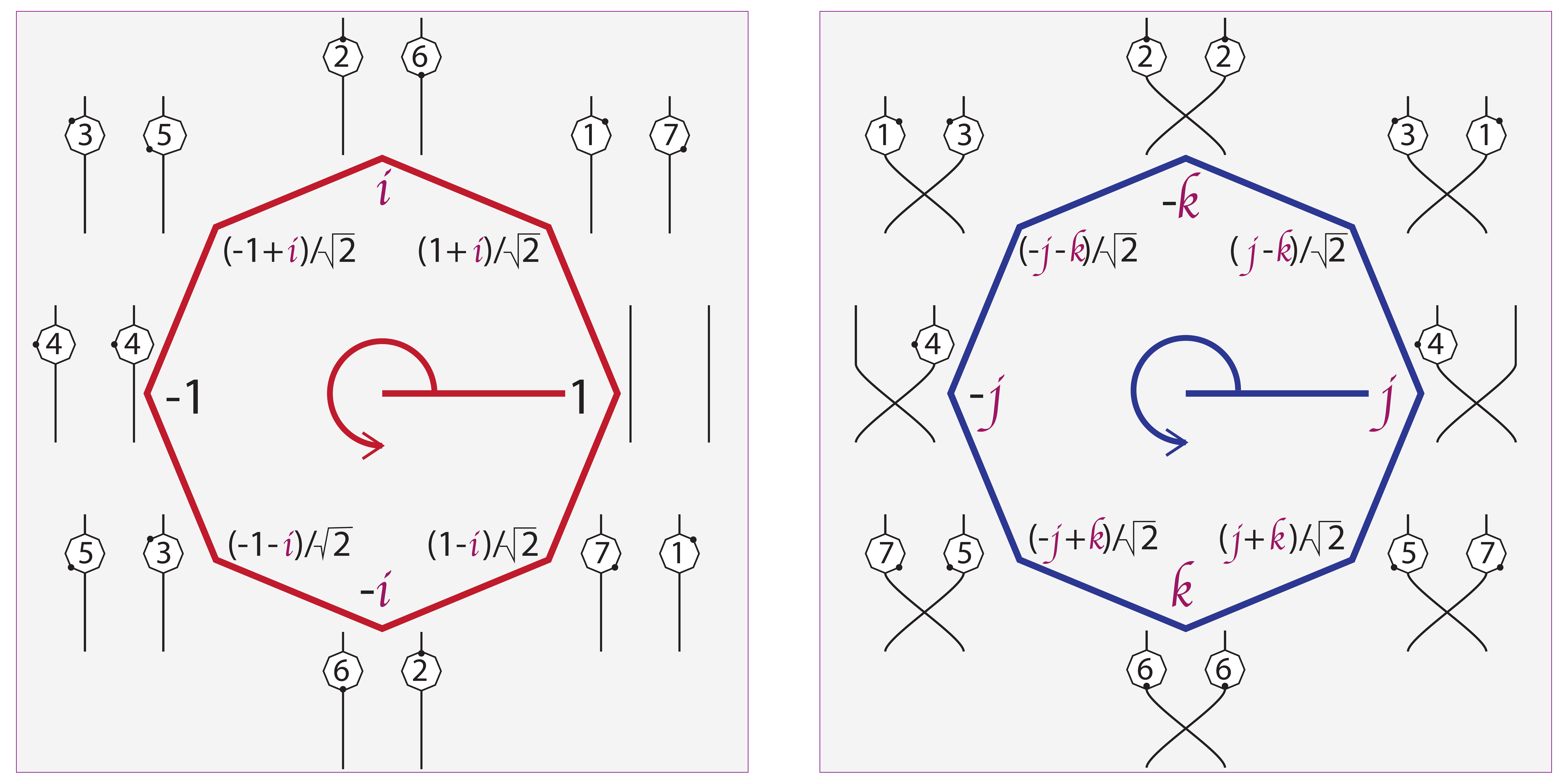}

Figure 19. The string-with-beads representations of the elements of the dicyclic group ${\mbox{\rm Dic}}_4$
\end{center}

This completes the proof of item (4) of Theorem~\ref{main}.

\section{The binary tetrahedral  group}
\label{S:BinTet}

The first example that we understood from the point of view of strings-with-beads was the binary tetrahedral group, $\widetilde{A_4}$,  which  is a twofold extension of the alternating group on four elements. The group is given via the presentation $\widetilde{A_4}=\langle a,b:  (ab)^2=a^3=b^3 \rangle$. It is a subgroup of $S^3$ via $a=1/2 (1 + {\boldsymbol i} + {\boldsymbol j}+ {\boldsymbol k})$, and $b=1/2 (1 + {\boldsymbol i} + {\boldsymbol j}- {\boldsymbol k}).$  and it is isomorphic to $\SL_2(\Z/3)$. This section contains a plethora of details about this group.

Subsection~\ref{SS:SL} describes the elements of $\SL_2(\Z/3)$  and their actions upon the eight non-zero vectors in $\Z/3 \times \Z/3.$ These are numbered by integers in $\{0,1,2,3,4,5,6,7\}$ such that antipodal points differ by four. This labeling gives a permutation representation that patently projects to the alternating group $A_4$. Subsection~\ref{SS:3}, describes a regular tetrahedron that has its vertices on diagonally opposed pairs of points on the cube $[-1,1]^3.$ This depiction allows for a representation of the group of orientation preserving symmetries into the group of $(3\times3)$-signed permutation matrices. 
In Section~\ref{SS:BinTet}, the binary tetrahedral group is described in terms of $a$ and $b$ as given above.  In the same subsection, the orders of the various element are catalogued and expressed as permutations on a set  of eight cosets. 

In Section~\ref{SS:Rib}, the depiction of the elements of $\widetilde{A_4}$  as strings-with-beads.  The composition is via vertical stacking, beads are allowed to move through crossings, and pairs of beads upon a single string cancel.  In this way the projection of elements in the binary extension $\widetilde{A_4}$ to the alternating group $A_4$ is achieved by ignoring the beads upon the strings. In Figs.~22 and ~23 the elements in $\widetilde{A_4}$ are illustrated as permutations, as elements of $\SL_2(\Z/3)$, and as strings-with-beads.

In subsection~\ref{SS:BTAlt}, some alternative permutation representations are presented. One resembles the dicyclic representation by examining the orbits of the ordered subgroup $S\cup{\boldsymbol j}S=[ (-1, -{\boldsymbol i}, 1, {\boldsymbol i}), (-{\boldsymbol j}, {\boldsymbol k},  {\boldsymbol j}, -{\boldsymbol k}) ].$ The subgroup $S\cup{\boldsymbol j}S=Q_8$ is, of course, the quaternion group. 
In subsection~\ref{SS:BTfin},  the cyclic subgroup $A=(1,a,a^2,-1,-a, -a^2)$ is used to give a different four-string representation.

\subsection{$\SL_2(\Z/3)$}
\label{SS:SL}
The group of determinant 1, $(2\times2)$-matrices with entries in the   $3$-element field is known as $\SL_2(\Z/3)$. It acts upon the vector space $\Z/3 \times \Z/3$. The entries in such a matrix can be taken from the set $\{-1,0,1\}$, where $1+1=-1$ and $(-1)^2=1=(-1)+(-1)$. Any matrix fixes the origin, and so the matrix can be thought of as a permutation of the remaining eight points. These vectors are written as rows and  numbered  follows:
\begin{center}
\begin{tabular}{||l||l||}\hline
$0 \leftrightarrow (1,-1)$ & $4 \leftrightarrow (-1,1)$ \\
$1 \leftrightarrow (0,-1)$  & $5 \leftrightarrow (0,1)$ \\
$2 \leftrightarrow (1,1)$ & $6 \leftrightarrow (-1,-1)$ \\
$3 \leftrightarrow (-1,0)$ & $7 \leftrightarrow (1,0)$ \\ \hline \end{tabular}
\end{center}
\vspace{.5cm}

The labels upon these vectors were chosen so that the isomorphism between $\SL_2(\Z/3)$ and the binary tetrahedral group $\widetilde{A_4}$ respects a permutation representation upon the latter that is induced via choosing a specific subgroup.

Here and elsewhere, in order to specify a group of symmetries as a permutation of a set, the elements of the set are labeled with non-negative (or sometimes strictly positive) integers. The group is generated by the four matrices that are given below. The matrix is specified, multiplied by a $(4\times 2)$ matrix whose columns correspond to $0,1,2,$ and $3$ respectively, and  the corresponding permutation of $\{0,1, \ldots, 7\}$ is indicated. The permutation is obtained since the points that are labeled $0$ and $4$ are antipodal, as are  the pairs $[1,5]$, $[2,6]$, and $[3,7]$; the transformations are linear.  

\begin{eqnarray*}
   \left[\begin{array}{rr} 1 & 0 \\ -1 & 1 \end{array}\right] \cdot 
   \left[\begin{array}{rrrr} 1 & 0 & 1 & -1 \\-1 & -1 & 1 & 0   \end{array}\right] 
   =  \left[\begin{array}{rrrr}1&0&1&-1 \\ 1 & -1 & 0 & 1  \end{array}\right] & \leftrightarrow & (027)(346);\\
      \left[\begin{array}{rr} 1 & 1 \\ 0 & 1 \end{array}\right] \cdot 
 \left[\begin{array}{rrrr} 1 & 0 & 1 & -1 \\-1 & -1 & 1 & 0   \end{array}\right]  
      =  \left[\begin{array}{rrrr} 0&-1&-1&-1 \\ -1 & -1 & 1 & 0  \end{array}\right] & \leftrightarrow & (016)(245);\\
 \left[\begin{array}{rr} -1 & 1 \\ -1 & 0 \end{array}\right] \cdot 
 \left[\begin{array}{rrrr} 1 & 0 & 1 & -1 \\-1 & -1 & 1 & 0   \end{array}\right] 
 =  \left[\begin{array}{rrrr} 1 & -1 &0& 1 \\ -1 & 0 & -1 & 1  \end{array}\right] & \leftrightarrow & (132)(576); \\
  \left[\begin{array}{rr} 0 & 1 \\ -1 & -1 \end{array}\right] \cdot 
 \left[\begin{array}{rrrr} 1 & 0 & 1 & -1 \\-1 & -1 & 1 & 0   \end{array}\right] 
  =  
  \left[\begin{array}{rrrr} -1 & -1 & 1 & 0 \\0 & 1 &1& 1 \end{array}\right] & \leftrightarrow & (035)(147).
 \end{eqnarray*}

\begin{center}
\includegraphics[width=0.3\paperwidth]{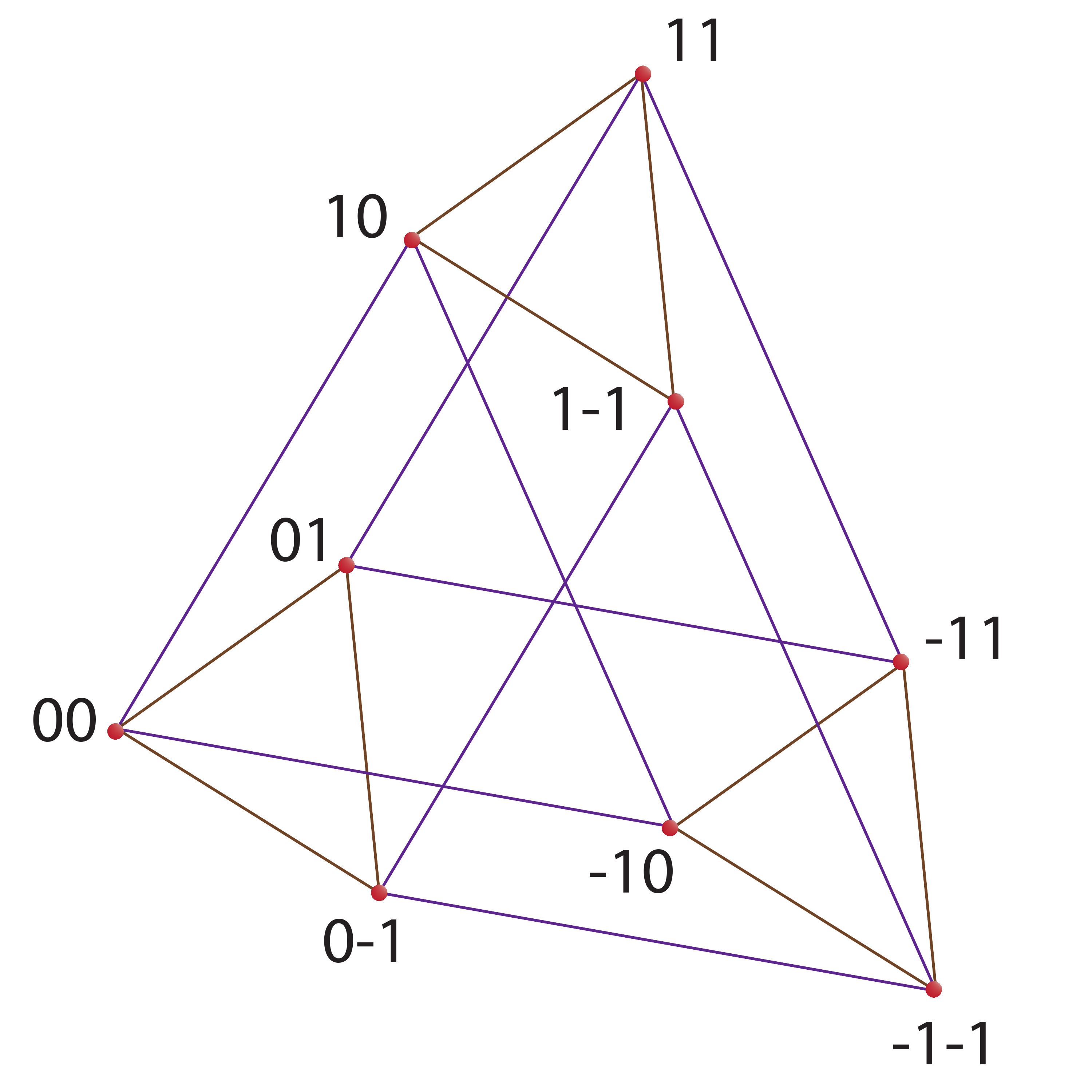}

Figure 20.The vector space $\Z/3 \times \Z/3$ as the vertices of  $\Delta \times \Delta$
\end{center}

\vspace{2cm}

The vector space $(\Z/3)^2$ consists of nine points. These are imagined as the vertices of the $4$-dimensional figure that is the cartesian product of a pair of equilateral triangles. An outline of the boundary of that $4$-dimensional figure is depicted in Fig.~20. Since the boundary of this space is a $3$-sphere, it can be decomposed as the union of two solid tori. The vertices lie on the common boundary of these. This single torus is unfolded onto the plane. In Fig.~21, the actions of half of the elements of $\SL_2(\Z/3)$ are depicted as permutations of the eight non-zero vectors. The matrices that are not shown are either inverses of those shown, or the identity matrix. Also within the figure, the permutations are expressed as products of cycles in the symmetric group on the set $\{0,1,\ldots, 7\}$.

The subgroup $\{I, -I\}$, where $I$ is the identity matrix, is normal. The quotient group ${\mbox{\rm PSL}}_2(\Z/3)=\SL_2(\Z/3)/\{I,-I\}$ is isomorphic to the alternating group on $4$-elements. There is a short exact sequence
$$0 \rightarrow \Z/2 \rightarrow \SL_2(\Z/3) \rightarrow A_4 \rightarrow 1.$$

\begin{center}
\includegraphics[width=0.6\paperwidth]{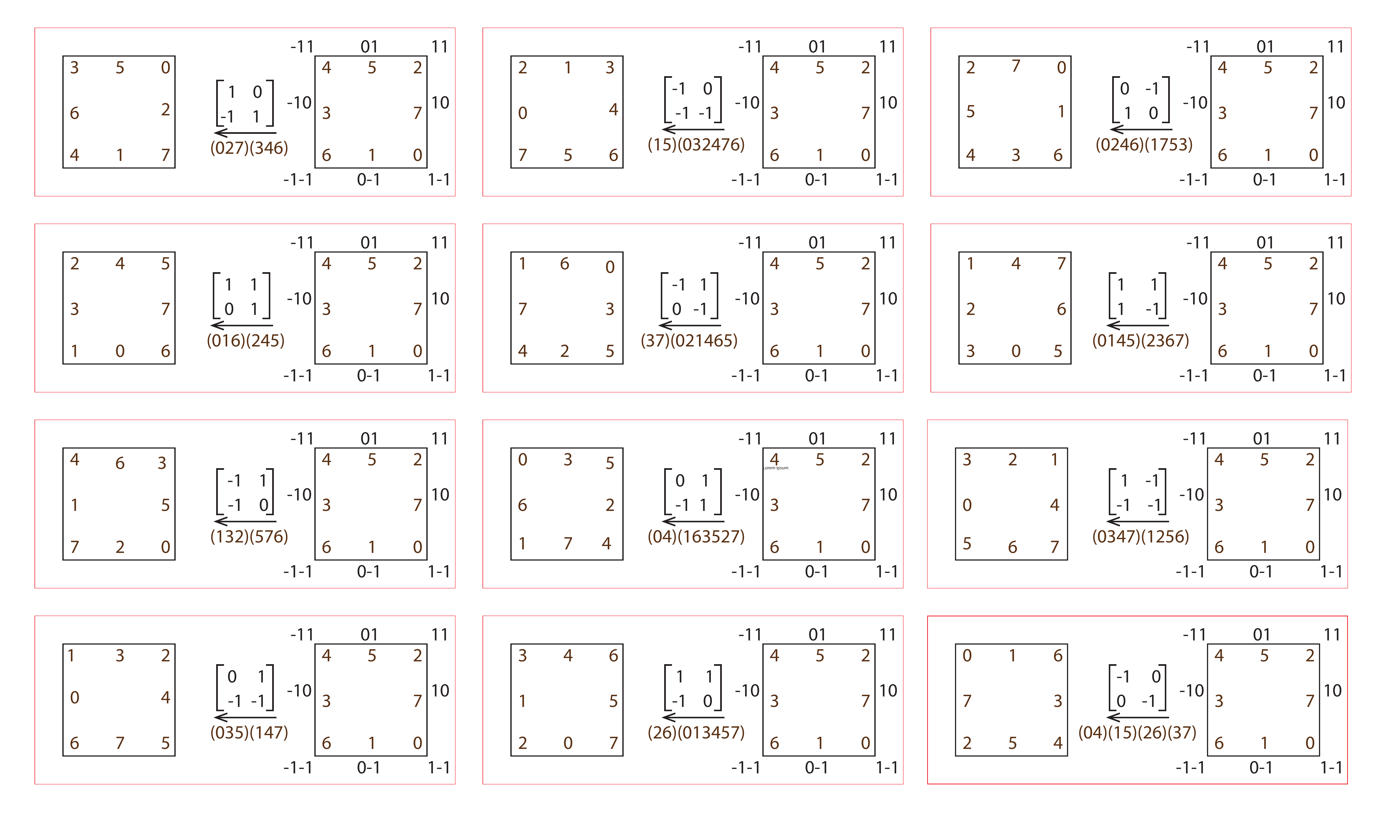}

Figure 21. The actions of half of the elements in $\SL_2(\Z/3)$
\end{center}

To see that ${\mbox{\rm PSL}}_2(\Z/3)$ is the alternating group, observe that the action of a $(2\times 2)$-matrix on the non-zero vectors in $(\Z/3)^2$ permutes the pairs of antipodal points. Specifically, the quotient map is given in terms of the permutation representation as follows. 
\begin{eqnarray*}
(0)(1)(2)(3) &\mapsfrom& \{(1), (04)(15)(26)(37) \}; \\
(01)(23) &\mapsfrom&  \{(0541)(2367), (0145)(2763) \}; \\ %
(02)(13) &\mapsfrom&  \{(0642)(1753), (0246)(1357) \}; \\ %
(03)(12) &\mapsfrom&  \{(0743)(1256), (0347)(1652) \}; \\ %
(021) &\mapsfrom&  \{(025)(146), (065421)(37) \}; \\ %
(012) &\mapsfrom&  \{(052)(164), (012456)(37)\}; \\ %
(013) &\mapsfrom&  \{(017)(345), (057413)(26) \}; \\ %
(031) &\mapsfrom&  \{(071)(354), (031475)(26) \} ;\\ %
(032) &\mapsfrom&  \{(036)(247), (076432)(15) \} ;\\  %
(023) &\mapsfrom&  \{(063)(274), (023467)(15) \}; \\ %
(123) &\mapsfrom&  \{(123)(567), (163527)(04) \}; \\ %
(132)&\mapsfrom&  \{(132)(576), (172536)(04) \} .\end{eqnarray*}

\begin{center}
\includegraphics[width=0.6\paperwidth]{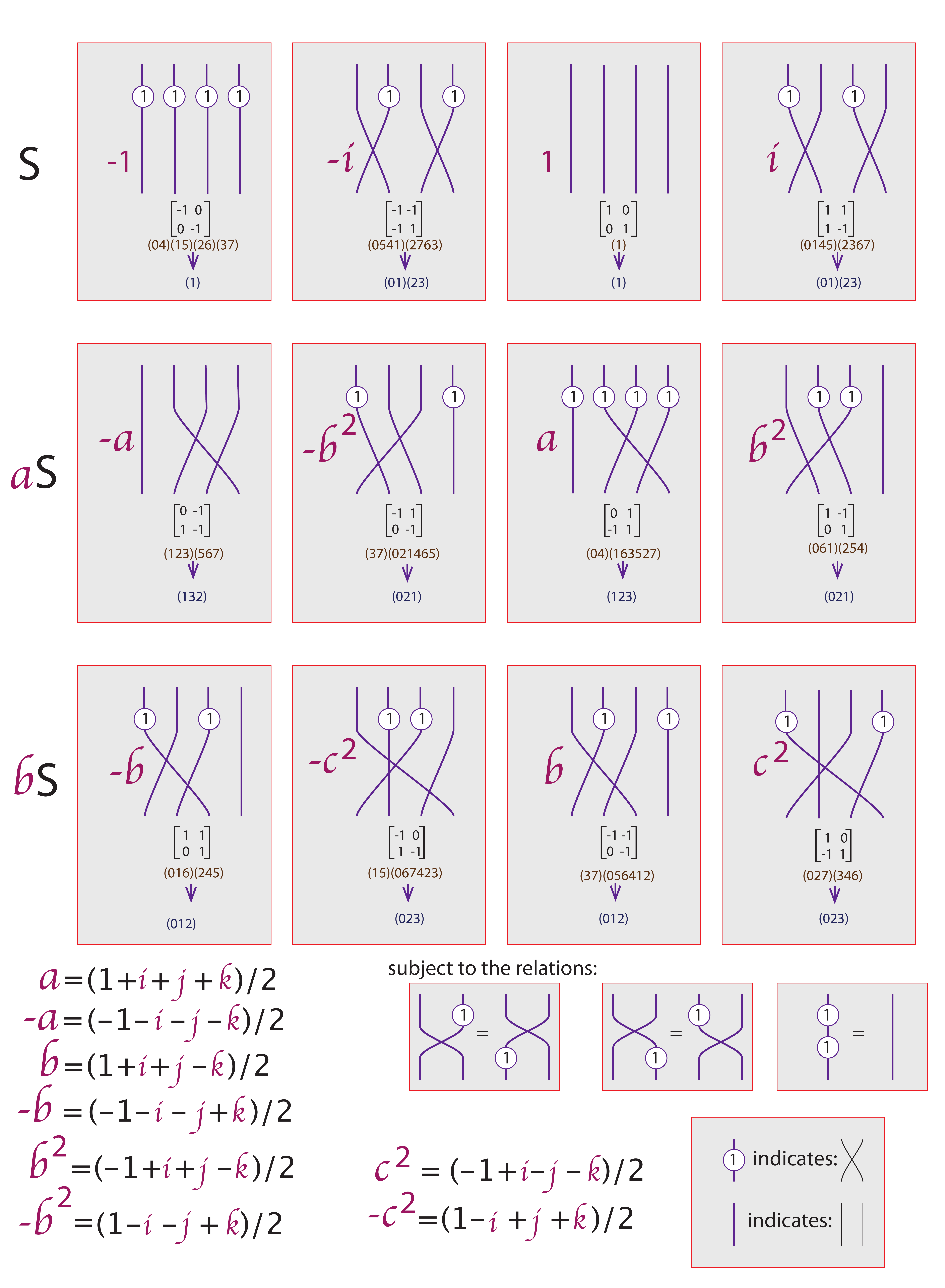}

Figure 22. The strings-with-beads representation of the elements of  $\SL_2(\Z/3),$ part 1
\end{center}

\newpage

\begin{center}
\includegraphics[width=0.6\paperwidth]{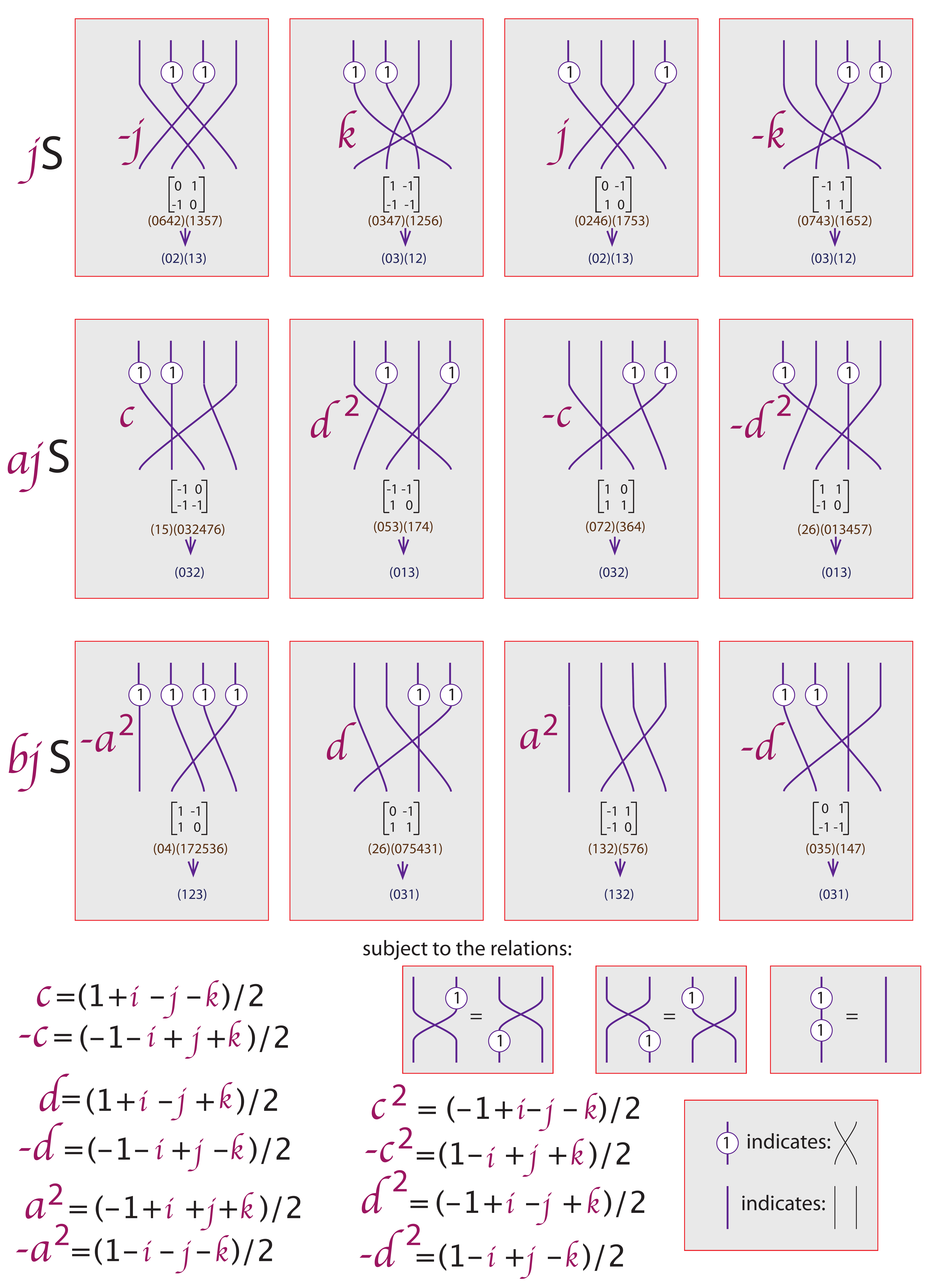}

Figure 23. The strings-with-beads representation of the elements of  $\SL_2(\Z/3),$ part 2
\end{center}

\newpage

A bit of book-keeping involves writing down both the inverse permutations and the inverse matrices to those depicted in Fig.~21.  These $24$  matrices and strings-with-beads representations are depicted in Fig.~22 and~23. That illustration also defines the elements of $\widetilde{A_4}$ and their correspondence to matrices in $\SL_2(\Z/3)$. Again, book-keeping demonstrates that matrices $M$ and $-M$ map to the same permutation in $A_4$.  Furthermore, by reducing the entries in the permutations modulo $4$ (for example, $(04)(172536) \mapsto (132)$) the quotient map above is obtained. 

\begin{center}
\rule{3in}{0.005in}
\end{center}

Each of the rows of Figs.~22 and~23 represents an ordered coset of the subgroup $S= (-1, - {\boldsymbol i}, 1, {\boldsymbol i})$ of the quaternions
\[Q_8= S\cup {\boldsymbol j}S =[ (-1, - {\boldsymbol i}, 1, {\boldsymbol i}), (-{\boldsymbol j},{\boldsymbol k}, {\boldsymbol j}, -{\boldsymbol k})]\] which was used to give a strings-with-beads presentation of $Q_8$ above.

\subsection{A regular tetrahedron and the alternating group}\label{SS:3}

A regular tetrahedron sits inside the cube $[-1,1]^3 \subset \R^3$ with vertices $0\leftrightarrow(1,1,1)$, $1\leftrightarrow(1,-1,-1)$, $2\leftrightarrow (-1,1,-1)$, and $3\leftrightarrow(-1,-1,1)$. The edges form diagonals of the faces of the cube, and therefore their lengths are each $\sqrt{2}.$ An illustration is given in Fig.~24. Under this labeling, a direct computation indicates the correspondences between  signed permutation matrices  of  determinant $1$  and elements of the alternating group indicated below.

\rule{0pt}{1em}

\begin{tabular}{llcrclcr} 
\rule{-30pt}{0in} &
$ \left[\begin{array}{rrr} 0 & 0 & 1 \\ 1 & 0 & 0 \\ 0 & 1 &0 \end{array}\right]^{\pm 1}$   
& $\leftrightarrow$ 
&   $(1,2,3)^{\pm 1};   $ 
& \rule{-12pt}{0in} 
&
$  \left[\begin{array}{rrr} 0 & 0 & -1 \\ -1 & 0 & 0 \\ 0 & 1 &0  \end{array}\right]^{\pm 1}$  
& $\leftrightarrow$ 
&   $(0,3,2)^{\pm 1}; $ \\
\rule{-30pt}{0in} &
$  \left[\begin{array}{rrr} 0 & 0 & 1 \\ -1 & 0 & 0 \\ 0 & -1 &0  \end{array}\right]^{\pm 1}$  
& $\leftrightarrow $ 
&  $(0,1,3)^{\pm 1};  $ 
&\rule{-12pt}{0in} 
&
$  \left[\begin{array}{rrr} 0 & 0 & -1 \\ 1 & 0 & 0 \\ 0 & -1 &0  \end{array}\right]^{\pm 1}$  
&$ \leftrightarrow$
&   $(0,2,1)^{\pm 1}; $  \\
\rule{-30pt}{0in} &
$  \left[\begin{array}{rrr} 1 & 0 & 0  \\ 0 & -1 & 0 \\ 0 & 0 &-1  \end{array}\right]$  
& $\leftrightarrow$ 
&   $(0,1)(2,3); $ 
&\rule{-12pt}{0in} 
&
$ \left[\begin{array}{rrr} -1 & 0 & 0  \\ 0 & 1 & 0 \\ 0 & 0 &-1  \end{array}\right]^{\rule{0in}{1ex}}$   
& $\leftrightarrow$ 
&  $(0,2)(1,3); $ \\
\rule{-30pt}{0in} &
$  \left[\begin{array}{rrr} -1 & 0 & 0  \\ 0 & -1 & 0 \\ 0 & 0 &1  \end{array}\right]^{\rule{0in}{1ex}}$  
& $\leftrightarrow$ 
&  $(0,3)(1,2);  $ 
& \rule{-12pt}{0in} 
&
$ \left[\begin{array}{rrr} 1 & 0 & 0  \\ 0 & 1 & 0 \\ 0 & 0 & 1  \end{array}\right]$   
&$\leftrightarrow$ 
&  $(0)(1)(2)(3). $
\end{tabular}
\vspace*{.5cm}

These signed permutation matrices are in the image of the restriction of the projection $p:\widetilde{A_4} \rightarrow \SO(3)$  that was given in Section~\ref{S:Dic}. 

\newpage

\begin{center}
\includegraphics[width=0.3\paperwidth]{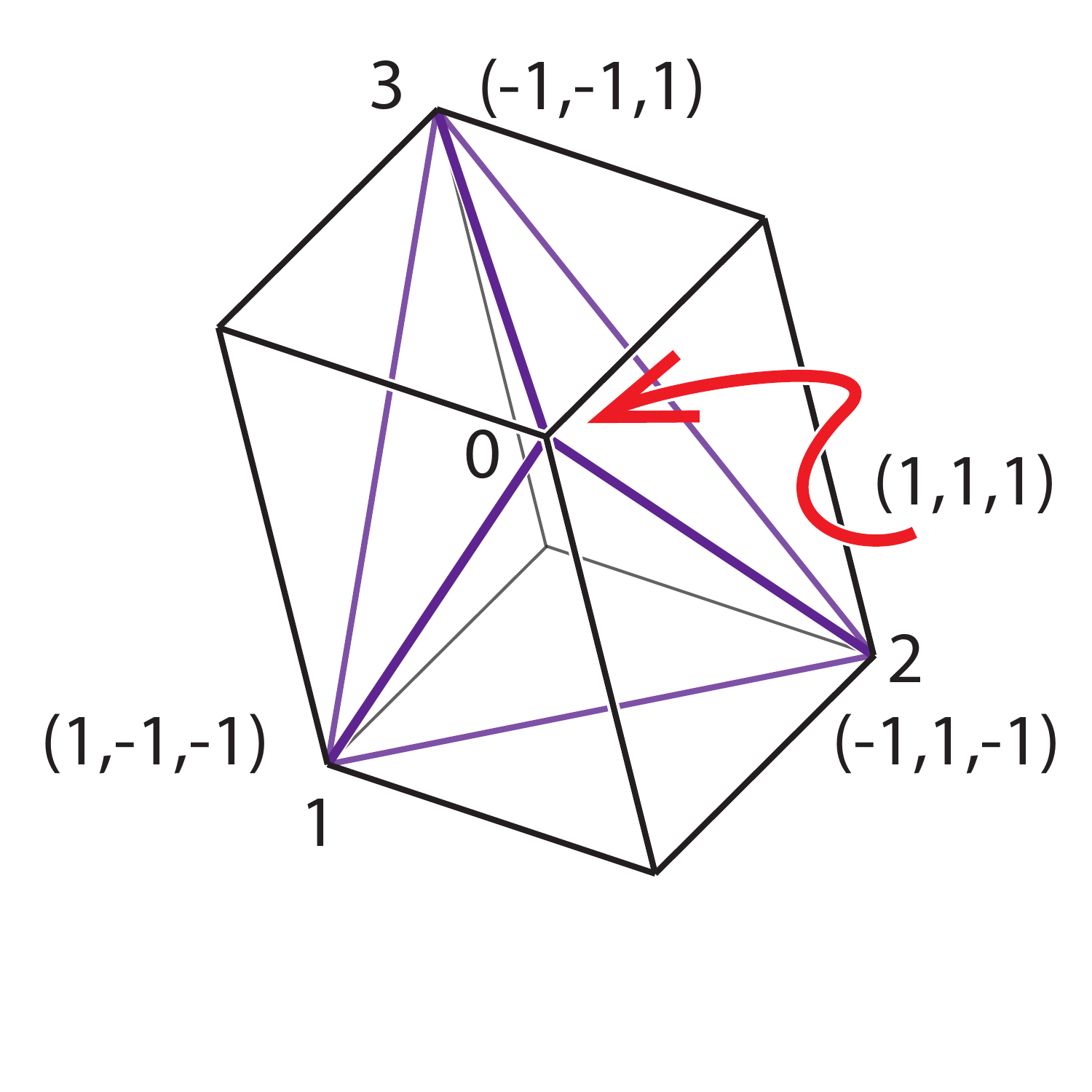}

Figure 24. A regular tetrahedron in the cube $[-1,1]\times[-1,1]\times[-1,1]$
\end{center}

\subsection{The binary tetrahedral group}
\label{SS:BinTet}

The binary tetrahedral group, $\widetilde{A_4}$, consists of the $24$ points in the $3$-dimensional sphere $S^3$ that are of the form
$$\pm 1, \pm  {\boldsymbol i}, \pm {\boldsymbol j}, \pm {\boldsymbol k}, \frac{1}{2}(\pm1\pm{\boldsymbol i}\pm {\boldsymbol j}\pm{\boldsymbol k}).$$
All  sign combinations are possible.  These are the vertices of a platonic hypersolid that has $24$ vertices, $96$ edges, $96$ triangular faces, and $24$ three dimensional faces. As such, it is self-dual. The $(2\times4)$-matrix
$\left[ \begin{array}{cccr} 3&2&0&2\\ 0&2&3&-2 \end{array}\right]$ projects the vertices and edges to the graph that is drawn in Fig.~25. Therein the vertices  that are labeled 
$$
a=\frac{1}{2} (1+{\boldsymbol i}+ {\boldsymbol j}+{\boldsymbol k}), \quad
b=\frac{1}{2} (1+{\boldsymbol i}+ {\boldsymbol j}-{\boldsymbol k}), $$
$$c=\frac{1}{2} (1+{\boldsymbol i}- {\boldsymbol j}-{\boldsymbol k}), \quad {\mbox{\rm and}} \quad 
d=\frac{1}{2} (1+{\boldsymbol i}- {\boldsymbol j}+{\boldsymbol k})$$
have order $6$. Within the figure $a$ is written as a vector $\frac{1}{2}(1,1,1,1)$, for example, to save space and visual clutter. One can compute directly that $a^3=b^3=c^3=d^3=-1$. The non-trivial powers of these elements are labeled within the figure. 
Let $N=(1,a^2,a^4)$; this is a subgroup of $\widetilde{A_4}$ that we have ordered as indicated. The set of cosets can be listed, ordered, and numbered as indicated below.

\vspace*{0.5cm}\begin{center}\begin{tabular}{|ll|}\hline
$\rule{0in}{1.25em} 0 \leftrightarrow N=(1,a^2,-a)$ & $4 \leftrightarrow aN=(a,-1,-a^2)$ \\
$1 \leftrightarrow {\boldsymbol i}N=({\boldsymbol i},-b,-d^2)$ & $5 \leftrightarrow {\boldsymbol i}aN=(d^2,-{\boldsymbol i},b)$\\
$2 \leftrightarrow {\boldsymbol j}N=({\boldsymbol j},c^2,-b^2)$ & $6 \leftrightarrow {\boldsymbol j}aN=(b^2,-{\boldsymbol j},-c^2)$\\  $3 \leftrightarrow  {\boldsymbol k}N=({\boldsymbol k},-d,c)$ &  $7 \leftrightarrow {\boldsymbol k}aN=(-c,-{\boldsymbol k},d)$ \\ \hline \end{tabular}\end{center}

\begin{center}
\includegraphics[width=0.6\paperwidth]{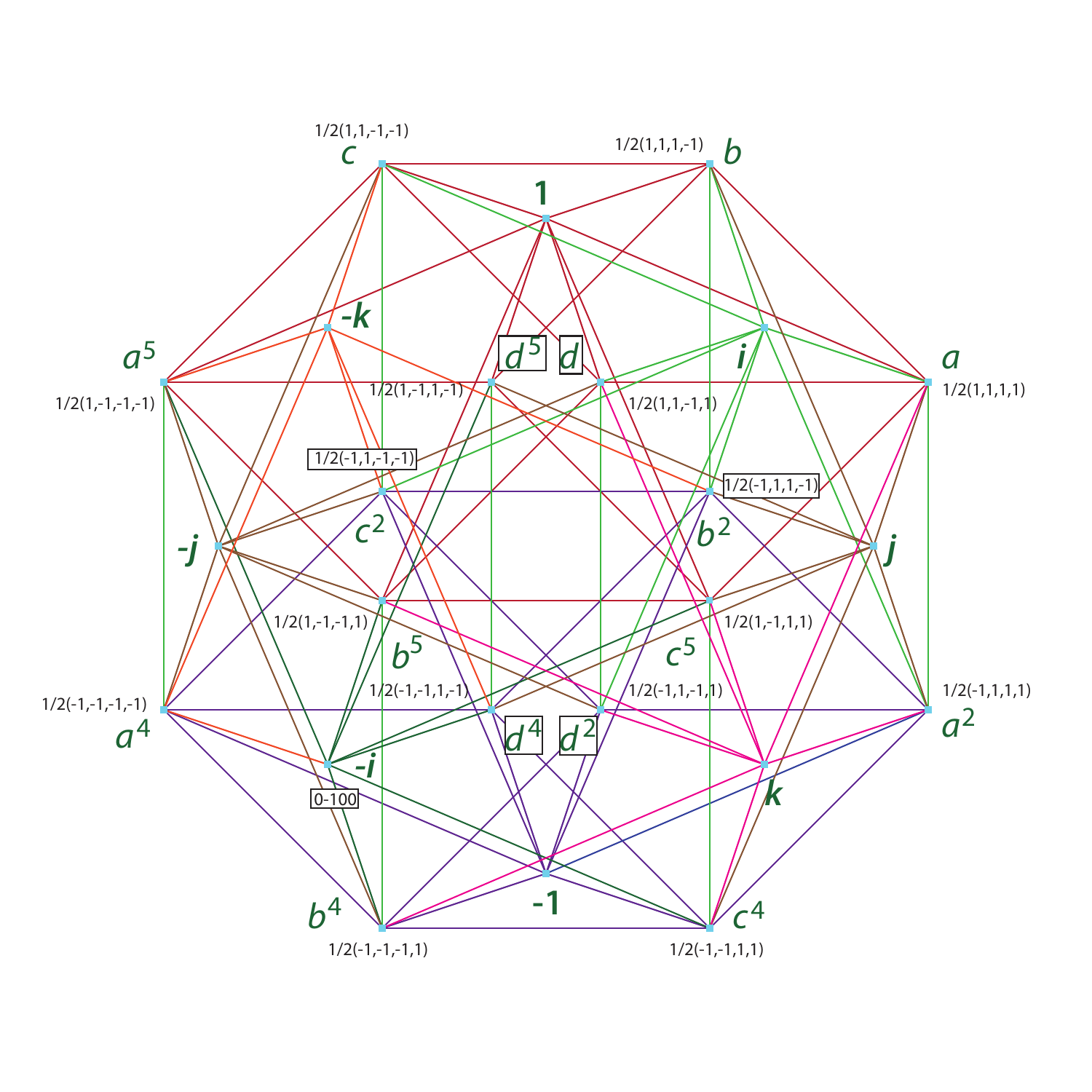}

Figure 25. A drawing of the $1$-skeleton of the $24$-cell
\end{center}

While we do not use the ordering here, we will use it below.
Under this labeling,  the actions of $a$ through $d$ upon the cosets are computed as follows:
\begin{eqnarray*}
a & \leftrightharpoons & (04)(163527); \\
b & \leftrightharpoons & (37)(056412); \\
c & \leftrightharpoons & (15)(032476); \\
d & \leftrightharpoons & (26)(075431).\end{eqnarray*}
Finally, one can compute directly (by reducing the integers mod 4),  or from the formula for the projection $p$ that is given above, that the actions of these elements upon the regular tetrahedron is given as
$$p'(a)=p'(-a)=(1,2,3), \quad p'(b)=p'(-b)=(0,1,2),$$  $$p'(c)=p'(-c)=(0,3,2), \quad {\mbox{\rm and}} \quad p'(d)=p'(-d)=(0,3,1).$$ Here $p'$ indicates the composition of the projection $p$ with the correspondence that indicates the action of signed permutation matrices upon the tetrahedron indicated in Fig.~24 Meanwhile, ${\boldsymbol i},{\boldsymbol j},{\boldsymbol k}$ map via $p$ to diagonal matrices under $p$, with $p'({\boldsymbol i})=(01)(23)$, $p'({\boldsymbol j})=(02)(13)$, and $p'({\boldsymbol k})=(03)(12)$.

\begin{center}
\includegraphics[width=0.6\paperwidth]{Gen2Els.pdf}

Figure 26. Some products in $\widetilde{A_4}$
\end{center}

\subsection{More about the strings-with-beads representation}
\label{SS:Rib}

Each of the elements in the binary tetrahedral group $\widetilde{A_4}\approx \SL_2(\Z/3)$ has a permutation representation in the  group, $\Sigma_8$, of symmetries of $\{0,1,\ldots, 7\}$. 
To construct the string-with-beads representation, these elements were 
arranged along a horizontal line segment in the order and grouping $(04)(15)(26)(37)$. Then the associated permutation was drawn in a braid-like fashion. In that way, each element in the group was first imagined as  consisting of four ribbons that have an even number of half-twists among them. Since permutations, rather than braids or ribbons, are being considered, two half-twists between a successive pair of strings cancel: they are of the form, for example, $(04)(04)$ which equals $1$ in the permutation group. Each ribbon (or pair of successive strings) is abstracted to a thick line, and the twists are abstracted to beads. Figs.~22 and~23 contain legends at the bottom.

The group operation, then,  is vertical juxtaposition of diagrams. Fig.~26 contains a number of sample computations that indicate how to write several of the elements in the group in terms of generators. For example, the product $ba$ is written with $b$ juxtaposed above $a$. A pair of beads upon a single string cancel (in this case upon the $2$\/nd and $4$\/th string). Redundant crossings are removed. Thus the $2$\/nd string at the bottom only crosses the first and the two crossings between it and the $3$\/rd are removed. The resulting simplified element is ${\boldsymbol i}$. All of ${\boldsymbol i}$, ${\boldsymbol j}$, ${\boldsymbol k}$, $c$, and $d$ can be seen to be expressed in terms of products of $a$ and $b$. 

It is worth pointing out that to compile these figures, the authors worked through products such as $ba = \frac{1}{2} (1+{\boldsymbol i}+ {\boldsymbol j}-{\boldsymbol k})\cdot
\frac{1}{2} (1+{\boldsymbol i}+ {\boldsymbol j}+{\boldsymbol k})$ by employing the distributive law repeatedly and simplifying  expressions  such as ${\boldsymbol i}{\boldsymbol j} = {\boldsymbol k}$. Such computations are subject to various typographic and transcription errors. Meanwhile, various note-taking and graphical softwares allow for easy implementations of the  associated diagrammatic computations.

This completes the description of $\widetilde{A_4}$ as a subgroup of $(\Z/2)^4 \rtimes A_4$ as given in item (5) of Theorem~\ref{main}.

\subsection{Some alternative depictions}
\label{SS:BTAlt}

\begin{center}
\includegraphics[width=0.5\paperwidth]{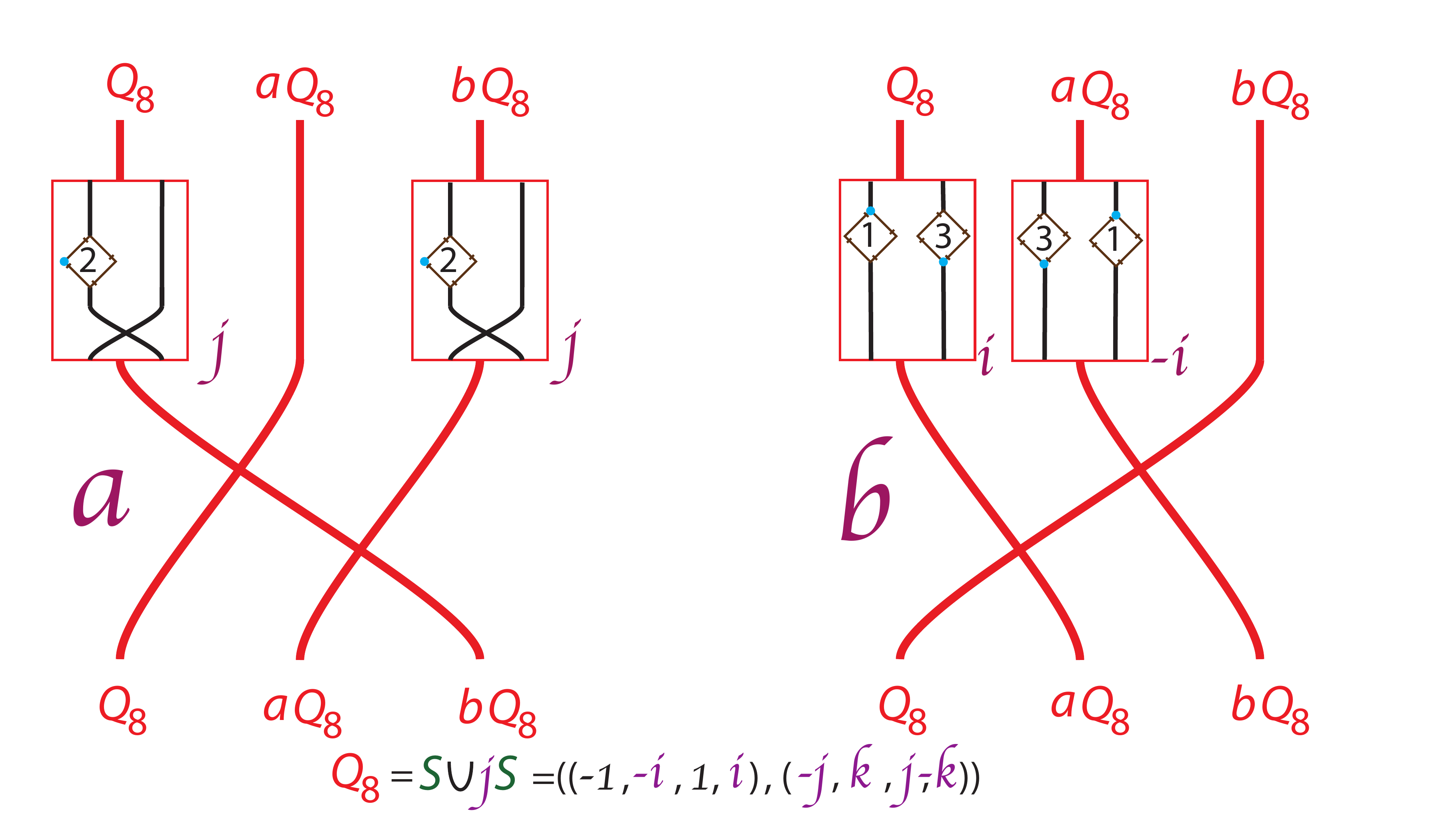}

Figure 27. The generators of $\widetilde{A_4}$ as elements in $(Q_8)^3 \rtimes \Z/3$
\end{center}

Since  $Q_8=\{ \pm 1,$ $ \pm {\boldsymbol i}, $ $ \pm {\boldsymbol j},$  $\pm {\boldsymbol k}\}$ is a subgroup of the binary tetrahedral group $\widetilde{A_4}$,  the unit quaternionic group $Q_8$ will be ordered (and partitioned) as  \[Q_8= S \cup {\boldsymbol j}S= [ (-1, -{\boldsymbol i}, 1, {\boldsymbol i}), (-{\boldsymbol j}, {\boldsymbol k}, {\boldsymbol j}, -{\boldsymbol k})].\] Then the binary tetrahedral group  will be written as the (ordered) set of cosets 
$\widetilde{A_4}= (Q_8,aQ_8,bQ_8)$.   The actions of $a$ and $b$ upon these ordered cosets will be computed, and an alternative  diagrammatic description  of those elements, depicted in Fig.~27, results from the computations below. An analogue of Fig.~26, is given in Figs.~28 and~29, and a compilation of the elements is presented in Figs.~30, and ~31.

  The cosets are ordered and computed as
\[\begin{array}{ll} S=(-1, -{\boldsymbol i}, 1, {\boldsymbol i}); &  {\boldsymbol j}S= (-{\boldsymbol j}, {\boldsymbol k}, {\boldsymbol j}, -{\boldsymbol k});\\
aS=(-a,-b^2,a,b^2); &a{\boldsymbol j}S=(c,d^2,-c,-d^2);\\  
 bS=(-b, -c^2,b,c^2);  &  b{\boldsymbol j}S=(-a^2,d,a^2,-d). \end{array}\]

Direct computations yield:
\[\begin{array}{ll} a^2(S)=(-a^2,d,a^2,-d); & a^2{\boldsymbol j}S=(b,c^2,-b,-c^2); \\ 
a(bS)={\boldsymbol j}S=( -{\boldsymbol j}, {\boldsymbol k}, {\boldsymbol j},-{\boldsymbol k}); & 
a(b{\boldsymbol j}S)= -S= (1, {\boldsymbol i}, -1, -{\boldsymbol i});\\
b(aS)={\boldsymbol i}S= (-{\boldsymbol i},1, {\boldsymbol i}, -1); & 
b(a{\boldsymbol j}S)= {\boldsymbol k}S=(-{\boldsymbol k}, -{\boldsymbol j},{\boldsymbol k}, {\boldsymbol j});\\
b^2S= (-b^2,a, b^2,-a); & 
b(b{\boldsymbol j}S)=(-d^2,c,d^2,-c). \end{array}\]
For notational convenience, we will write, for example, $b(aS)= [1/4, S]$ to indicate that as a set, $abS=S$, and there is a $1/4$ rotation in the cyclic order of the elements. The actions of $a$ and $b$ upon the cosets are given as:
\[\begin{array}{ll} a(S)=[0/4,aS]; & a({\boldsymbol j}S)= [0 /4,a{\boldsymbol j}S ];\\
a^2S=[0/4,b{\boldsymbol j}S ];  & a^2{\boldsymbol j}S=[2/4,bS]; \\
a(bS)=[0/4,{\boldsymbol j}S ]; &  a(b{\boldsymbol j}S)=[2/4,S].   \\
 b(S)=[ 0/4,bS ];  & b({\boldsymbol j}S)=[0/4 ,b{\boldsymbol j}S ] \\
b(aS)=[1/4,S]; & ba{\boldsymbol j}S ={\boldsymbol k}S=[3/4,S];  \\
b(bS)=[1/4,aS];  & b(b{\boldsymbol j}S)=[3/4,a{\boldsymbol j}S]. \end{array}
\]

\begin{center}
\includegraphics[width=0.65\paperwidth]{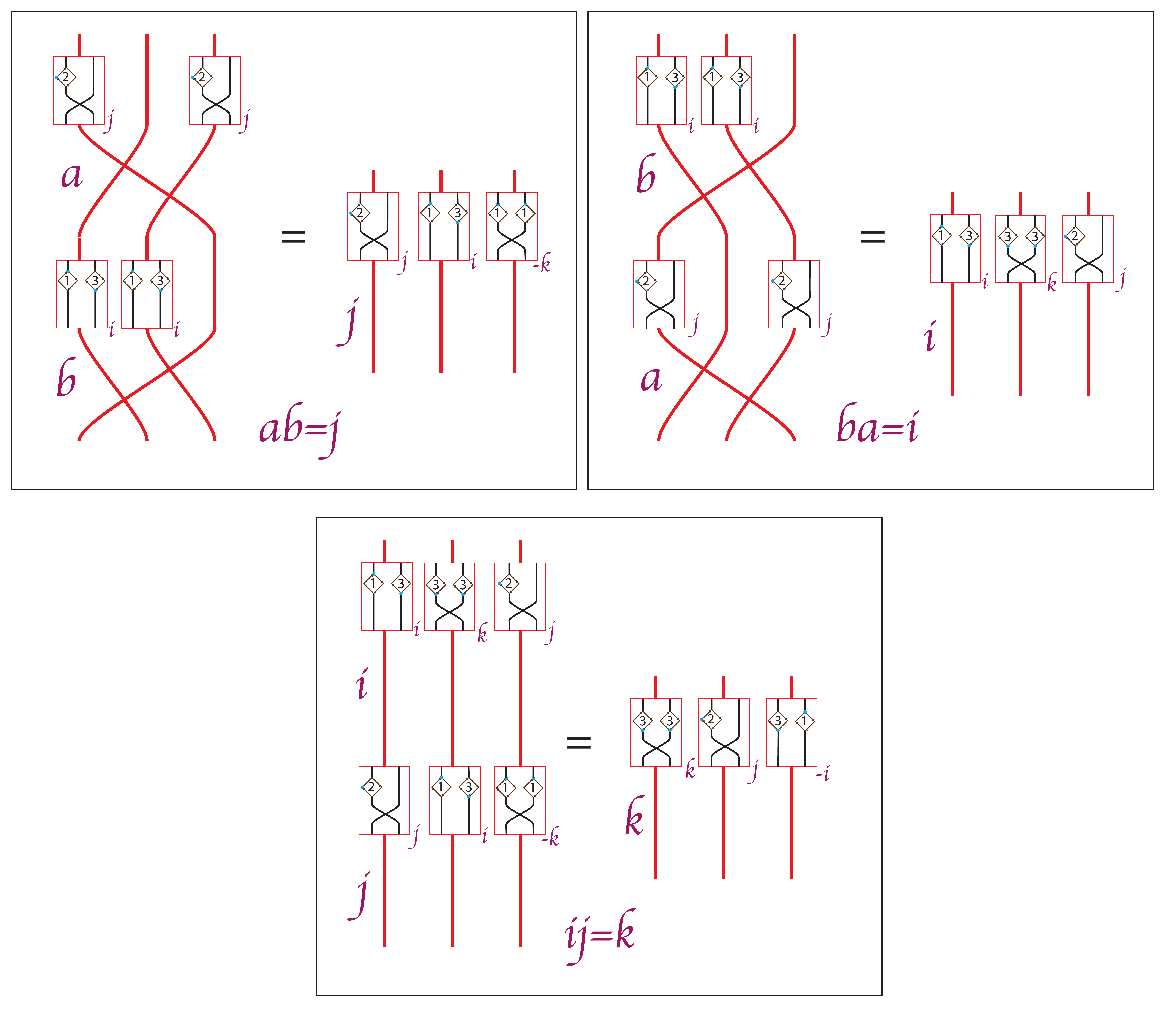}

Figure 28. The relations $ab={\boldsymbol j}$, $ba={\boldsymbol i}$, and ${\boldsymbol i}{\boldsymbol j}={\boldsymbol k}$.
\end{center}

\begin{center}
\includegraphics[width=0.65\paperwidth]{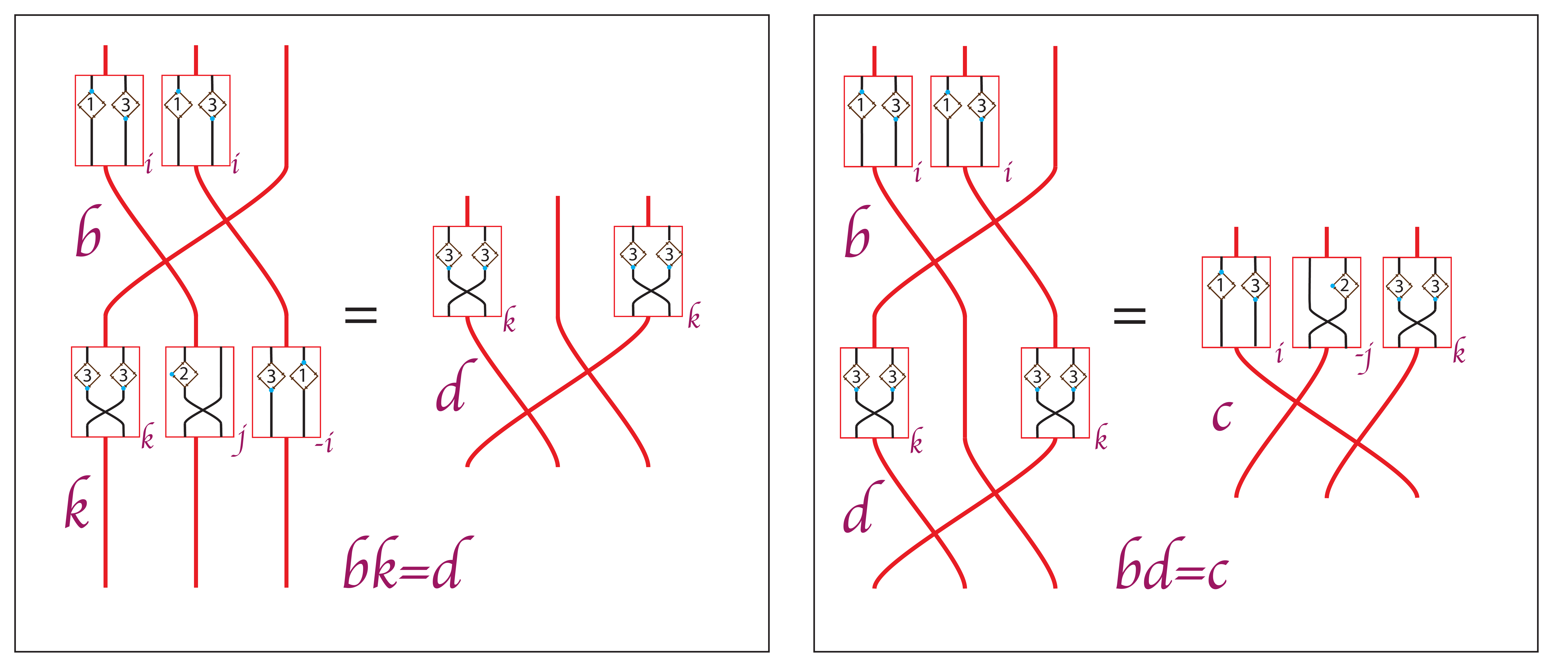}

Figure 29.  The relations $b{\boldsymbol k}=d$, and $bd=c$
\end{center}

In Figs.~30 and~31, each element is depicted as three strings-with-beads. The beads in this case are elements of the unit quaternions $Q_8$. Since $Q_8$ is a normal subgroup of $\widetilde{A_4}$ with quotient $\Z/3$, there is a further matrix representation of $\widetilde{A_4}$. The elements $a$, $b$, $c$ and $d$ can be written as follows:
\[ a \leftrightharpoons \left[ \begin{array}{ccc} 0 & 0 & {\boldsymbol j} \\ 1 & 0 & 0 \\ 0 & {\boldsymbol j} & 0  \end{array} \right];
\quad b \leftrightharpoons 
\left[ \begin{array}{ccc} 0 &  {\boldsymbol i } & 0  \\  0 & 0 &  {\boldsymbol i }  \\ 1 & 0 &  0 \end{array} \right]; \] 
\[  c \leftrightharpoons \left[ \begin{array}{ccc} 0 & 0 & {\boldsymbol i} \\ -{\boldsymbol j} & 0 & 0 \\ 0 & {\boldsymbol k} & 0  \end{array} \right];
\quad d \leftrightharpoons 
\left[ \begin{array}{ccc} 0 &  {\boldsymbol k } & 0  \\  0 & 0 &  1  \\ {\boldsymbol k} & 0 &  0 \end{array} \right]. 
\]
Furthermore, ${\boldsymbol i}, {\boldsymbol j},$ and ${\boldsymbol k}$ have the following representations:
\[ 
{\boldsymbol i}  \leftrightharpoons \left[ \begin{array}{ccc}  {\boldsymbol i}  & 0 & 0 \\  0 & {\boldsymbol k} & 0 \\ 0 & 0 & {\boldsymbol j}   \end{array} \right]; \quad
{\boldsymbol j}  \leftrightharpoons \left[ \begin{array}{ccc}  {\boldsymbol j}  & 0 & 0 \\  0 & {\boldsymbol i} & 0 \\ 0 & 0 & -{\boldsymbol k}   \end{array} \right]; \quad
 {\boldsymbol k}  \leftrightharpoons \left[ \begin{array}{ccc}  {\boldsymbol k}  & 0 & 0 \\  0 & {\boldsymbol j} & 0 \\ 0 & 0 & -{\boldsymbol i}   \end{array} \right]. \]

The reader is encouraged  to check that the relations depicted in Figs.~28 and~29  hold among these matrices. 

This completes the description of  $\widetilde{A_4} \approx Q_8 \rtimes (\Z/3)$ --- item (6) in Theorem~\ref{main}. See also \cite{wiki:BinTet}.

\begin{center}
\includegraphics[width=0.55\paperwidth]{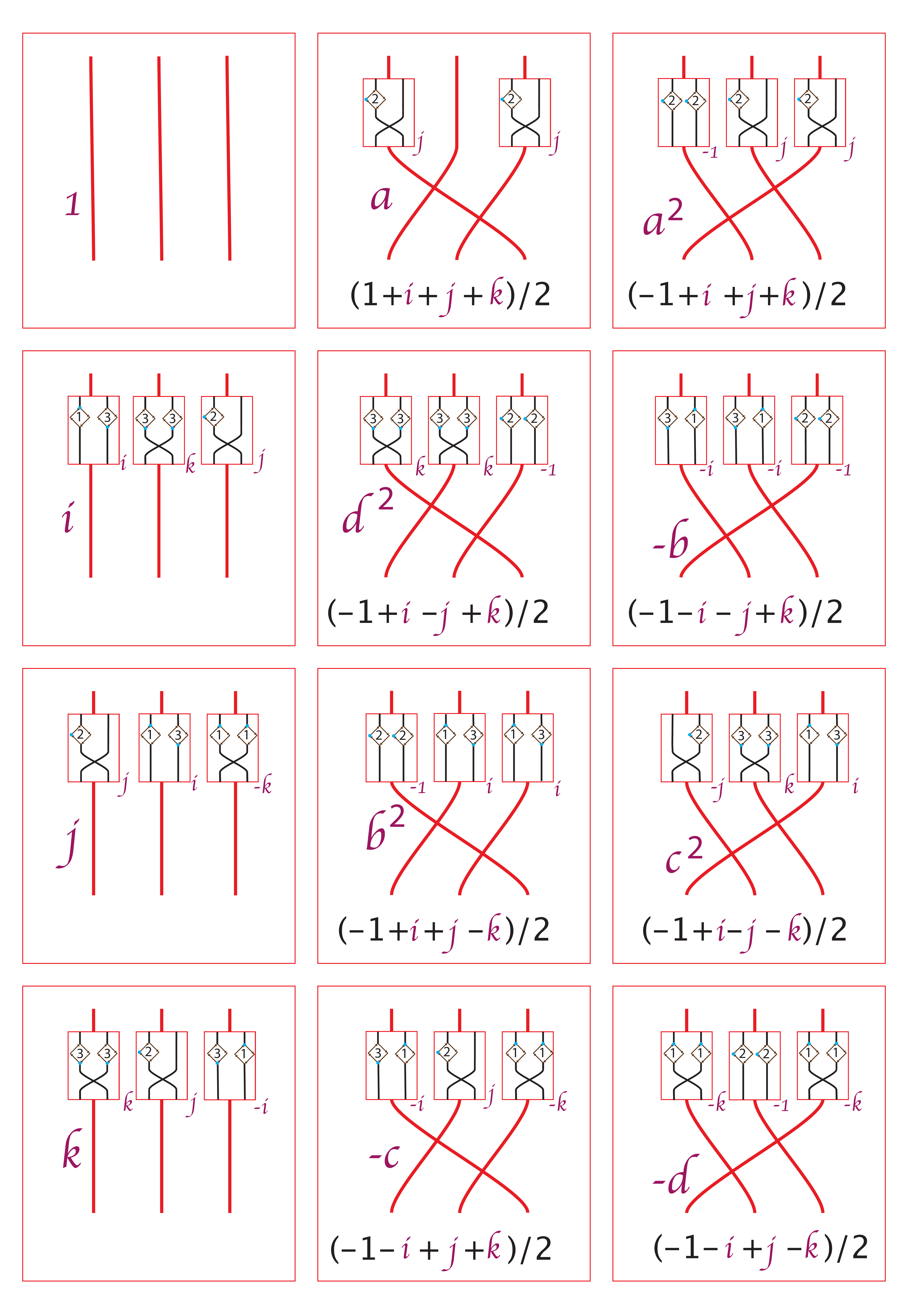}

Figure 30. The binary tetrahedral group  as it acts on cosets of $Q_8$, part 1
\end{center}

\begin{center}
\includegraphics[width=0.55\paperwidth]{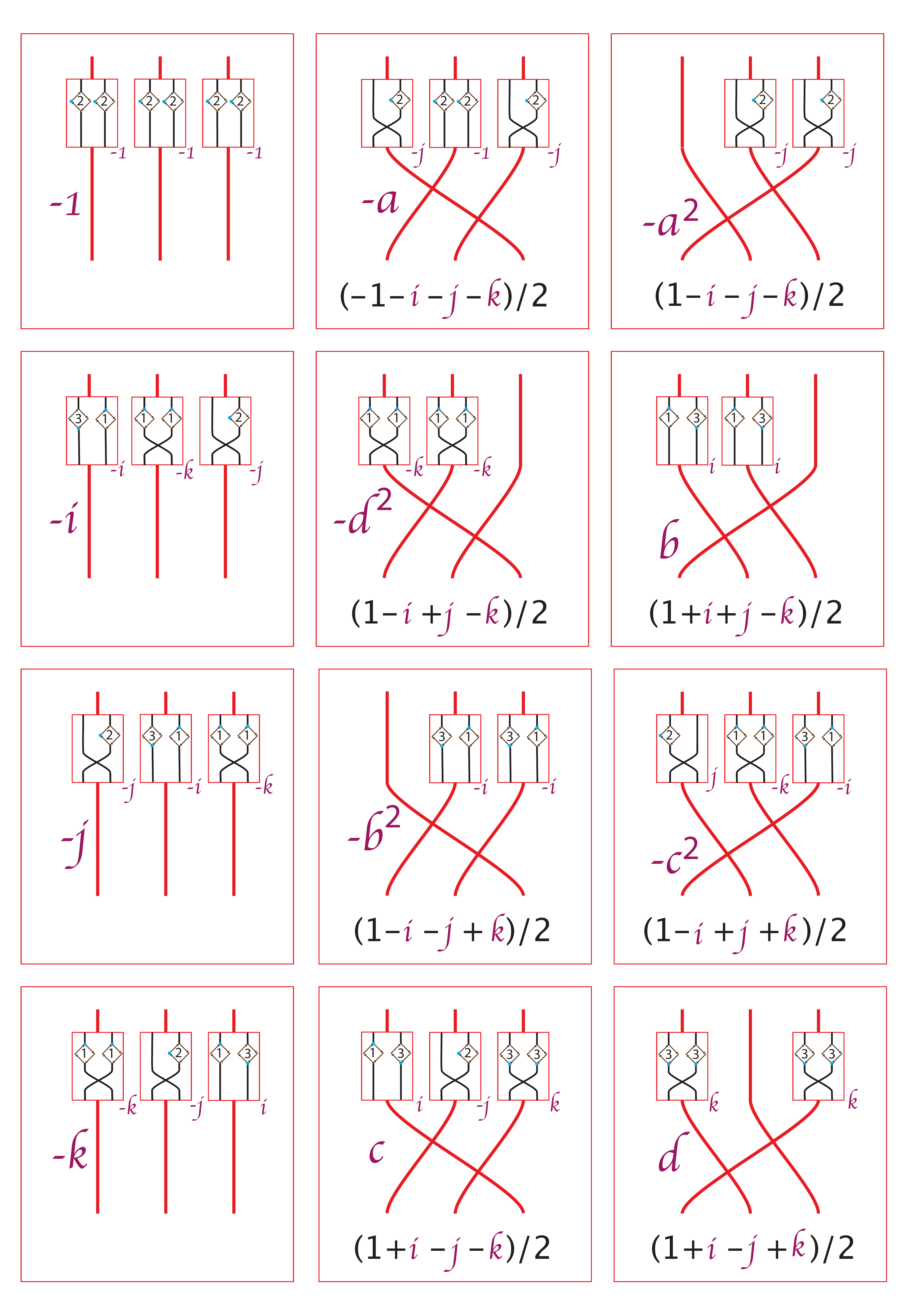}

Figure 31. The binary tetrahedral group  as it acts on cosets of $Q_8$, part 2
\end{center}

\subsection{The actions of $a$ and $b$ upon the cyclic subgroup generated by $a$}
\label{SS:BTfin}

The subgroup $N=(1,a^2,a^4)$   and coset $aN=(a,-1,a^5)$, when interwoven as an ordered coset form $A=(1,a,a^2,-1,-a,-a^2)$ --- a cyclic subgroup of $\widetilde{A_4}$. This will lead to a string-with-beads representation of the binary tetrahedral group $\widetilde{A_4}$ in which there are four strings that correspond to the ordered set of ordered cosets $(A, {\boldsymbol i}A, 
{\boldsymbol j}A, {\boldsymbol k}A)$. Two beads along a single string compose as  elements of $\Z/6$.

\begin{center}
\includegraphics[width=0.55\paperwidth]{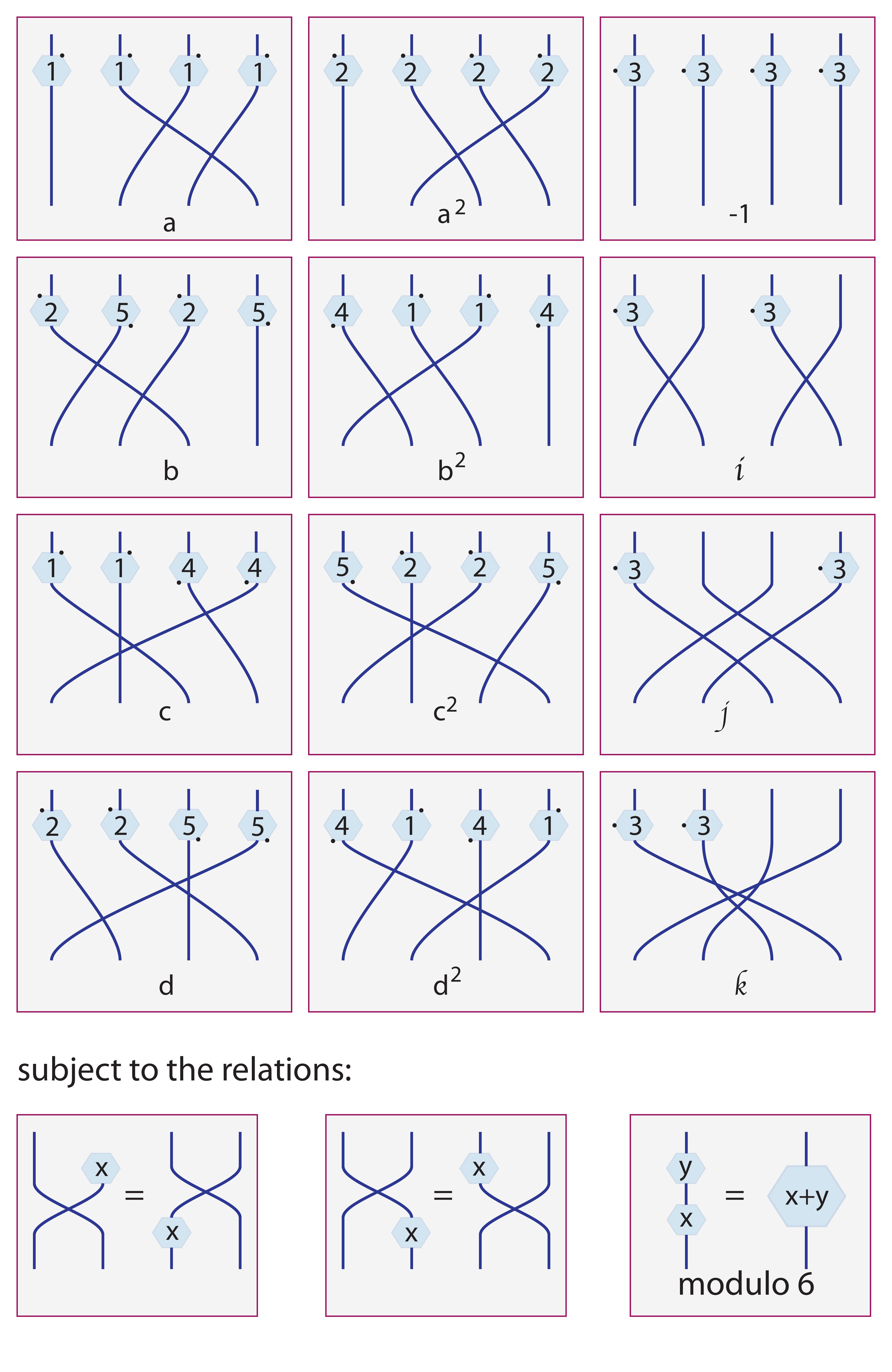}

Figure 32. A different 4 string representation
\end{center}

Explicitly, the cosets are ordered as follows: 
\begin{eqnarray*} 
A &=& (1,a,a^2,-1,-a,-a^2)\\
{\boldsymbol i}A &=& ({\boldsymbol i}, d^2, -b,-{\boldsymbol i}, -d^2, b ) \\ 
{\boldsymbol j}A &=& ({\boldsymbol j}, b^2, c^2,-{\boldsymbol j}, -b^2, -c^2 ) \\
{\boldsymbol k}A &=& ({\boldsymbol k}, -c, -d,-{\boldsymbol k}, c,d ).\end{eqnarray*}

The action of $a$ upon these ordered sets is as follows:

\begin{eqnarray*} 
aA &=& (a,a^2,-1,-a,-a^2,1) = [1/6,A]\\
a({\boldsymbol i}A) &=& (b^2,c^2,-{\boldsymbol j}, -b^2,-c^2,{\boldsymbol j})=[1/6,{\boldsymbol j}A] \\ 
a({\boldsymbol j}A) &=& (-c,-d,-{\boldsymbol k}, c,d,{\boldsymbol k})= [1/6,{\boldsymbol k}A] \\
a({\boldsymbol k}A) &=& (d^2,-b,-{\boldsymbol i}, -d^2,b,{\boldsymbol i} )= [1/6,{\boldsymbol i}A].\end{eqnarray*}

And here is the action of $b$:
\begin{eqnarray*} 
bA &=& ({b,\boldsymbol i}, d^2, -b,-{\boldsymbol i}, -d^2) = [5/6,{\boldsymbol i}A]\\
b({\boldsymbol i}A) &=& ( c^2,-{\boldsymbol j}, -b^2, -c^2,{\boldsymbol j}, b^2)=[1/3,{\boldsymbol j}A] \\ 
b({\boldsymbol j}A) &=& (a^2,-1,-a,-a^2,1,a)= [1/3,A] \\
b({\boldsymbol k}A) &=& (d,{\boldsymbol k}, -c,-d,-{\boldsymbol k},c)= [5/6,{\boldsymbol k}A].\end{eqnarray*}

 An illustration of half the elements is given in Fig.~32. The remaining half of the elements are negatives of these.  This completes the proof of item (7) of Theorem~\ref{main}.

\newpage 

\section{${\mbox{\rm GL}}_2(\Z/3)$}
\label{S:GL}

A common mistake is one that we made in a preliminary version of this paper. The group of invertible $(2\times 2)$ matrices defined over $\Z/3$ (the group ${\mbox{\rm GL}}_2(\Z/3)$) is a $2$-fold extension of $\Sigma_4$, but {\it it is not isomorphic to the binary octahedral group.} In this section, the strings-with-beads representation of ${\mbox{\rm SL}}_2(\Z/3)$ that is presented in Figs.~22
and~23 is extended to ${\mbox{\rm GL}}_2(\Z/3)$.

As in subsection~\ref{SS:SL},  $(2\times2)$ matrices with entires in $\Z/3$ acting upon the eight non-zero vectors in $\Z/3 \times \Z/3$ are considered. We recall that 
\begin{center}
\begin{tabular}{||l||l||}\hline
$0 \leftrightarrow (1,-1)$ & $4 \leftrightarrow (-1,1)$ \\
$1 \leftrightarrow (0,-1)$  & $5 \leftrightarrow (0,1)$ \\
$2 \leftrightarrow (1,1)$ & $6 \leftrightarrow (-1,-1)$ \\
$3 \leftrightarrow (-1,0)$ & $7 \leftrightarrow (1,0)$ \\ \hline \end{tabular}
\end{center}
\vspace{.5cm}

The details of  the computation that are analogous to those in Fig.~21 are left to the reader who can obtain  the permutation representations for the matrices of determinants  $-1$ that follow.

\begin{eqnarray*}
   \left[\begin{array}{rr} -1 & 0 \\ 0 & 1 \end{array}\right] 
    \leftrightarrow  (0,6)(2,4)(3,7); & &
    \left[\begin{array}{rr} 1 & 0 \\ 0 & -1 \end{array}\right] 
       \leftrightarrow  (0,2)(1,5)(4,6);\\
 \left[\begin{array}{rr} 0 & 1 \\ 1 & 0 \end{array}\right]
 \leftrightarrow  (0,4)(1,3)(5,7); & &
  \left[\begin{array}{rr} 0 & -1 \\ -1 & 0 \end{array}\right]
   \leftrightarrow (1,7)(2,6)(3,5); \\
    \left[\begin{array}{rr} 1 & 0 \\ 1 & -1 \end{array}\right]
 \leftrightarrow  (1,5)(2,7)(3,6); & &
  \left[\begin{array}{rr} 1 &0 \\ -1 & -1 \end{array}\right]
   \leftrightarrow (0,7)(1,5)(3,4); \\
    \left[\begin{array}{rr} -1 & 0 \\ 1 & 1 \end{array}\right]
 \leftrightarrow  (0,3)(2,6)(4,7); & &
  \left[\begin{array}{rr} -1 &0 \\ -1 & 1 \end{array}\right]
   \leftrightarrow (0,4)(2,3)(6,7); \\
    \left[\begin{array}{rr} 1 & 1 \\ 1 & 0 \end{array}\right]
 \leftrightarrow  (0,5,7,2,4,1,3,6); & &
  \left[\begin{array}{rr} 1 &-1 \\ -1 & 0 \end{array}\right]
   \leftrightarrow (0,6,5,3,4,2,1,7); \\
    \left[\begin{array}{rr} -1 & 1 \\ 1 & 0 \end{array}\right]
 \leftrightarrow  (0,2,5,7,4,6,1,3); & &
  \left[\begin{array}{rr} -1 &-1 \\ -1 & 0 \end{array}\right]
   \leftrightarrow (0,1,7,6,4,5,3,2); \end{eqnarray*}
    
   \begin{eqnarray*}
    \left[\begin{array}{rr} 1 & 1 \\ 0 & -1 \end{array}\right]
 \leftrightarrow  (0,5)(1,4)(2,6); & &
  \left[\begin{array}{rr} 1 &-1 \\ 0 & -1 \end{array}\right]
   \leftrightarrow (0,4)(1,2)(5,6); \\
    \left[\begin{array}{rr} -1 & 1 \\ 0 & 1 \end{array}\right]
 \leftrightarrow  (1,6)(2,5)(3,7); & &
  \left[\begin{array}{rr} -1 &-1 \\ 0 & 1 \end{array}\right] 
   \leftrightarrow (0,1)(3,7)(4,5); \\
    \left[\begin{array}{rr} 0 & 1 \\ 1 & 1 \end{array}\right]
 \leftrightarrow  (0,3,1,6,4,7,5,2); & &
  \left[\begin{array}{rr} 0 &1 \\ 1 & -1 \end{array}\right]
   \leftrightarrow (0,6,3,1,4,2,7,5); \\
    \left[\begin{array}{rr} 0 & -1 \\ -1 & 1 \end{array}\right]
 \leftrightarrow  (0,2,3,5,4,6,7,1); & &
  \left[\begin{array}{rr} 0 &-1 \\ -1 & -1 \end{array}\right]
   \leftrightarrow (0,7,1,2,4,3,5,6); \\
       \left[\begin{array}{rr} 1 & 1 \\ -1 & 1 \end{array}\right] 
 \leftrightarrow  (0,5,2,3,4,1,6,7); & &
  \left[\begin{array}{rr} 1 &-1 \\ 1 & 1 \end{array}\right]
   \leftrightarrow (0,3,6,5,4,7,2,1); \\
    \left[\begin{array}{rr} -1 & 1 \\ -1 & -1 \end{array}\right]
 \leftrightarrow  (0,7,6,1,4,3,2,5); & &
  \left[\begin{array}{rr} -1 &-1 \\ 1 & -1 \end{array}\right]
   \leftrightarrow (0,1,2,7,4,5,6,3).
 \end{eqnarray*}
 
 The strings-with-beads representations are given in Fig.~33 and~34.

\begin{center}
\includegraphics[width=0.6\paperwidth]{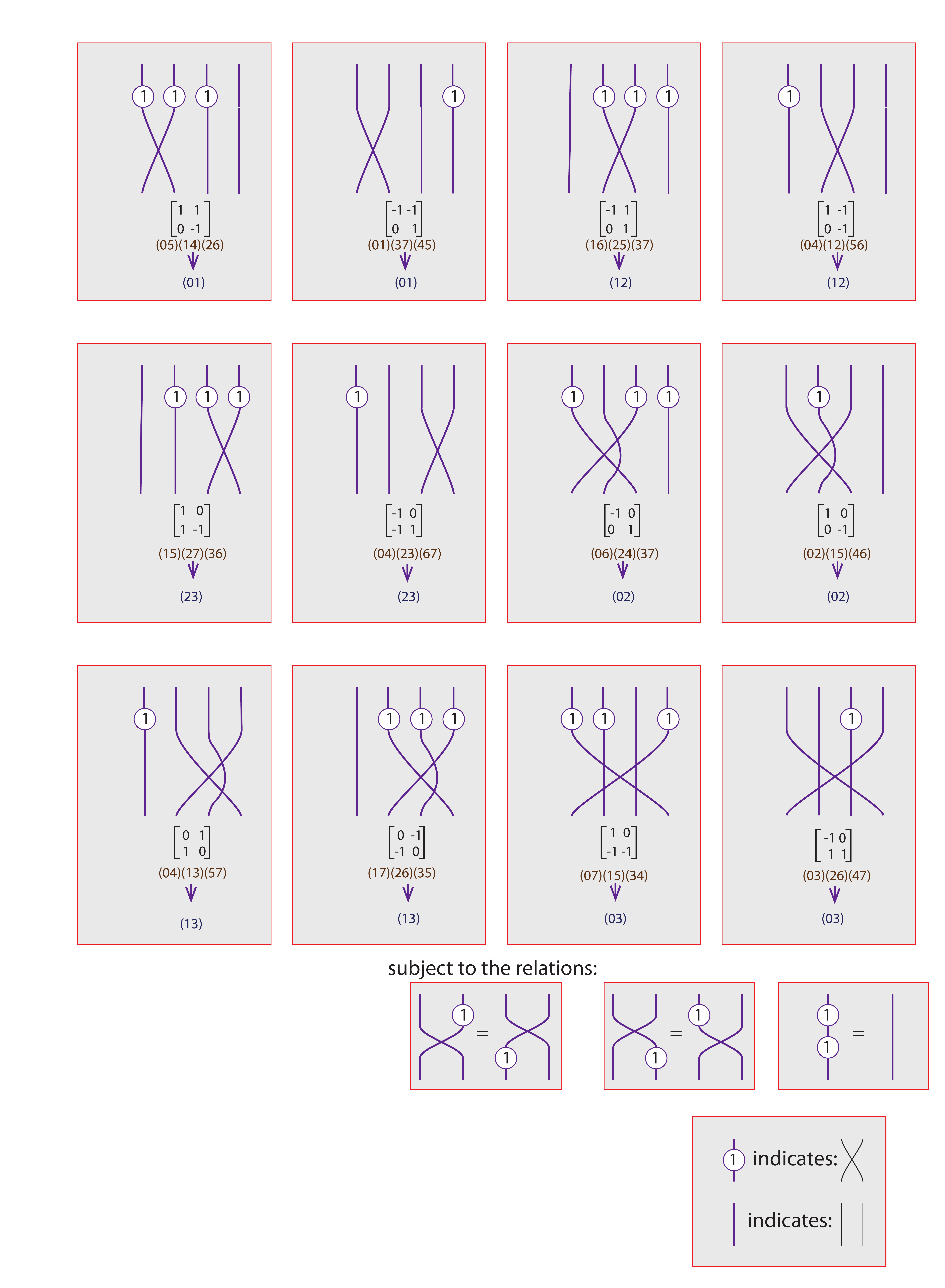}
Figure 33. 
Strings-with-beads representations of  ${\mbox{\rm GL}}_2(\Z/3)$, part 1
\end{center}

\begin{center}
\includegraphics[width=0.6\paperwidth]{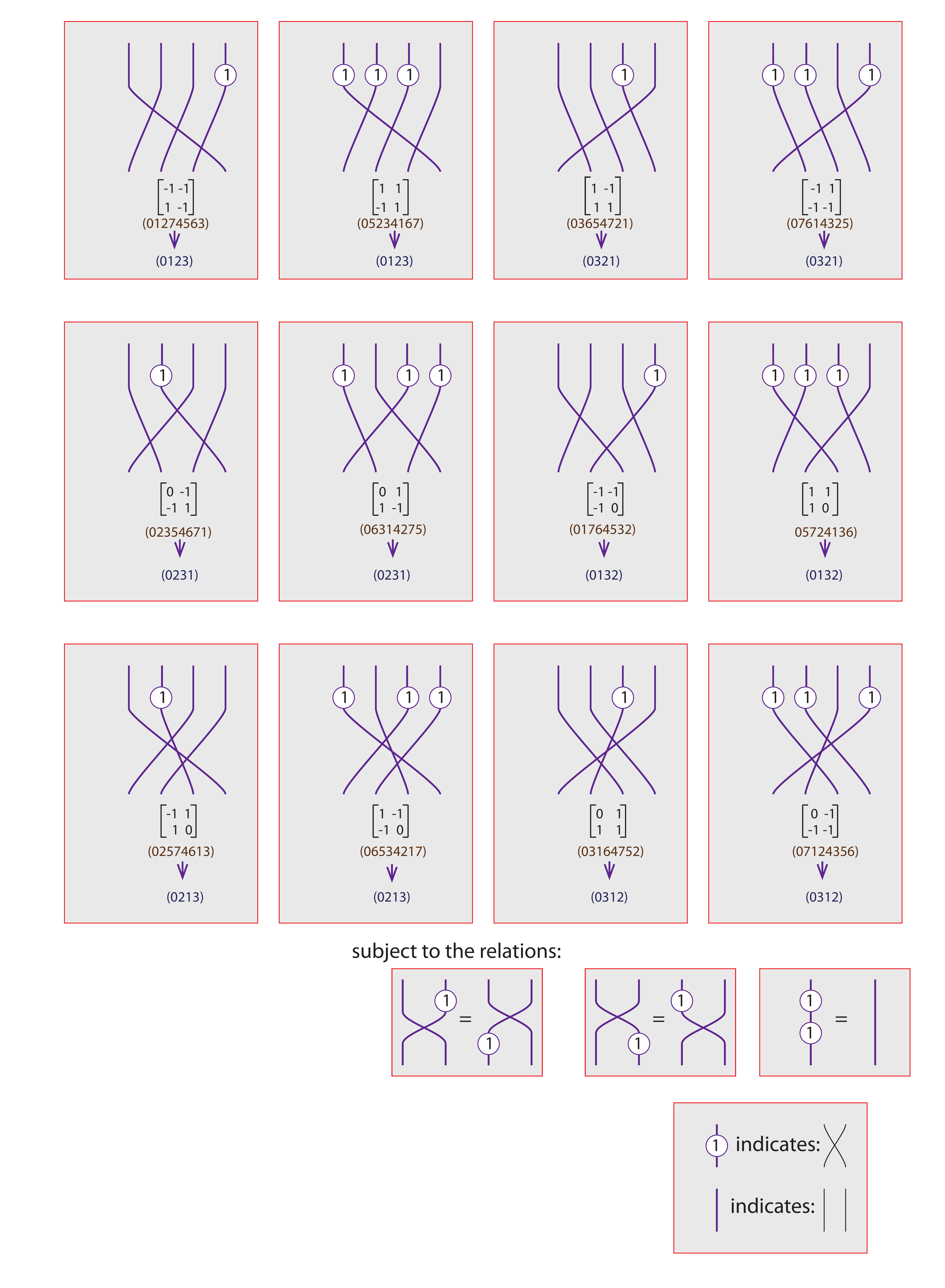}

Figure 34. Strings-with-beads representations of   ${\mbox{\rm GL}}_2(\Z/3)$,  part 2
\end{center}

This completes the proof of item (8) of Theorem~\ref{main}.

\newpage

 \begin{center}
\includegraphics[width=0.5\paperwidth]{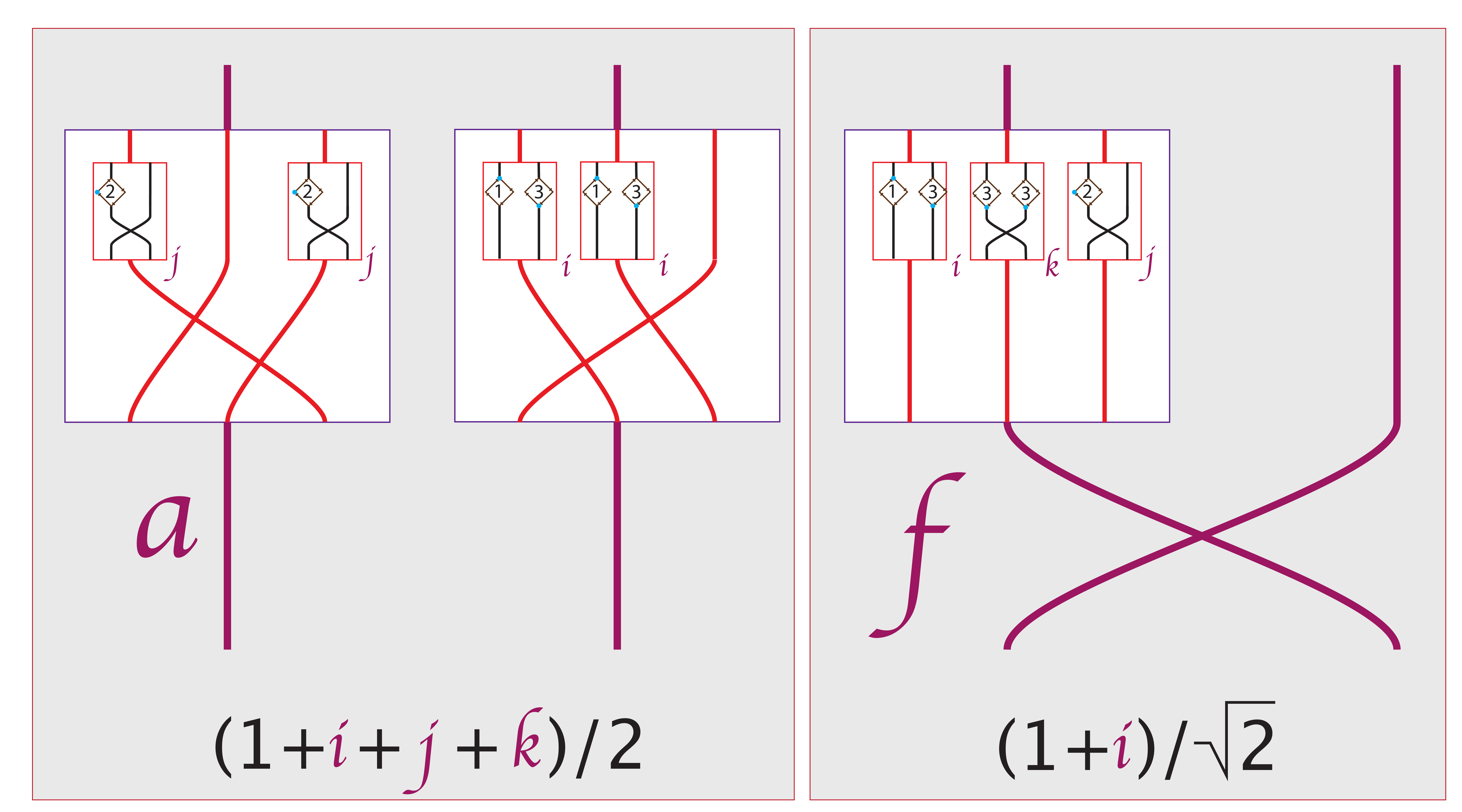}

Figure 35. A strings-with-beads representations of $a$ and $f$
\end{center}

\section{The binary octahedral group}
\label{S:BinOct}

The binary octahedral group, $\widetilde{\Sigma_4}$ is a $2$-fold extension of the permutation group $\Sigma_4$. It is given via the presentation 
$$\widetilde{\Sigma_4}= \langle a,f: a^3=f^4=(af)^2 \rangle.$$
The generators are identified with the elements $a=1/2 (1 + {\boldsymbol i} + {\boldsymbol j}+ {\boldsymbol k})$ and  $f=(1+{\boldsymbol i})/\sqrt{2}$ in $S^3$. See \cite{wiki:BinOct}. 
We give several representations that depend upon the subgroup chosen and the ordering of its cosets.

To begin, the actions of the  generators $a$ and $f$ on the (ordered) cosets of the subgroup $S=(-1,- {\boldsymbol i}, 1,  {\boldsymbol i})$ are computed. From these actions, strings-with-beads representations of $a$ and $f$ are obtained that have three strings and for which the beads are isomorphic to elements of the  ordered subgroup $Q_8=S \cup {\boldsymbol j}S$. The representations of $a$ and $f$ are presented in Fig.~35. Sample computations of products are presented in Figs.~36,~37, and~38. Then in anticipation of the next rerpresentation, half of the elements are presented in Figs.~39,~40,~41,~42,~43, and~44. The remaining twenty-four elements are negatives of these. 

The subgroup $P=(1,f,{\boldsymbol i}, f^3, -1, -f, -{\boldsymbol i}, -f^3)$ that consists of powers of $f$ and its coset ${\boldsymbol j}P$ form the dicyclic group ${\mbox{\rm Dic}}_4$ that has sixteen elements. The elements of this subgroup and those of the cosets $a(P\cup {\boldsymbol j}P)$ and $b(P\cup {\boldsymbol j}P)$ are represented in Figs.~46,~47,~48,~49,~50, and ~51. From these pictures the matrix representations, for example, 
\[ a \leftrightharpoons \left[ \begin{array}{ccc} 0 & 0 & {\boldsymbol j} \\ 1 & 0 & 0 \\ 0 & {\boldsymbol j} & 0  \end{array} \right] \quad {\mbox{\rm and}}
\quad b \leftrightharpoons 
\left[ \begin{array}{ccc} 0 &  {\boldsymbol i } & 0  \\  0 & 0 &  {\boldsymbol i }  \\ 1 & 0 &  0 \end{array} \right]; \] 
are extended to give the representation
\[f\leftrightharpoons 
\left[ \begin{array}{ccc} \frac{(1+{\boldsymbol i})}{\sqrt{2}}  &  0 & 0  \\  0 & 0  &  \frac{(1+{\boldsymbol i})}{\sqrt{2}}  \\ 0 & \frac{({\boldsymbol j}+{\boldsymbol k})}{\sqrt{2}}  &  0 \end{array} \right] \] 
 where the non-zero entries of the remaining $(3\times 3)$-matrices are  elements of the dicyclic group ${\mbox{\rm Dic}}_4$. 
 
 A final set of strings-with-beads representations of the binary octahedral group is obtained by observing that the subgroup $A=(1,a,a^2,-1,-a,-a^2)$ and its coset $\frac{{\boldsymbol i}-{\boldsymbol j}}{\sqrt{2}} A$ form a subgroup, $D= A \cup \frac{{\boldsymbol i}-{\boldsymbol j}}{\sqrt{2}} A $ that is isomorphic to the dicyclic group ${\mbox{\rm Dic}}_3$ of order twelve. There are four cosets, $D$, ${\boldsymbol i}D$, ${\boldsymbol j}D$, and ${\boldsymbol k}D$. So the strings-with-beads representation has four strings, and the beads are elements of $D$. Alternatively, the cosets of $A$ are considered. In this case, there are eight strings as was the case with the binary tetrahedral group $\widetilde{A_4}$, but the individual strings are bundles of four strings. They undergo quarter twists and sometimes cross pairwise. 
 
  \begin{center}
\includegraphics[width=0.5\paperwidth]{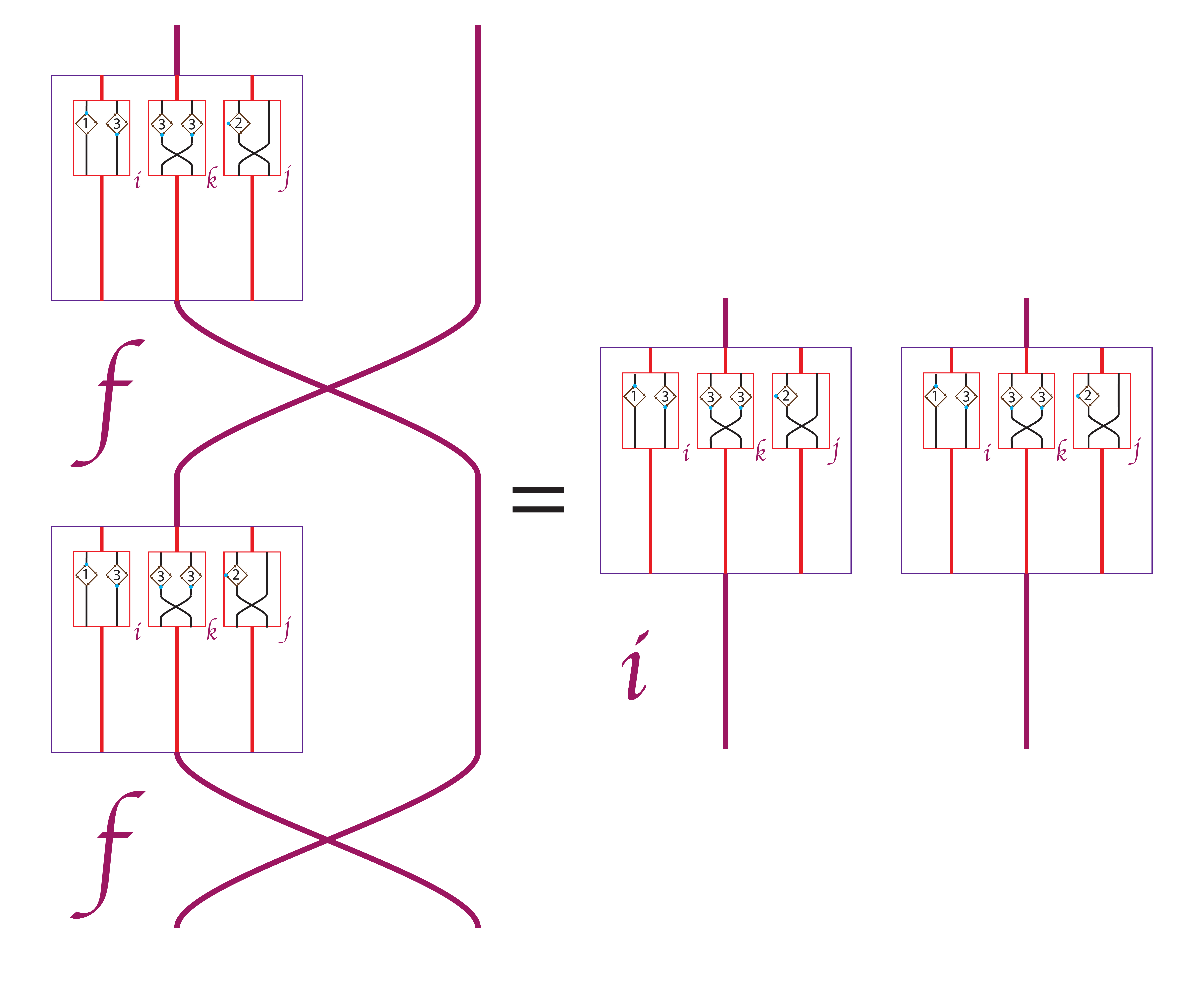}

Figure 36. The relationship $f^2= {\boldsymbol i}$ 
\end{center}

   \begin{center}
\includegraphics[width=0.5\paperwidth]{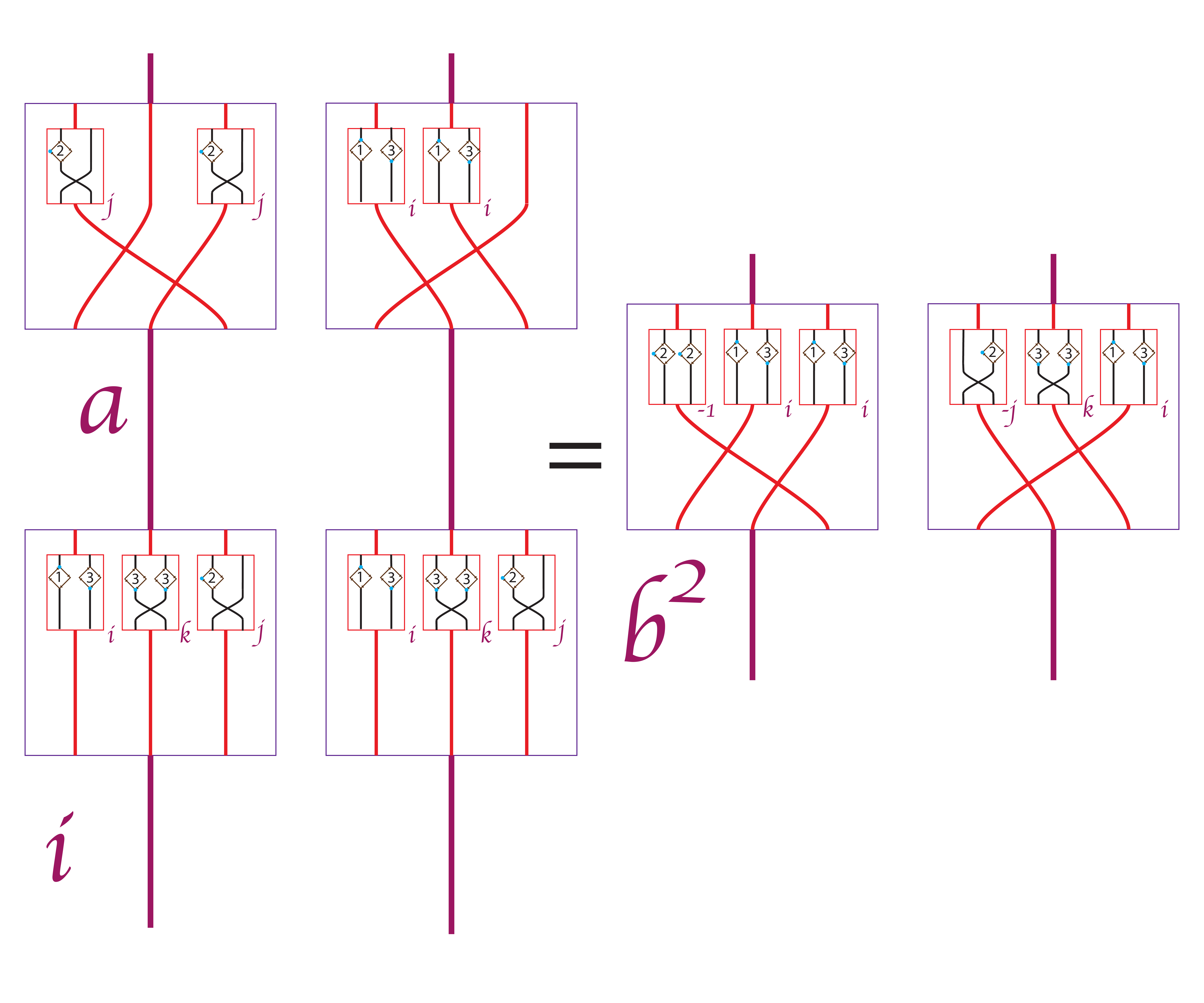}

Figure 37. The relationship $a {\boldsymbol i}=b^2$ 
\end{center}

\begin{center}
\includegraphics[width=0.5\paperwidth]{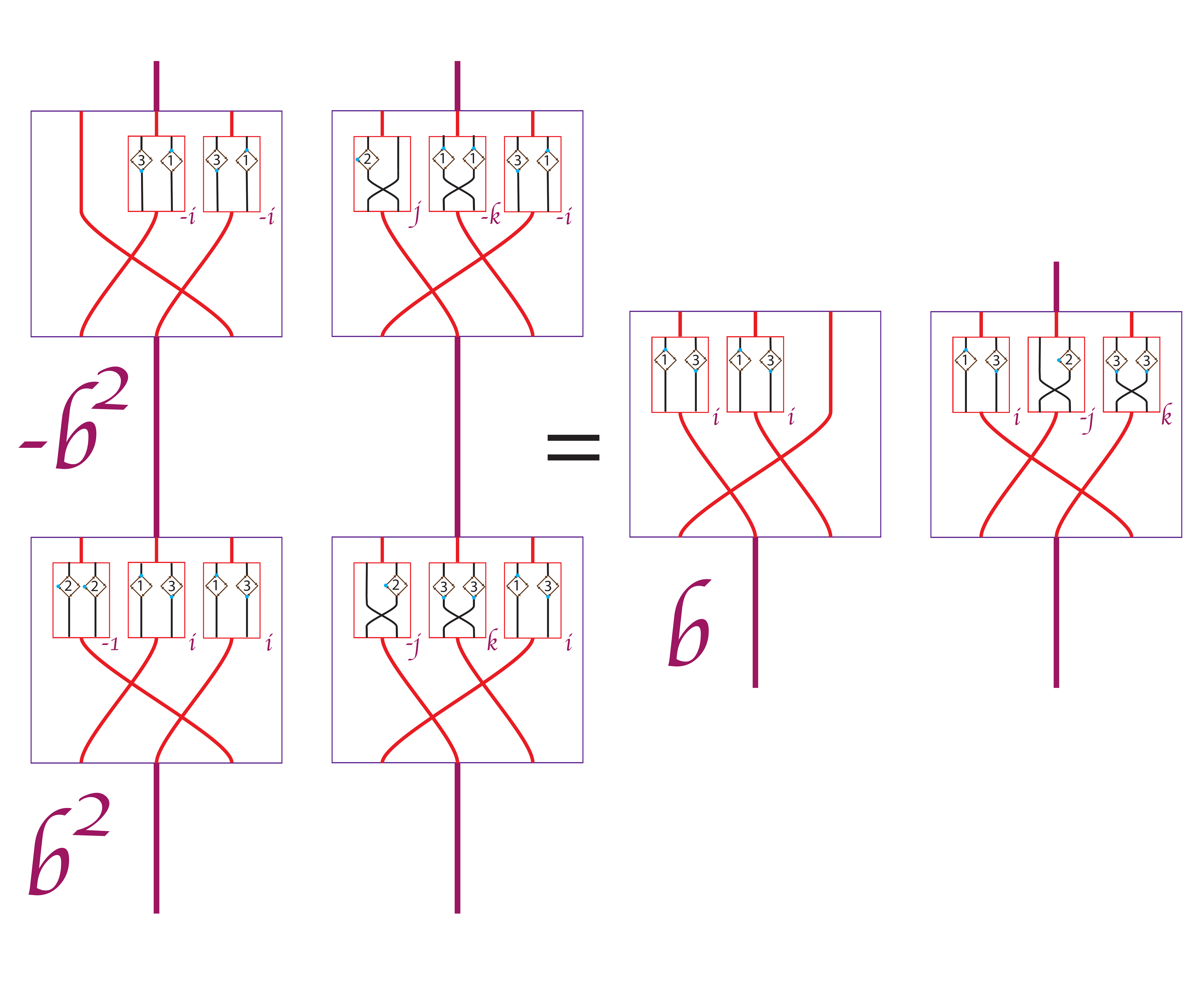}

Figure 38. The relationship $-b^4=b$ 
\end{center}

  \begin{center}
\includegraphics[width=0.6\paperwidth]{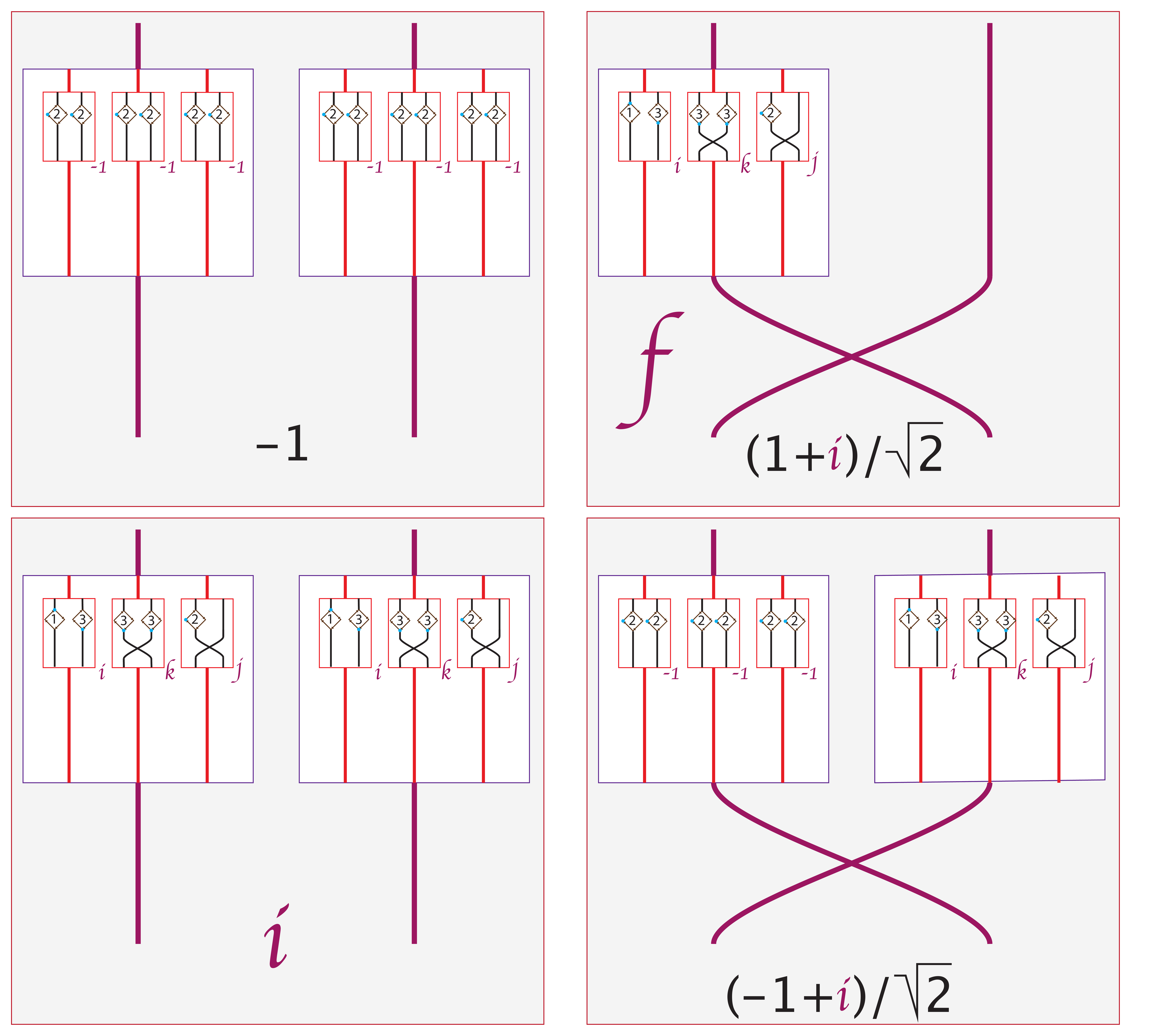}

Figure 39. The elements $-1$, $f= (1+{\boldsymbol i})/\sqrt{2}$, ${\boldsymbol i}$, and $f^3$ form half of the subgroup $P= (1,f,{\boldsymbol i},f^3, -1, -f, -{\boldsymbol i}, -f^3 ).$
\end{center}

     \begin{center}
\includegraphics[width=0.6\paperwidth]{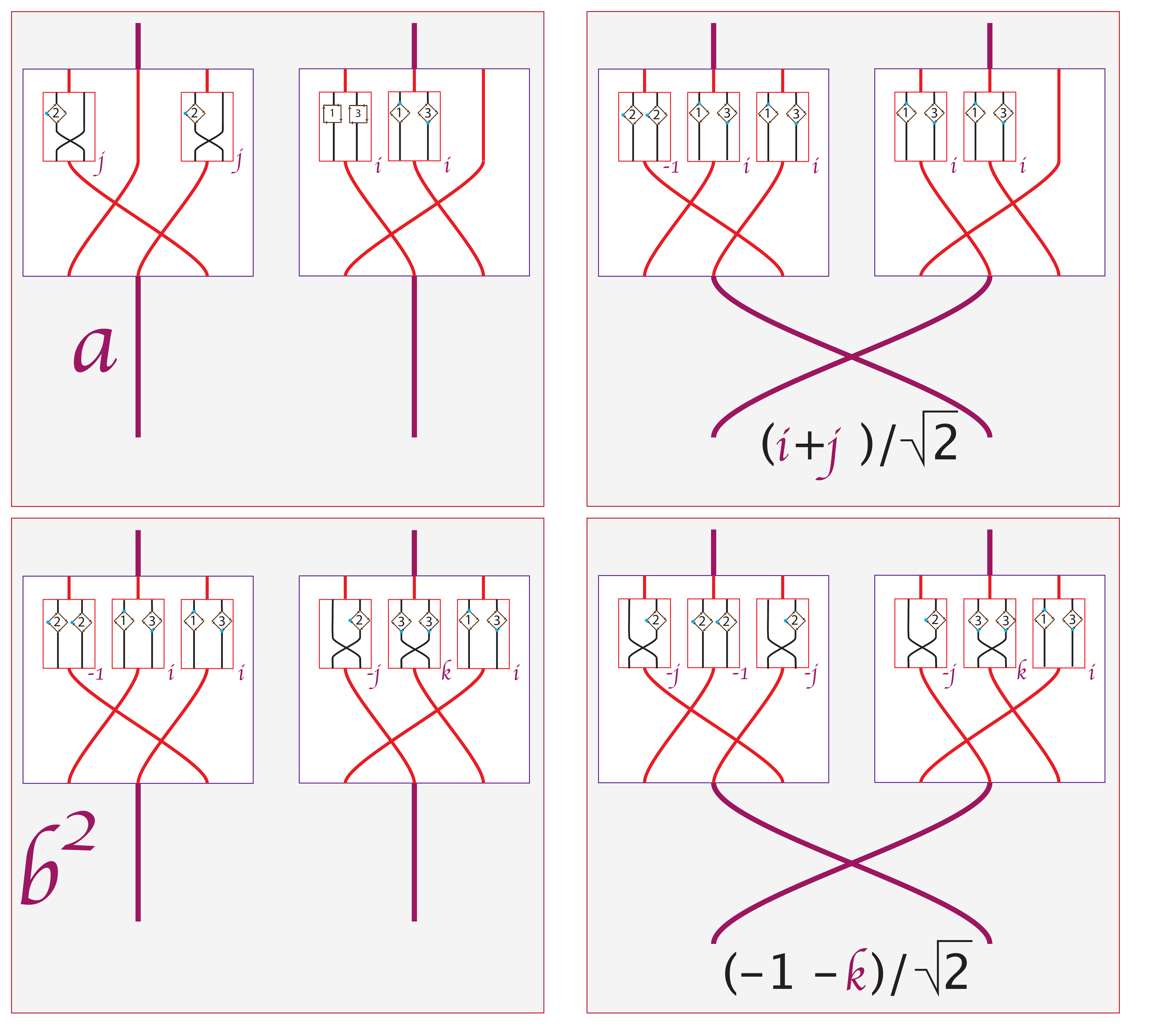}

Figure 40. Half of the elements in $aP$
\end{center}

\begin{center}
\includegraphics[width=0.6\paperwidth]{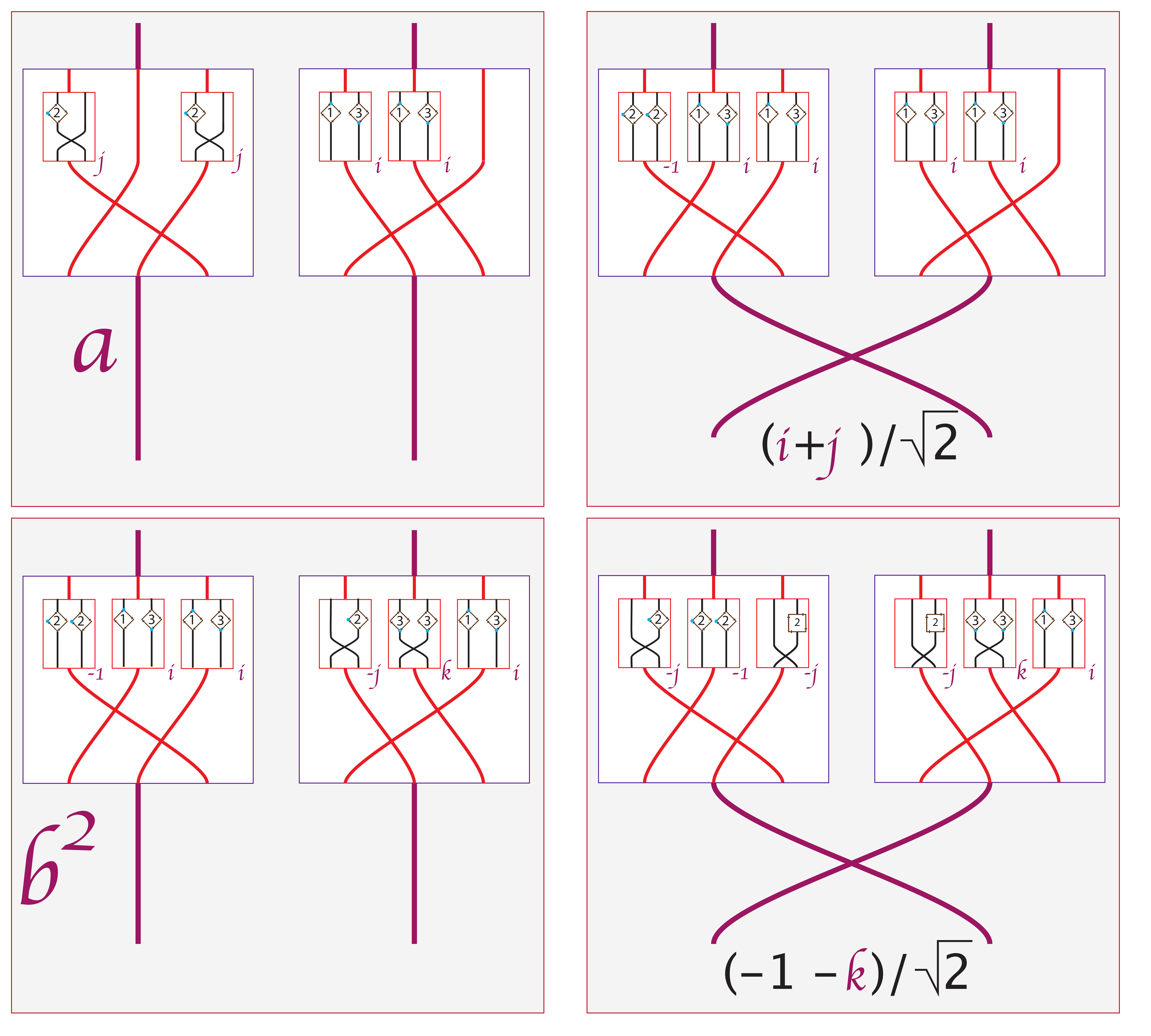}

Figure 41. Half of the elements in $bP$
\end{center}

      \begin{center}
\includegraphics[width=0.6\paperwidth]{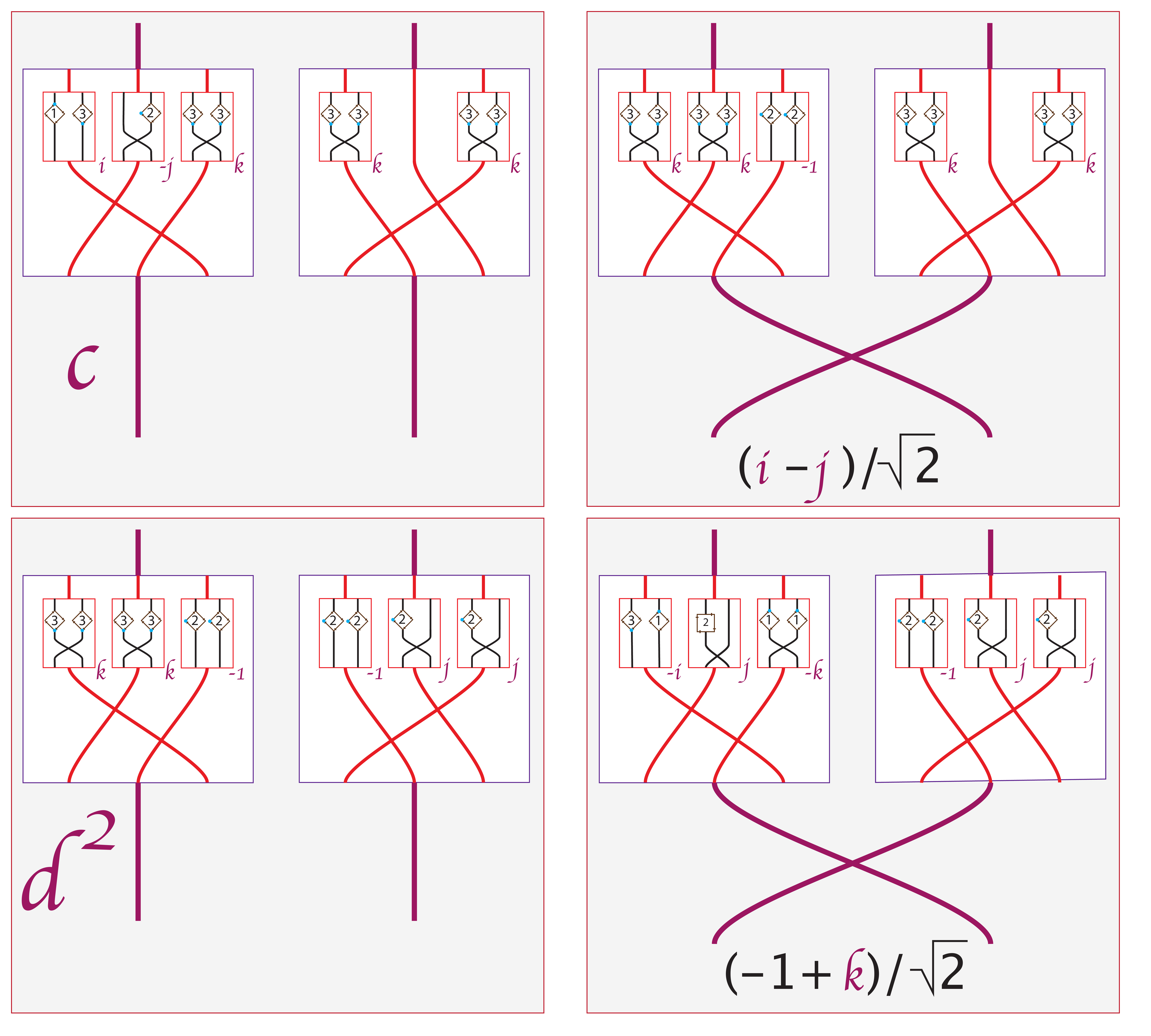}

Figure 42. Half of the elements in $cP$
\end{center}

 \begin{center}
\includegraphics[width=0.6\paperwidth]{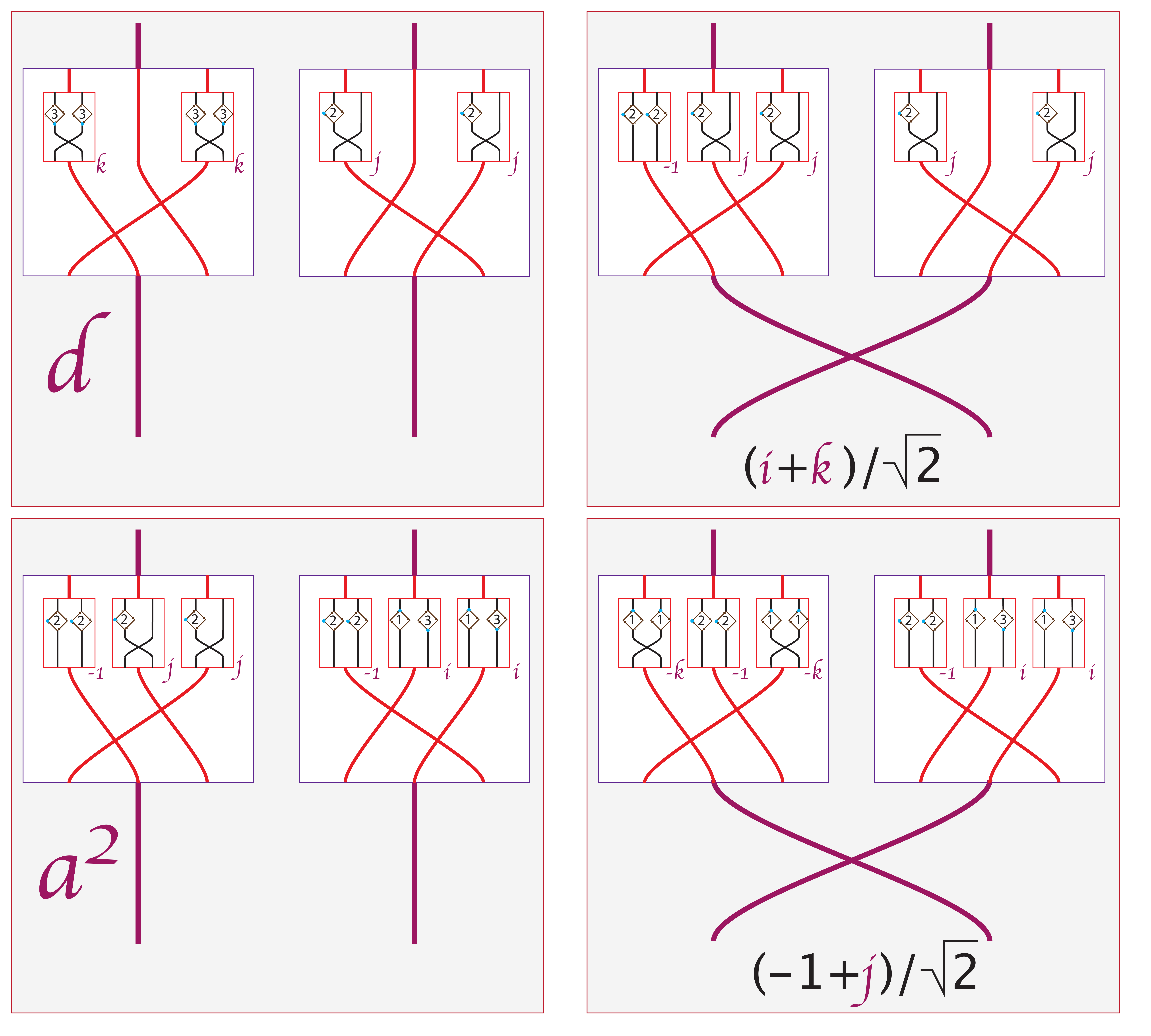}

Figure 43. Half of the elements in $dP$
\end{center}

   \begin{center}
\includegraphics[width=0.6\paperwidth]{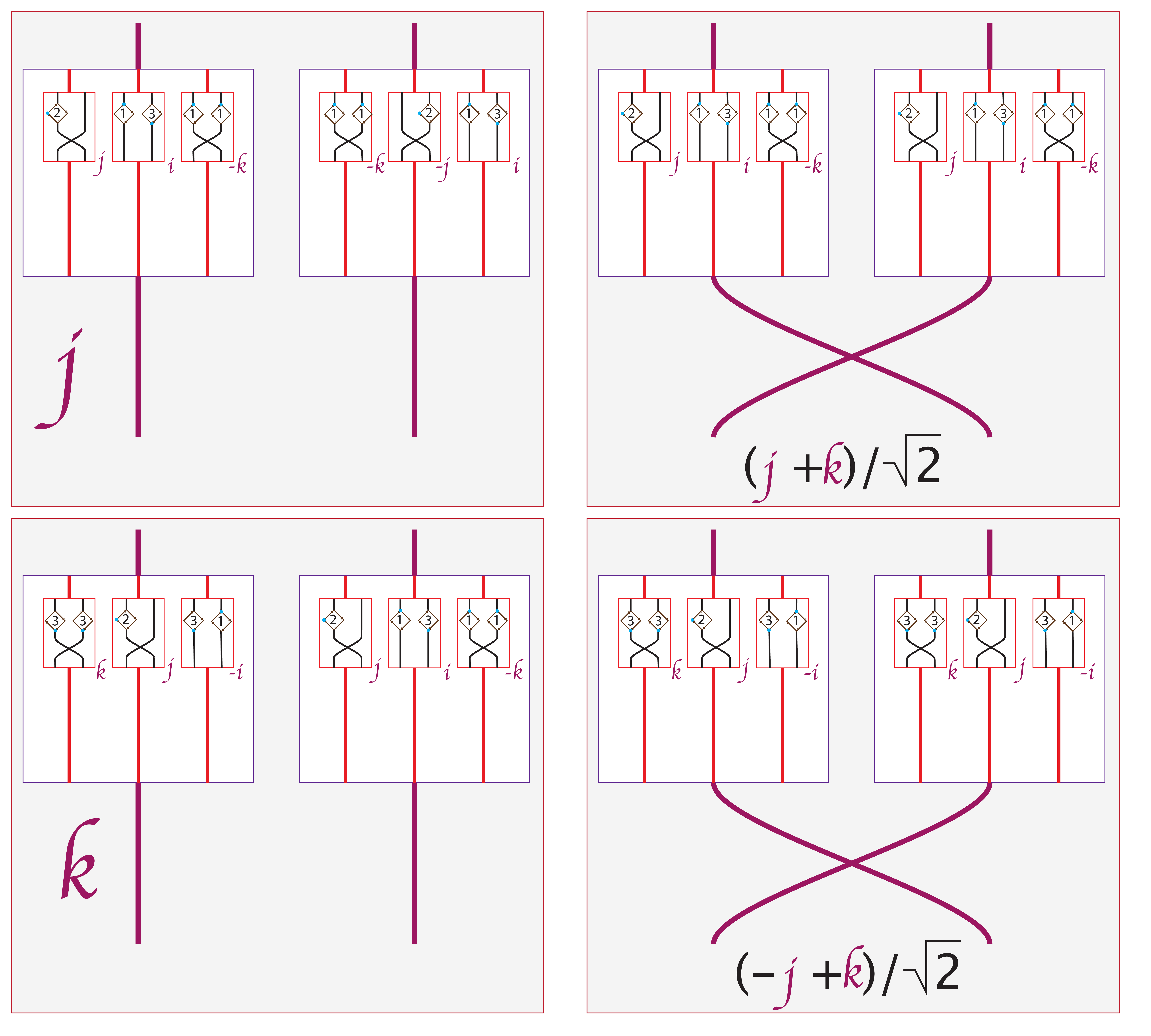}

Figure 44. Half of the elements in ${\boldsymbol j}P$
\end{center}

\subsection{Cosets of   $S=(-1,- {\boldsymbol i}, 1,  {\boldsymbol i})$}
\label{SS:poj}
 In the spirit of subsection~\ref{SS:BTAlt}, we decompose the binary octahedral group $\widetilde{\Sigma_4}$ into (cyclically) ordered cosets of the subgroup, $S=(-1,- {\boldsymbol i}, 1,  {\boldsymbol i})$. Of course the cosets ${\boldsymbol i}S$, $aS$, $a{\boldsymbol i}S$, $bS$ and $b{\boldsymbol i}S$ are as before. In $\widetilde{\Sigma_4}$, there are the six additional cosets: 
\[fS=
((-1-{\boldsymbol i}), (1-{\boldsymbol i}),  (1+{\boldsymbol i}), (-1+{\boldsymbol i}) )/\sqrt{2}, \]%
\[ f{\boldsymbol j}S=
((-{\boldsymbol j}-{\boldsymbol k}), (-{\boldsymbol j}+{\boldsymbol k}) ,({\boldsymbol j}+{\boldsymbol k}),({\boldsymbol j}-{\boldsymbol k}))/\sqrt{2},\] %
\[faS=((-{\boldsymbol i}-{\boldsymbol k}),(1-{\boldsymbol j}), ({\boldsymbol  i}+{\boldsymbol k}), (-1+{\boldsymbol j}) )/\sqrt{2}, \]
\[fa{\boldsymbol j}S
=(({\boldsymbol i}-{\boldsymbol k}) , (-1-{\boldsymbol j}),  (-{\boldsymbol i}+{\boldsymbol k}),(1+{\boldsymbol j}))
/\sqrt{2},\]
\[fbS=
((-{\boldsymbol i}-{\boldsymbol j}), (1+{\boldsymbol k}),   ({\boldsymbol  i}+{\boldsymbol j}), (-1-{\boldsymbol k}) )
/\sqrt{2}, \]
and
\[ fb{\boldsymbol j}S= ( (1-{\boldsymbol k}),  ({\boldsymbol i}-{\boldsymbol j}),  (-1+{\boldsymbol k}) ,(-{\boldsymbol i}+{\boldsymbol j}))/
\sqrt{2}.\]

 We write $x(yS)=[m/4,zS]$ where $x,y,z \in \widetilde{\Sigma_4}$, $m\in \Z/4$, and $S=(-1,-{\boldsymbol i},1,{\boldsymbol i})$ to indicate the equality of cosets $xyS= zS$, and that as cyclically ordered objects, the ordering is incremented by $m$. 
 For example, since  $f^2= {\boldsymbol i}$, we compute that $f(fS)=[1/4,S]$ and $f(f{\boldsymbol j})S=[3/4,{\boldsymbol j}S]$. 
 The action of $a$  on $\widetilde{A_4}= [S, {\boldsymbol j}S, aS, a{\boldsymbol j}S,  bS, b{\boldsymbol j}S]$ was computed earlier.
The  identities  $af =  ({\boldsymbol i}+ {\boldsymbol j})/\sqrt{2},$ $afa = (-1+ {\boldsymbol i})/\sqrt{2}$,  and $afb = (-1+ {\boldsymbol j})/\sqrt{2}$ are also very helpful and
lead to  the direct calculations below.
 \[\begin{array}{ll}
  a(fS)=[0/4,fbS], & 
  a(f{\boldsymbol j}S)=[0/4,fb{\boldsymbol j}S], \\
  a(faS)=[1/4,fS], & 
  a(fa{\boldsymbol j}S)=[3/4, f {\boldsymbol j}S], \\
  a(fbS)=[1/4,faS], &
  a(fb{\boldsymbol j}S)=[3/4,fa{\boldsymbol j}S], \\
  f(fS)=[1/4,S], & 
  f(f{\boldsymbol j}S)=[3/4, {\boldsymbol j}S], \\
  f(faS)=[3/4,a{\boldsymbol j}S], &
  f(fa{\boldsymbol j}S)=[3/4,aS] ,\\
  f(fbS)=[0/4,b{\boldsymbol j}S], &
  f(f(b{\boldsymbol i}S))=[2/4,bS]. \end{array}\]
  
The string-with-beads diagrams for $a$ and $f$, in which the beads are elements of the subgroup $Q_8$, are illustrated in Fig.~35. Sample computations appear in Figs.~36,~37, and~38. A compendium of half the elements (the remaining half are negatives of these) is presented in the Figs.~39,~40,~41,~42,~43, and~44. These diagrams are manifestations of the inclusion of normal subgroups
\[ 0 \rightarrow \Z/4 \rightarrow Q_8 \rightarrow \widetilde{A_4} \rightarrow \widetilde{\Sigma_4} \]
with quotient groups
\[ \widetilde{\Sigma_4}/\widetilde{A_4} \cong \Z/2; \quad  \widetilde{A_4}/Q_8 \cong \Z/3; \quad {\mbox{\rm and}} \  Q_8/(\Z/4) \cong \Z/2.\]

In order to present diagrams for the remaining elements of the binary octahedral group, many elements need to be written in terms of the generators. Fig.~36 indicates the identity $f^2 = {\boldsymbol i}$. Fig.~37  indicates that $a{\boldsymbol i} = b^2$, and Fig.~38 depicts the identity $-b^4=b$. The depictions of the identities $ab={\boldsymbol j}$, ${\boldsymbol i}{\boldsymbol j} = {\boldsymbol k}$, $b {\boldsymbol k} = d,$ and $bd=c$ are left for the reader.  This completes teh proof of item (9) of Theorem~\ref{main}.

\subsection{Cosets of   the dicyclic group ${\mbox{\rm Dic}}_4$ of order $16$.}
The union 
$S\cup {\boldsymbol j}S \cup fS \cup f{\boldsymbol j}S$ is the dicyclic group ${\mbox{\rm Dic}}_4$. See for example, Fig.~17. So  the elements of $S$ and $fS$ are interwoven and reordered into the ordered subgroup $P=(1,f,{\boldsymbol i}, f^3, -1, -f, -{\boldsymbol i}, -f^3)$ as are the elements of ${\boldsymbol j}S$ and $f{\boldsymbol j}S$ into the coset ${\boldsymbol j}P=  ({\boldsymbol j},{\boldsymbol j}f,{-\boldsymbol k}, {\boldsymbol j}f^3, -{\boldsymbol j}, -{\boldsymbol j}f, {\boldsymbol k}, -{\boldsymbol j}f^3)$.

Consider the ordered cosets $[(P\cup {\boldsymbol j}P), a(P\cup {\boldsymbol j}P) , b(P\cup {\boldsymbol j}P)]$ and the actions of $a$ and $f$ upon them. Write, for example, $f(P)=[1/8,P]$ and $f( {\boldsymbol j}P)=[7/8,  {\boldsymbol j}P]$ to indicate that $f$ rotates the cyclically ordered set $P$ through an angle of $45$ degrees, and it rotates $ {\boldsymbol j}P$ through an angle of $315$ degrees. The resulting computations appear below:

\begin{eqnarray*}
P & = & \left(1, 
\frac{1 + {\boldsymbol i}}{\sqrt{2}}, 
{\boldsymbol i},
\frac{-1+ {\boldsymbol i}}{\sqrt{2}}, 
-1, 
\frac{-1 - {\boldsymbol i}}{\sqrt{2}}, 
-{\boldsymbol i},
\frac{1 - {\boldsymbol i}}{\sqrt{2}}\right); \\
{\boldsymbol j}P &=&
\left({\boldsymbol j}, 
 \frac{{\boldsymbol j}-{\boldsymbol k}}
 {\sqrt{2}}, -{\boldsymbol k}, 
 \frac{-{\boldsymbol j}-{\boldsymbol k}}{\sqrt{2}},
 -{\boldsymbol j},
\frac{-{\boldsymbol j}+{\boldsymbol k}}{\sqrt{2}}, {\boldsymbol k}, \frac{{\boldsymbol j}+{\boldsymbol k})}{\sqrt{2}};\right)\\
aP & = & \left(a, 
\frac{{\boldsymbol i}+ {\boldsymbol j}}
{\sqrt{2}}, 
b^2,
\frac{-1- {\boldsymbol k}}{\sqrt{2}}, 
-a, 
\frac{-{\boldsymbol i} - {\boldsymbol j}}{\sqrt{2}}, 
-b^2,
\frac{1 + {\boldsymbol k}}{\sqrt{2}}\right);\\
a{\boldsymbol j}P & = & \left(
-c, 
\frac{-{\boldsymbol i} +{\boldsymbol j}}{\sqrt{2}},  
 -d^2,
\frac{1 - {\boldsymbol k}}{\sqrt{2}}, 
c, 
\frac{{\boldsymbol i} - {\boldsymbol j}}{\sqrt{2}},
d^2,
\frac{-1+ {\boldsymbol k}}{\sqrt{2}} \right); \end{eqnarray*} \begin{eqnarray*}
bP & = & \left(b, 
\frac{{\boldsymbol i} - {\boldsymbol k}}{\sqrt{2}}, 
c^2,\frac{-1- {\boldsymbol j}}{\sqrt{2}}, 
-b, 
\frac{-{\boldsymbol i} +{\boldsymbol k}}{\sqrt{2}}, 
-c^2,\frac{1 + {\boldsymbol j}}{\sqrt{2}}\right);\\
b{\boldsymbol j}P 
 & = & \left(
 a^2,
 \frac{-1+{\boldsymbol j}}{\sqrt{2}}, 
 -d, 
\frac{-{\boldsymbol i} - {\boldsymbol k}}{\sqrt{2}}, 
-a^2,
\frac{1 - {\boldsymbol j}}{\sqrt{2}},
d, 
\frac{{\boldsymbol i} +{\boldsymbol k}}{\sqrt{2}}, 
\right);
\end{eqnarray*}

\begin{eqnarray*} a(P) = [0/8, aP] ; & \rule{2ex}{0ex} & a( {\boldsymbol j}P)  = [0/8,  a{\boldsymbol j}P]; \\
a(aP) = [0/8, b{\boldsymbol j}P];  & \rule{2ex}{0ex} &  a(a  {\boldsymbol j}P)  =  [4/8,  b P]; \\
a( b P) = [0/8,{\boldsymbol j}P]; & \rule{2ex}{0ex} &  a(b  {\boldsymbol j}P) = [4/8, P]; \\
 f(P) =[1/8, P] ; & \rule{2ex}{0ex} &  f( {\boldsymbol j}P)= [7/8,  {\boldsymbol j}P]; \\
f(aP) = [7/8, b  {\boldsymbol j}P]; & \rule{2ex}{0ex} &  f(a  {\boldsymbol j}P) = [5/8,  b P]; \\
f( b P) =[1/8,a P]; & \rule{2ex}{0ex} & f(b  {\boldsymbol j}P) = [7/8, a  {\boldsymbol j}P] .
\end{eqnarray*}

 \begin{center}
\includegraphics[width=0.65\paperwidth]{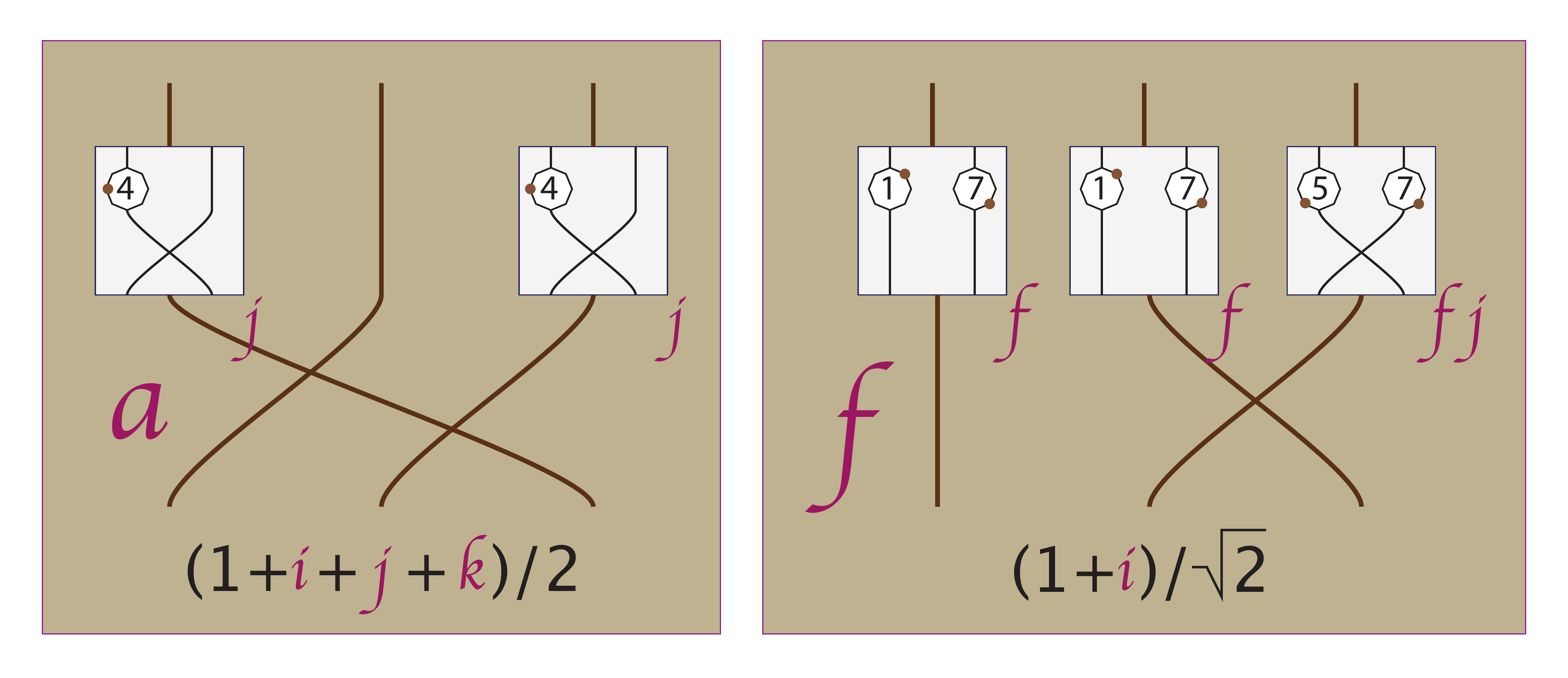}

Figure 45. The strings-with-beads representation of $a$ and $f$ with the cosets of $P\cup {\boldsymbol j}P$
\end{center}

 \begin{center}
 
\includegraphics[width=0.6\paperwidth]{CosetP.pdf}

Figure 46. The strings-with-beads representation of the powers of $f$ in $P$
\end{center}

\newpage

 \begin{center}
 
\includegraphics[width=0.6\paperwidth]{CosetjP.pdf}

Figure 47. The strings-with-beads representation of  the coset ${\boldsymbol j}P$
\end{center}

\newpage

\begin{center}
\includegraphics[width=0.6\paperwidth]{CosetaP.pdf}

Figure 48. The strings-with-beads representation of the  coset $aP$
\end{center}

\newpage

\begin{center}
\includegraphics[width=0.6\paperwidth]{CosetajP.pdf}

Figure 49. The strings-with-beads representation of the  coset $a{\boldsymbol j}P$
\end{center}

\newpage

\begin{center}
\includegraphics[width=0.6\paperwidth]{CosetbP.pdf}

Figure 50. The strings-with-beads representation of the  coset $bP$
\end{center}

\newpage

\begin{center}
\includegraphics[width=0.6\paperwidth]{CosetbJP.pdf}

Figure 51.The strings-with-beads representation of the  coset $b{\boldsymbol j}P$
\end{center}


\newpage

This gives the embedding $\widetilde{\Sigma_4} \subset \left[ {\mbox{\rm Dic}}_4\right]^3 \rtimes \Sigma_3 $ of Theorem~\ref{main} item (10).

The strings-with-beads representations for $a$ and $f$ are depicted in Fig.~45.  The elements of $P$, ${\boldsymbol j}P$, and those of the cosets $a(P\cup {\boldsymbol j}P)$ and $b(P\cup {\boldsymbol j}P)$ are represented in Figs.~46,~47,~48,~49,~50, and ~51.

\bigskip

\subsection{Cosets of $A=(1,a,a^2,-1,-a,-a^2)$}
The subgroup $A$ and its translate $xA$,
 where $x=\frac{({\boldsymbol i} - {\boldsymbol j})}{\sqrt{2}}$, will be considered. These are written below together with their cosets. They form a  form a subgroup of the binary octahedral group $\widetilde{\Sigma_4}$ that is isomorphic to a dicyclic group ${\mbox{\rm Dic}}_3$ of order $12$. Write $D= A\cup xA$ with the usual understanding that both $A$ and $xA$ are cyclically ordered as indicated above. The remaining ordered cosets are
 \[A=(1,a,a^2,-1,-a,-a^2),\]    
 \[ {\boldsymbol i}A = ({\boldsymbol i}, d^2, -b, -{\boldsymbol i}, -d^2,b), \] 
\[ {\boldsymbol j}A = ({\boldsymbol j}, b^2, c^2, -{\boldsymbol j}, -b^2,-c^2), \]
\[ {\boldsymbol k}A = ({\boldsymbol k}, -c, -d, -{\boldsymbol k}, c,d). \]
And
 \[ xA=\left(\frac{({\boldsymbol i} - {\boldsymbol j})}{\sqrt{2}}, 
\frac{(-{\boldsymbol j} +{\boldsymbol k})}{\sqrt{2}},
\frac{(-{\boldsymbol i} + {\boldsymbol k})}{\sqrt{2}},
\frac{(-{\boldsymbol i} + {\boldsymbol j})}{\sqrt{2}},
\frac{({\boldsymbol j} - {\boldsymbol k})}{\sqrt{2}},
\frac{({\boldsymbol i} - {\boldsymbol k})
}{\sqrt{2}}\right),\]
\[ 
x{\boldsymbol i}A = 
\left( \frac{-(1+ {\boldsymbol k})}{\sqrt{2}}, 
\frac{-({\boldsymbol j} +{\boldsymbol k})}{\sqrt{2}},
\frac{(1 - {\boldsymbol j})}{\sqrt{2}},
\frac{(1+ {\boldsymbol k})}{\sqrt{2}},
\frac{({\boldsymbol j}+{\boldsymbol k})}{\sqrt{2}},
\frac{(-1+{\boldsymbol j})}{\sqrt{2}} \right),
 \]

 \[
x{\boldsymbol j}A=
\left( \frac{ (1- {\boldsymbol j})}{\sqrt{2}}, 
\frac{(1 +{\boldsymbol i})
}{\sqrt{2}},
\frac{({\boldsymbol i} + {\boldsymbol k})}
{\sqrt{2}},
\frac{(-1 + {\boldsymbol k})}
{\sqrt{2}},
\frac{-(1+{\boldsymbol i})}
{\sqrt{2}},
\frac{-({\boldsymbol i} +{\boldsymbol k}))}{\sqrt{2}}\right),
 \] 
 \[
x{\boldsymbol k}A= 
\left( \frac{({\boldsymbol i} + {\boldsymbol j})}{\sqrt{2}},
\frac{(-1 + {\boldsymbol i})}{\sqrt{2}}, 
\frac{-(1+{\boldsymbol j})}{\sqrt{2}},
\frac{-({\boldsymbol i} + {\boldsymbol j})}{\sqrt{2}},
\frac{(1- {\boldsymbol i})}{\sqrt{2}},
\frac{(1+{\boldsymbol j}))}{\sqrt{2}} \right). \]

The actions of $f=(1+{\boldsymbol i})/\sqrt{2}$ and $a=(1+{\boldsymbol i} + {\boldsymbol j}+{\boldsymbol k})/2$ upon the  (ordered) cosets of the subgroup $A$ is as follows. 
We write $aA=[1,A]$ to indicate that $A$ has been rotated by $\pi/3$. 
\[ aA=[1,A], \quad  a{\boldsymbol i} A =[1,{\boldsymbol j} A], \quad  a{\boldsymbol j} A =[1,{\boldsymbol k} A] \quad a{\boldsymbol k} A =[1,{\boldsymbol i} A],\]
\[axA=[5,xA],  \quad  a{\boldsymbol i}x A =[5,{\boldsymbol j}x A], \quad a{\boldsymbol j}x A =[5,{\boldsymbol j}x A],  \quad a{\boldsymbol i}x A =[5,{\boldsymbol i}x A],\]
\[ fA= [1, {\boldsymbol j}xA], \quad f{\boldsymbol i}A= [1, {\boldsymbol k}xA], \quad f{\boldsymbol j}A= [4, {\boldsymbol i}xA], \quad f{\boldsymbol k}A= [1, xA],\]
\[ f(xA) = [2, {\boldsymbol j}A], 
 \quad  f(x{\boldsymbol i}A) = [2, {\boldsymbol k}A],  \quad f(x{\boldsymbol j}A) = [5, {\boldsymbol i}A], 
 \quad f(x{\boldsymbol k}A) = [1, A]. \]
The representation of these elements as strings-with-beads, where the beads are elements of the dicyclic group $A\cup xA$ (and $x=(i-j)\sqrt{2}$) is presented in Fig.~52. The elements in the various cosets encircle the central hexagons in Figs.~53 through~60. This completes the proof of Theorem~\ref{main}, item (11).

       \begin{center}
       
\includegraphics[width=0.6\paperwidth]{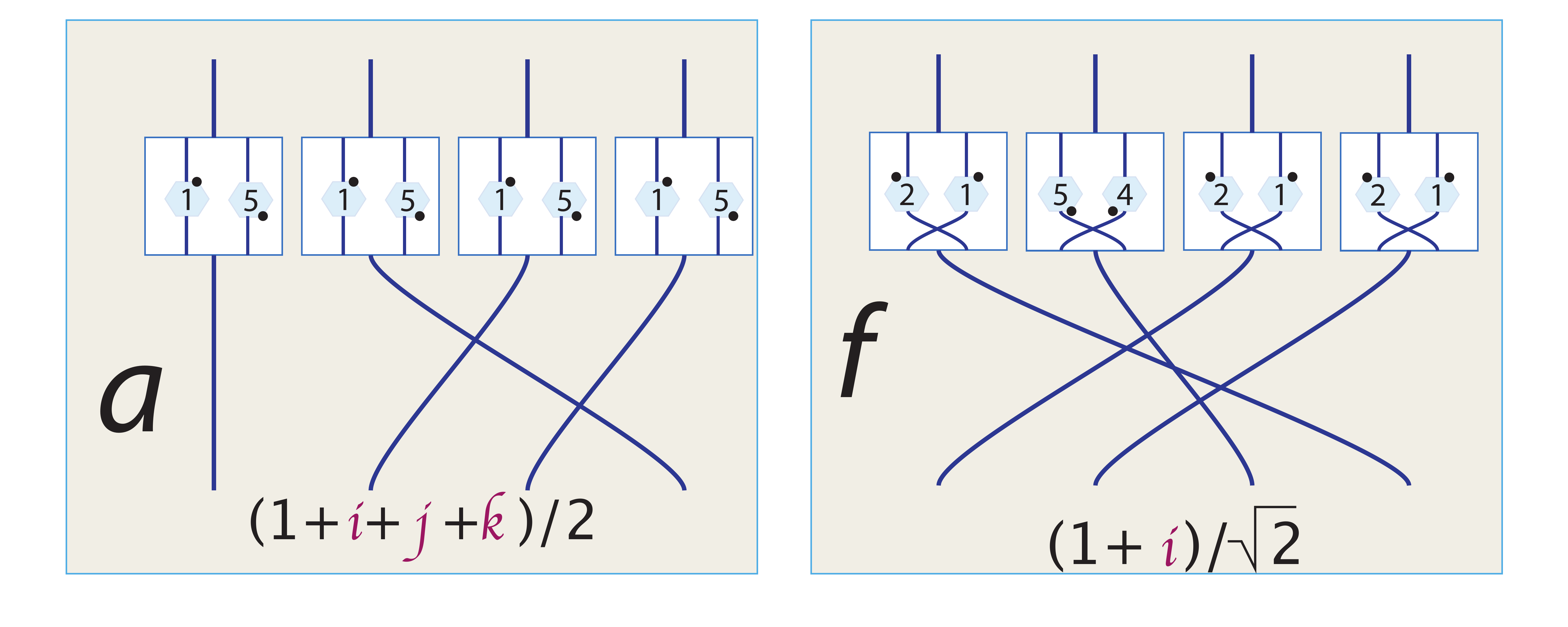}

Figure 52. The representations of $a$ and $f$ using the subgroup $A\cup xA$
\end{center}

\newpage

       \begin{center}
\includegraphics[width=0.65\paperwidth]{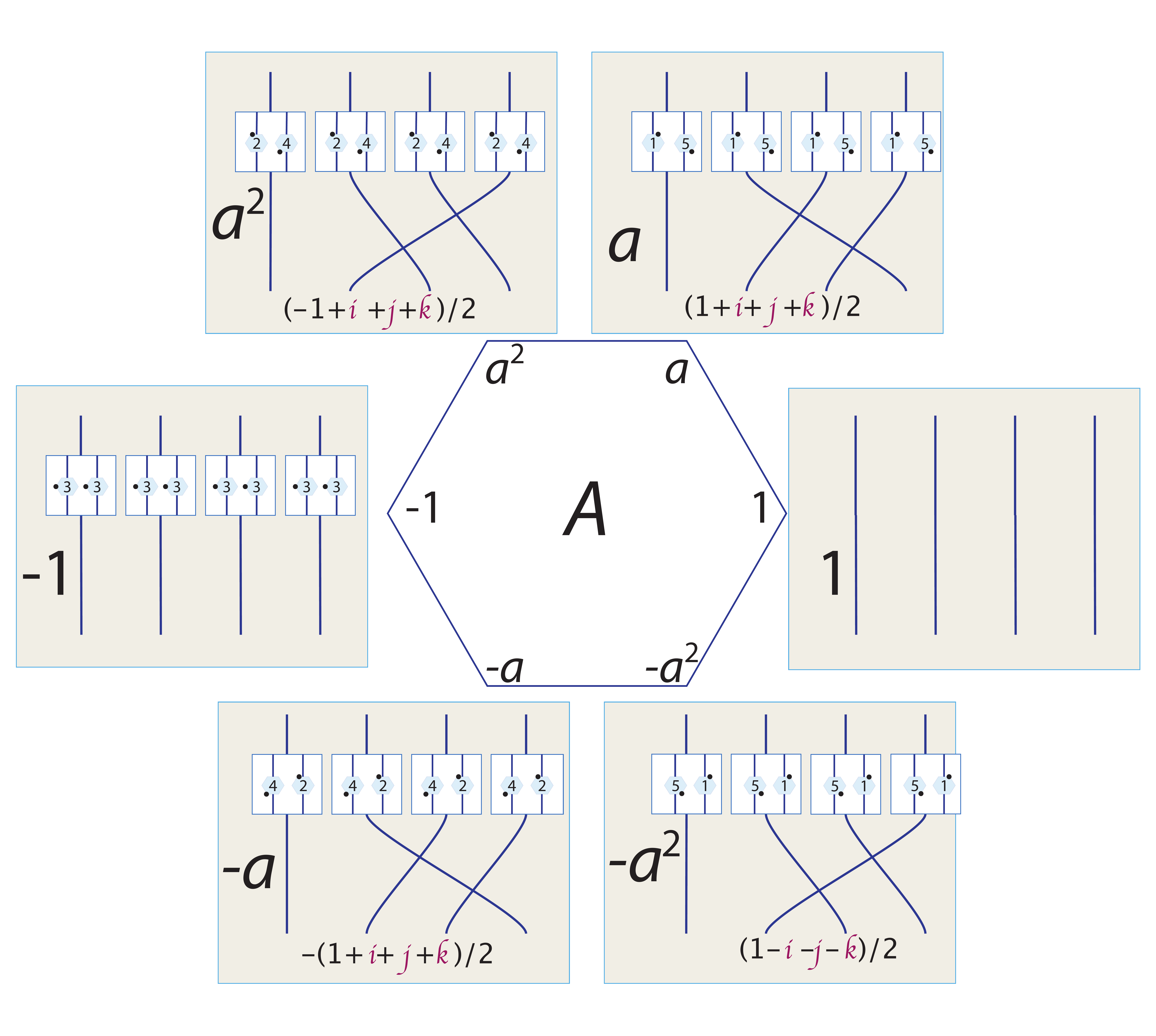}

Figure 53. The coset $A$
\end{center}

\newpage

       \begin{center}
\includegraphics[width=0.65\paperwidth]{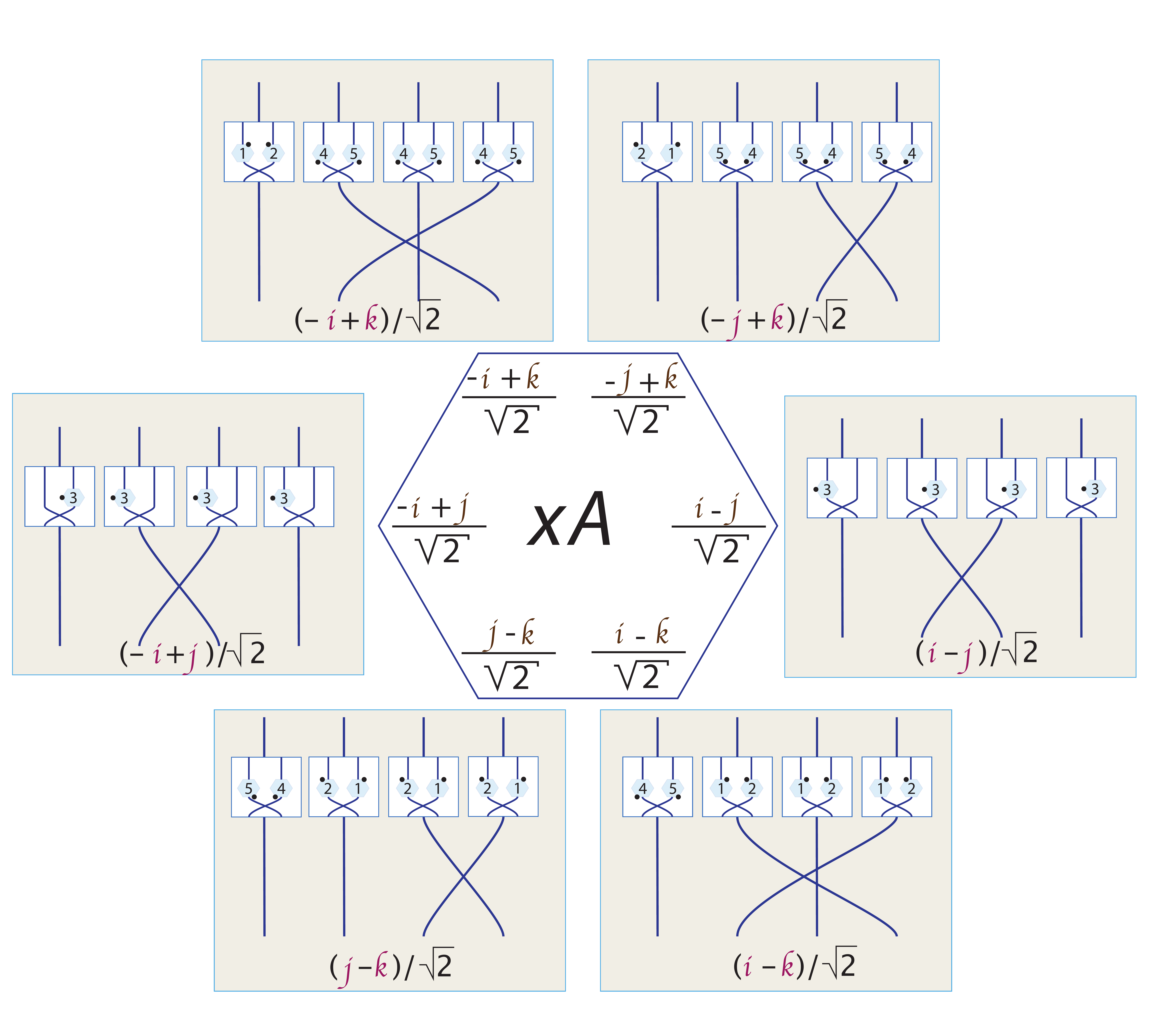}

Figure 54. The coset $xA$
\end{center}

\newpage 

     \begin{center}
    \includegraphics[width=0.65\paperwidth]{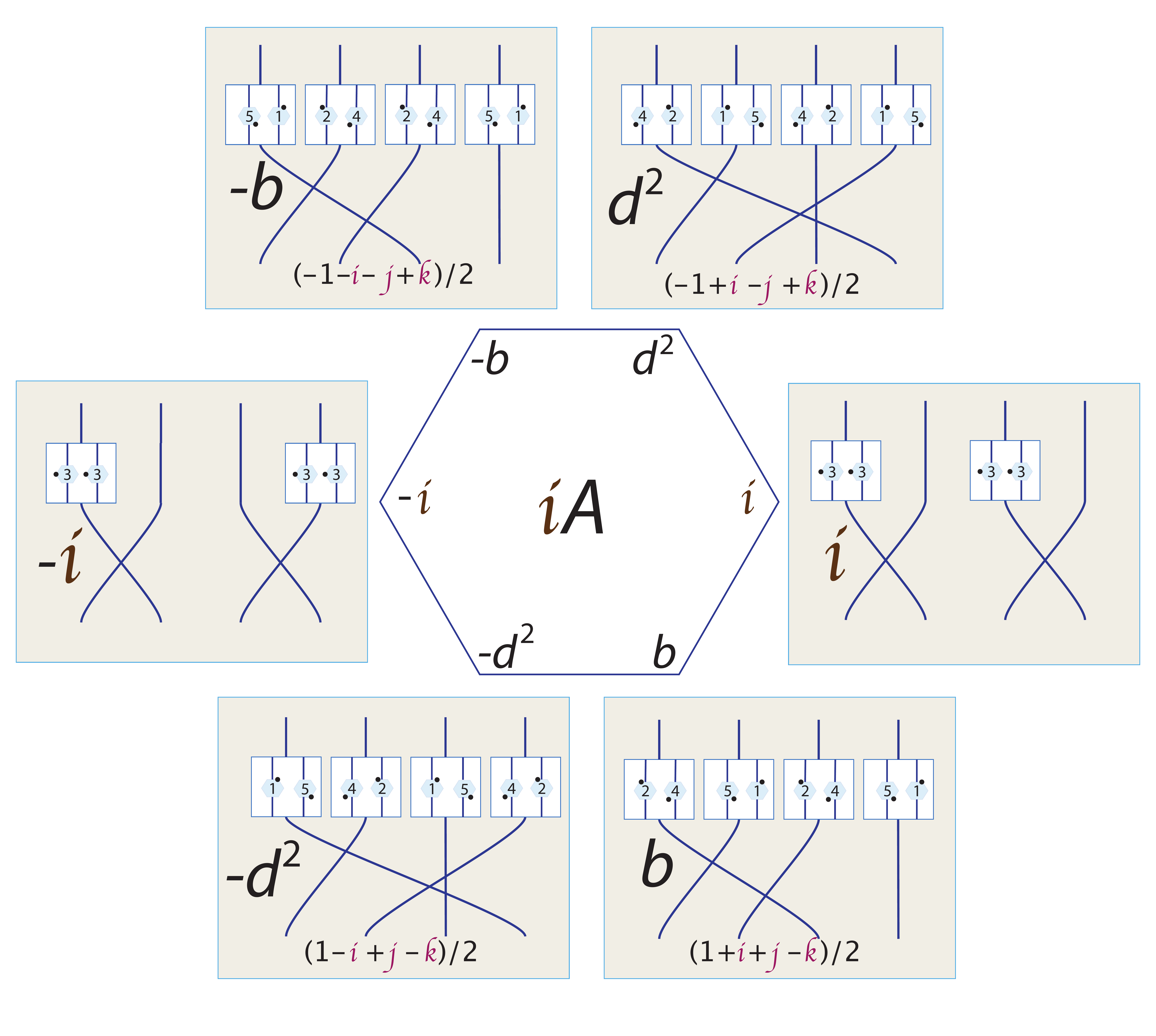}

Figure 55. The coset ${\boldsymbol i}A$
\end{center}

\newpage

       \begin{center}
\includegraphics[width=0.65\paperwidth]{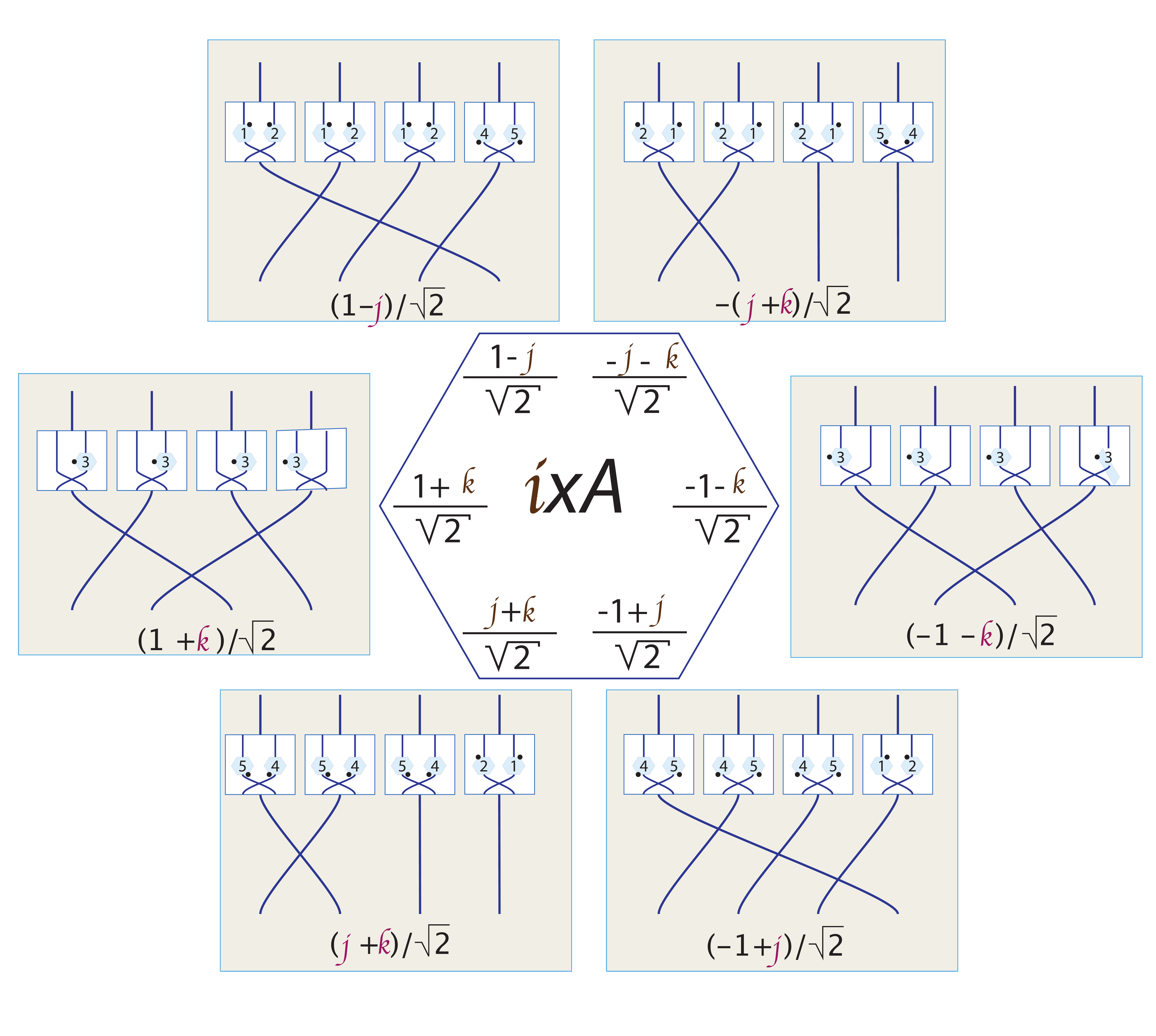}

Figure 56. The coset $x{\boldsymbol i}A$
\end{center}

\newpage

    \begin{center}
\includegraphics[width=0.5\paperwidth]{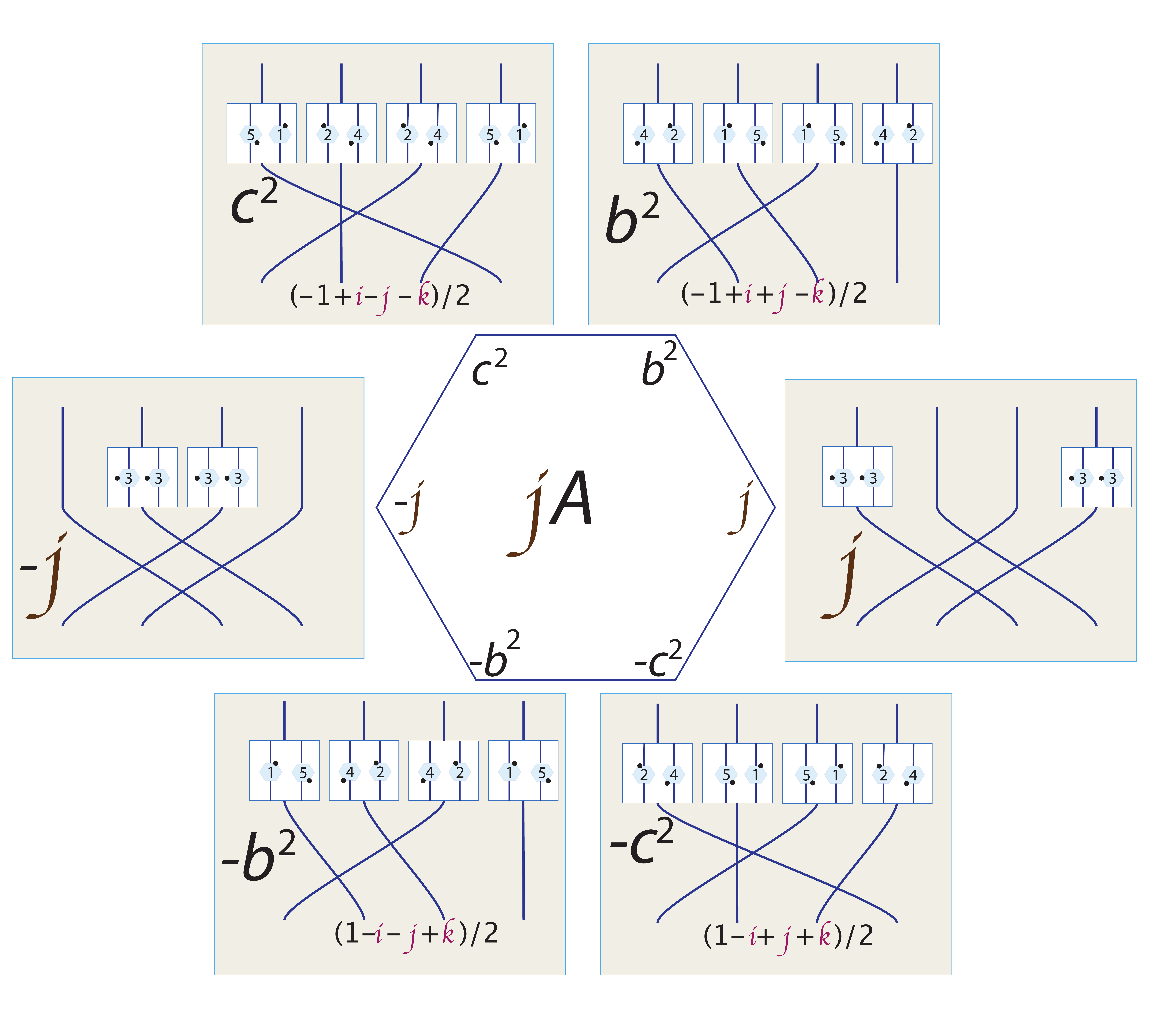}

Figure 57. The coset ${\boldsymbol j}A$
\end{center}

\newpage

   \begin{center}
\includegraphics[width=0.65\paperwidth]{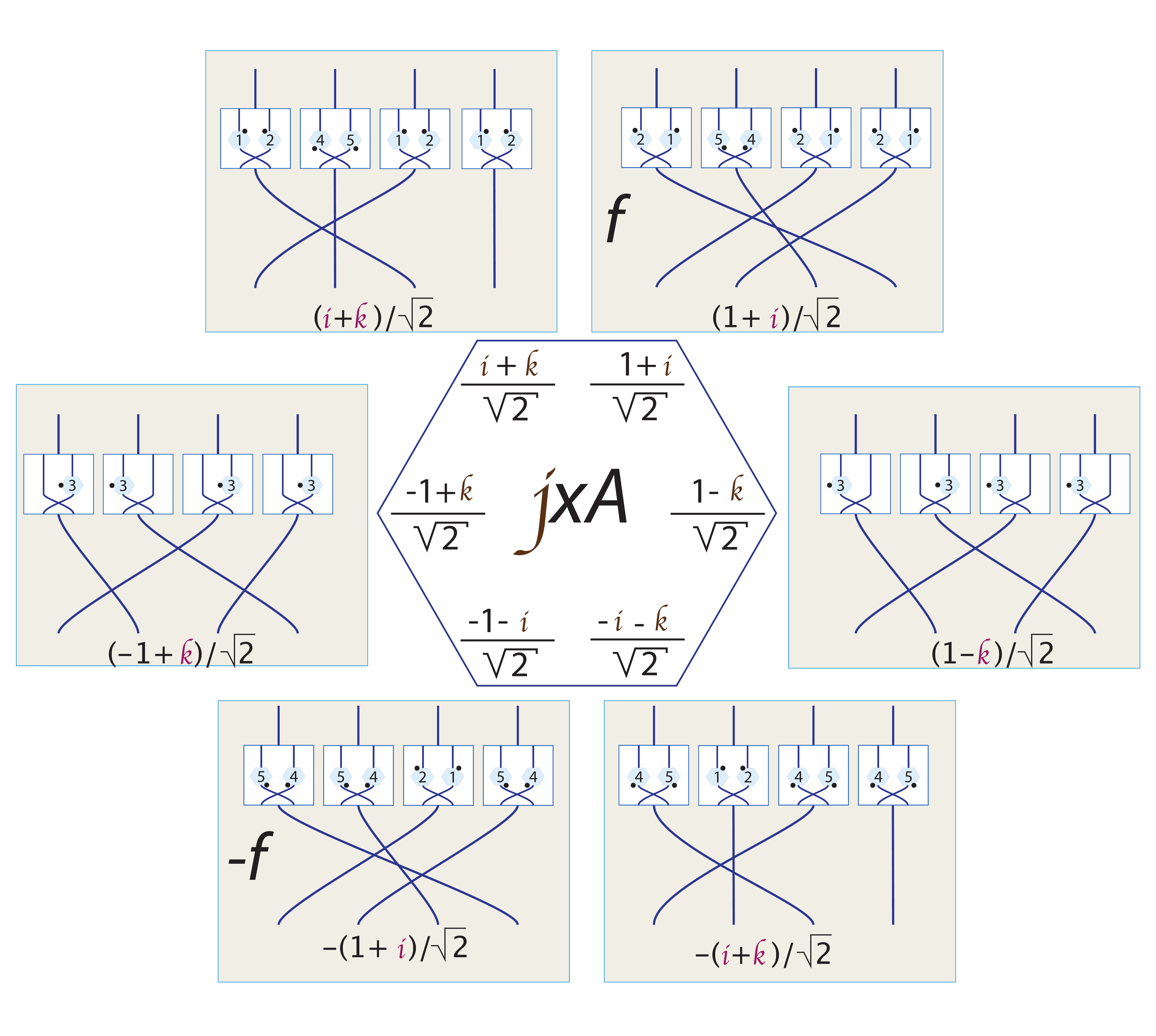}

Figure 58. The coset $x{\boldsymbol j}A$
\end{center}

\newpage

     \begin{center}
\includegraphics[width=0.65\paperwidth]{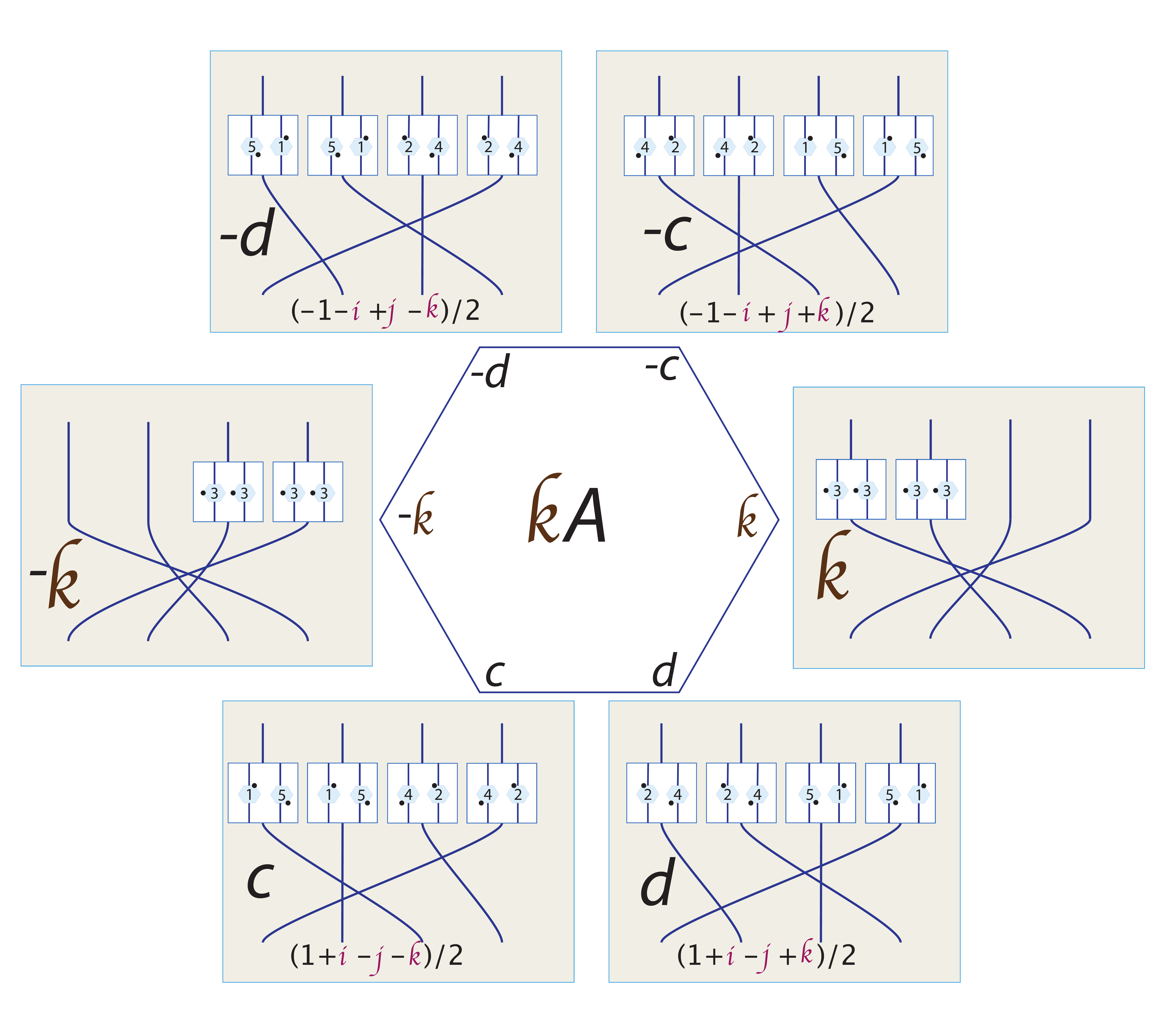}

Figure 59. The coset ${\boldsymbol k}A$
\end{center}

\newpage

        \begin{center}
\includegraphics[width=0.65\paperwidth]{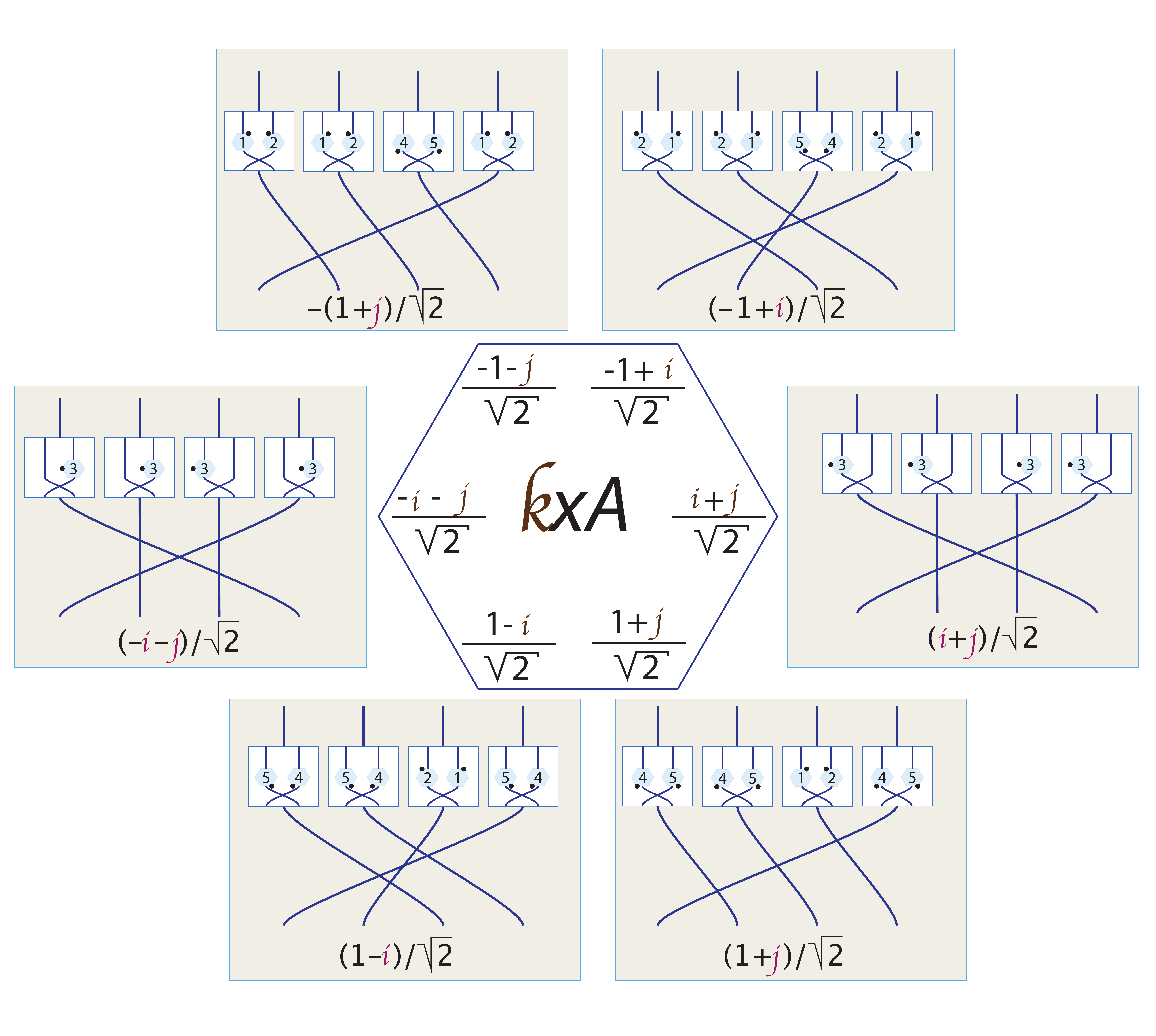}
Figure 60. The coset $x{\boldsymbol k}A$
\end{center}

\newpage

\section{The binary icosahedral group}
\label{S:BinIco}

The binary icosahedral group is given by the presentation:
$$\widetilde{A_5}= \langle a,t: (at)^2=a^3=t^5 \rangle.$$
It is a $2$-fold extension of the alternating group on $5$-elements, and the generators correspond to the  elements 
$a=1/2 (1 + {\boldsymbol i} + {\boldsymbol j}+ {\boldsymbol k}) $ and $t=1/2( \phi +(\phi -1){\boldsymbol i}+{\boldsymbol j})$ $\in S^3$ where $\phi=\frac{1+\sqrt{5}}{2}$ denotes the golden ratio. We remind the reader that $1/\phi=\phi-1$; also $\phi^2=\phi+1$, and $(1/\phi)^2=2-\phi.$  

Two pairs of diagrams will be used to describe the generators $t$ and $a$. One uses the cosets of the binary tetrahedral group   as well as its ordered decomposition into the  ordered cosets  $\widetilde{A_4} = (A, {\boldsymbol i}A, {\boldsymbol j}A, {\boldsymbol k}A)$ where $A=(1,a,a^2,-1,-a,-a^2)$ is the subgroup generated by $a$. The other uses the cosets of the subgroup $T=\langle t \rangle$ that is cyclic of order $10$. In both cases, figures will be used to demonstrate that $(at)^2=a^3=t^5=-1$.

Write
\begin{eqnarray*}
 \left[0\right] & = & 
 \left[ A, {\boldsymbol i}A, {\boldsymbol j}A, {\boldsymbol k}A\right];  \\ 
 \left[1\right] & = &  
 = \left[ tA, t{\boldsymbol i}A, t{\boldsymbol j}A, t{\boldsymbol k}A\right];  \\
\left[2\right] & = &  
 = \left[ t^2A, t^2{\boldsymbol i}A, t^2{\boldsymbol j}A, t^2{\boldsymbol k}A\right]; \\
 \left[3\right] & = &  
 \left[ t^3A, t^3{\boldsymbol i}A, t^3{\boldsymbol j}A, t^3{\boldsymbol k}A\right];  \\
\left[4\right] & = & 
\left[ t^4A, t^4{\boldsymbol i}A, t^4{\boldsymbol j}A, t^4{\boldsymbol k}A\right]. 
\end{eqnarray*}
Observe that $t[0]=[1], t[1]=[2],$ and $t[3]=[4]$. Since $t^5=-1$, we obtain that $t^5A=A$, but $A$ is rotated three steps. So in mimicing previous notation, we have $t(t^4A)=[A,3/6]$, and similarly  $t(t^4{\boldsymbol i}A)=[{\boldsymbol i}A,3/6]$, $t(t^4{\boldsymbol j}A)=[{\boldsymbol j}A,3/6]$,  and $t(t^4{\boldsymbol k}A)=[{\boldsymbol k}A,3/6].$

Figure~61 illustrates both the actions of $a$ and of $t$ on the cosets $[0]$ through $[4]$ by means of the thick strings. For example, $a$ corresponds to the cyclic permutation $(142)$ on the elements $\{[0],[1],[2],[3],[4] \}$ (where brackets and commas have been omitted from the expression of the $3$-cycle). At the top of each of these thick strings are the corresponding elements in $\widetilde{A_4}$ depicted using their actions on the cosets, $(A, {\boldsymbol i}A, {\boldsymbol j}A, {\boldsymbol k}A)$. From left-to-right in the illustration these correspond to the elements $a, {\boldsymbol j}, b^2, b, -a^2$. 

The details that go into creating these illustrations are analogous to those of the previous sections. The elements $a$ and $t$ permute the twenty cosets that constitute $[0]$ through $[4]$. The permutation action of $t$ is straight-forward, and the computation of the action of $a$ is facilitated by computing the products $at^n{\boldsymbol x}$ for $n=0,\ldots, 4$ and ${\boldsymbol x} = {\boldsymbol i}, {\boldsymbol j},$ or ${\boldsymbol k}.$ Then the location of that entry in the corresponding coset is determined. Rather than given any further outline, we will demonstrate that the strings-with-beads satisfy the relations that define $\widetilde{A_5}$.

Figure~62 indicates the composition $a^3$. At the bottom of  Fig.~63, the elements upon each of the five thick strings are composed, and then the elements on the 20 thinner strings are compiled  at the top of Fig.~63. Along the five thick strings of the illustration of Fig.~62, the gray rectangles indicate the products $a^3$, ${\boldsymbol j}b^2(-a^2)$, $b^2(-a^2){\boldsymbol j}$, $b^3$, and $(-a^2){\boldsymbol j}b^2.$ The representations of these elements coincides with those given in Fig.~32. Each product results in a representation of $(-1)$ in $\widetilde{A_4}$. That fact can be seen in Fig.~63. 
At the top of the figure the integral sum of the rotations is illustrated. Since $9$ and $15$ are congruent to $3$ modulo $6$, the identity $a^3=-1$ is demonstrated as the composition of the strings-with-beads.

   \begin{center}
\includegraphics[width=0.6\paperwidth]{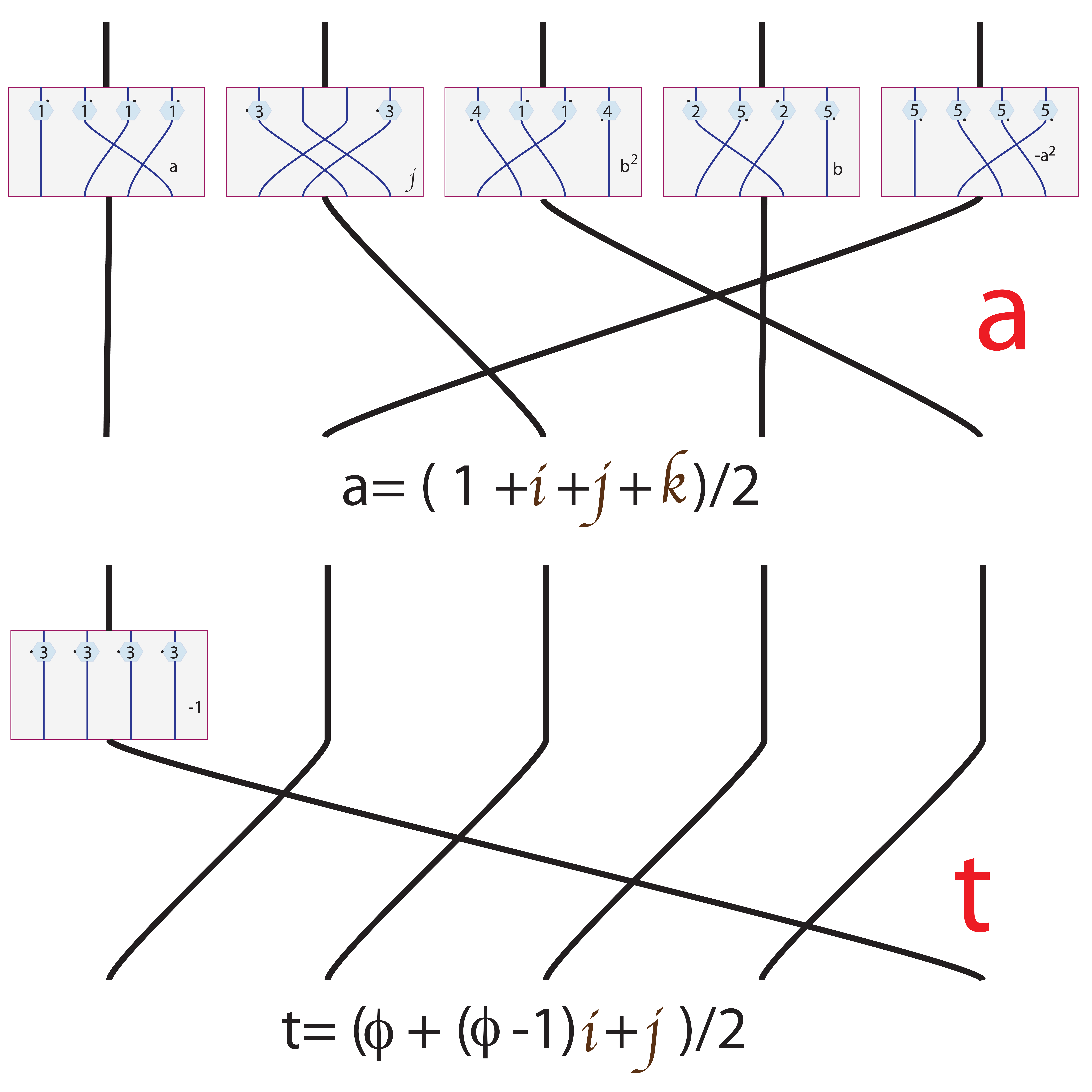}

Figure 61. Representing the elements $a$ and $t$ using cosets of the binary tetrahedral group
\end{center}

The identity $t^5=-1$ can be read from Fig.~64 directly. Each gray rectangle contains an illustration of $-1$ as quadruple of $180$ rotations on each of the cosets of $A$. That $t$ projects to the $5$-cycle $(01234)$ also is visually apparent. 

In Fig.~65, the composition $(at)^2$ is illustrated. In Fig.~66, the rectangular ``beads" are compiled and composed. Then the hexagonal beads on each of these $20$ strings are stacked, and the integral sum of $60$ degree twists is written at the top of the illustration.

   \begin{center}
\includegraphics[width=0.65\paperwidth]{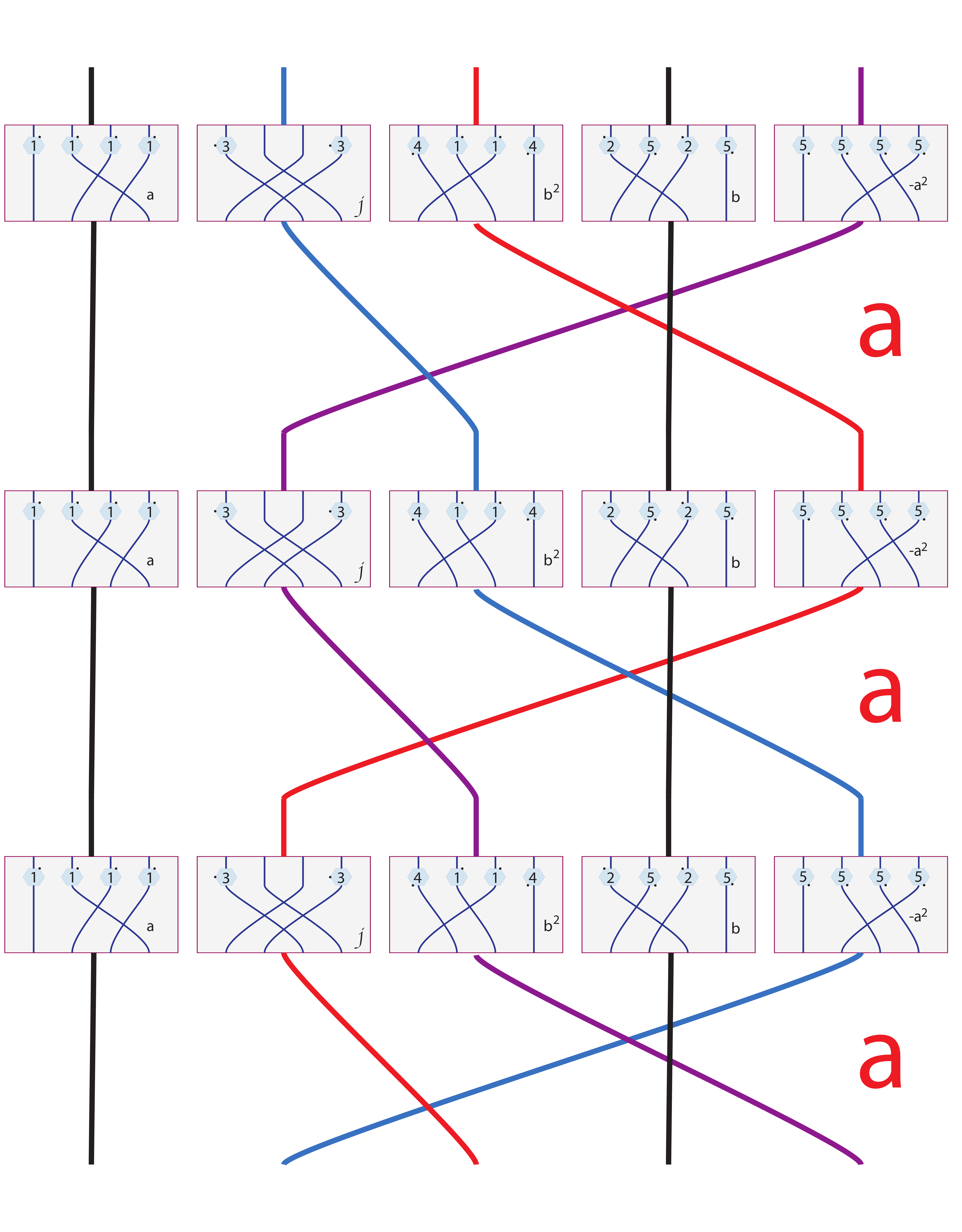}

Figure 62. The relation $a^3=-1$ in the binary icosahedral group
\end{center}

   \begin{center}
\includegraphics[width=0.65\paperwidth]{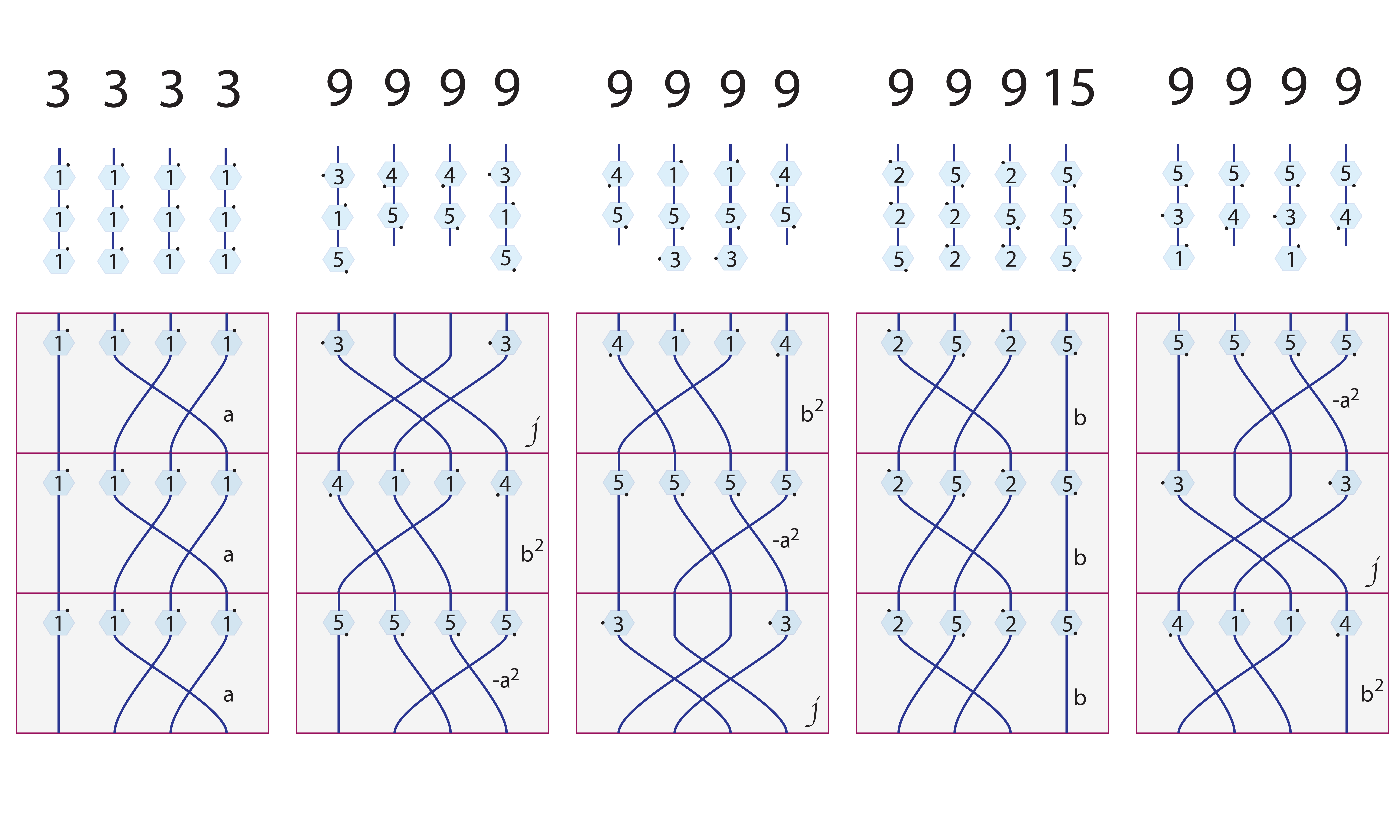}

Figure 63. The actions on the $20$ cosets of $A$
\end{center}

 Since  $a^3=-1=t^5$,  the order of the  generator $t$ divides $10$. Since $t^2= \left((\phi-1) +{\boldsymbol i} + \phi {\boldsymbol j}\right)/2$, the order of $t$ is exactly $10$. We compute the remaining powers of $t$ as follows. The reciprocals in $S^3$ are computed $(w + x{\boldsymbol i} + y{\boldsymbol j}+ z{\boldsymbol k})^{-1}= (w - x{\boldsymbol i} - y{\boldsymbol j}- z{\boldsymbol k})$.
Therefore, $t^9=\left(\phi + (1-\phi){\boldsymbol i}-{\boldsymbol j}\right)/2.$ We obtain $t^8=t^{-2}$, $t^7=t^5t^2=-t^2$,  $t^6=-t$, $t^4=t^{14}=-t^9$, and $t^3=t^{-7}$.

Let $T=\langle t \rangle$ denote the ordered subgroup of $\widetilde{A_5}$ that consists of the powers of $t$. Its cosets are ordered as \[(T, {\boldsymbol i}T, {\boldsymbol j}T, {\boldsymbol k}T, aT,d^2T, b^2T, cT, a^2T, bT, c^2T, dT).\]
This ordering was chosen because the action of $a$ preserves the groupings 
\[\left[ (T, {\boldsymbol i}T, {\boldsymbol j}T, {\boldsymbol k}T), (aT,d^2T, b^2T, cT),( a^2T, bT, c^2T, dT)\right].\]
The actions of $a$ and $t$ upon these cosets is illustrated in Fig.~67, and therein one can see that $a$ takes each quadruple in the grouping to the cyclically next quadruple. Unfortunately, the action of $t$ is not so well patterned. 
The results of some rather detailed and tedious calculations are also contained in these illustrations. The positive integers encircled by decagons indicate the amount of cyclic twisting induced by either $a$ or $t$ in the given cosets. The remaining Figs.~69-72 indicate that for these representatives the relations $a^3=t^5=(at)^2=-1$ hold. Figures~70 and~72 compile the total twists on each of the $12$ thick strings.  At the risk of either pedantry or having only one item that is completely understandable in the discussion, the sums along these $12$ strings are compiled. Each is congruent to $5$ modulo $10$. 

   \begin{center}
\includegraphics[width=0.4\paperwidth]{teatothefifth.pdf}

Figure 64. The relation $t^5=-1$ in the binary icosahedral group
\end{center}

   \begin{center}
\includegraphics[width=0.55\paperwidth]{EhTeaEhTeaA.pdf}

FIgure 65. The relation $(at)^2=-1$ in the binary icosahedral group
\end{center}

   \begin{center}
\includegraphics[width=0.6\paperwidth]{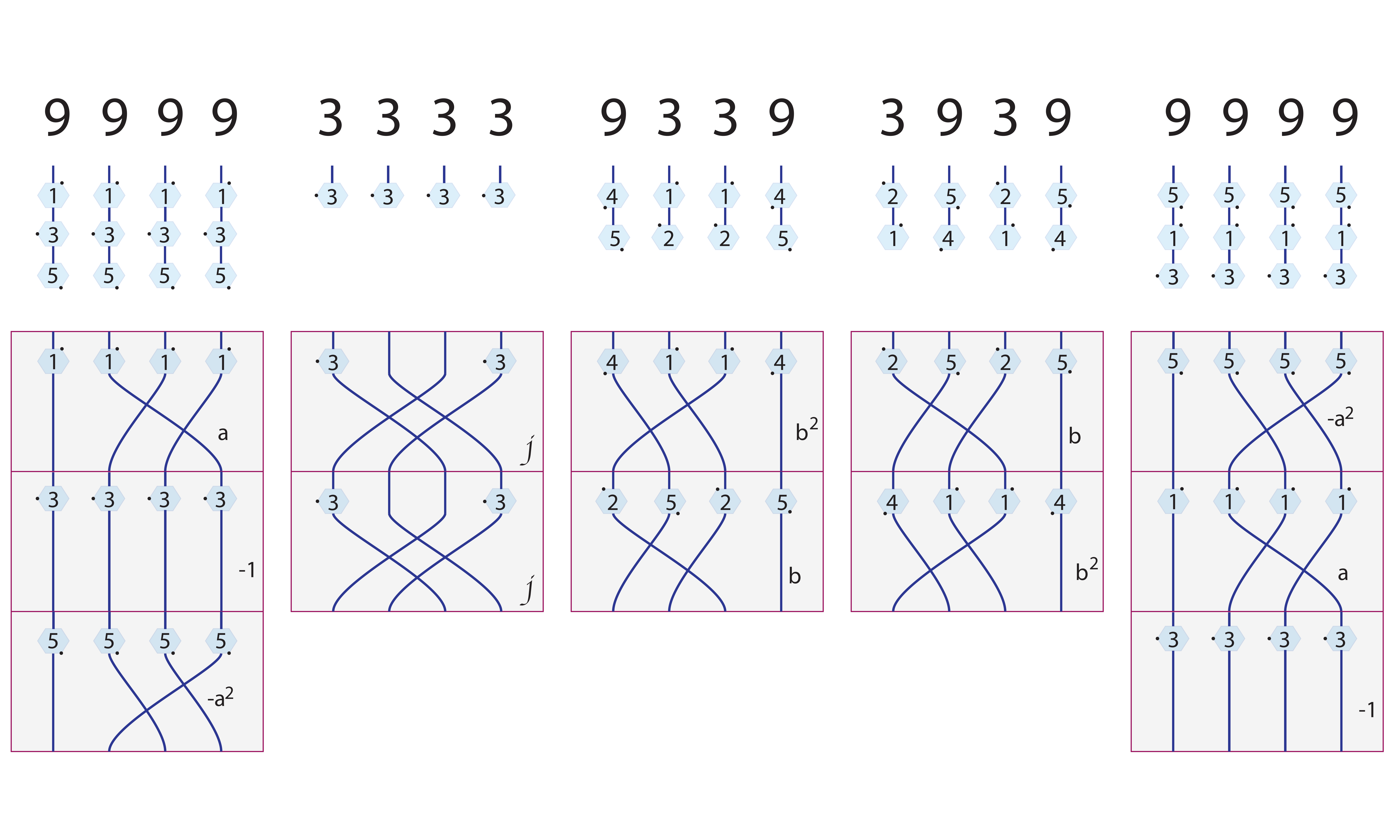}

Figure 66. Adding the total twists on the twenty strings
\end{center}

Items (12) and (13) of Theorem~\ref{main} are complete. And thus the proof of Theorem~\ref{main} is complete.

   \begin{center}
\includegraphics[width=0.6\paperwidth]{tanda2.pdf}

Figure 67. Representing the elements $a$ and $t$ using cosets of ${\mathcal T}=\langle t \rangle$
\end{center}

   \begin{center}
\includegraphics[width=0.6\paperwidth]{EhCubed2.pdf}

Figure 68. The relation $a^3=-1$ in the binary icosahedral group using cosets of ${\mathcal T}$
\end{center}

   \begin{center}
\includegraphics[width=0.45\paperwidth]{t2the5thA.pdf}

Figure 69. The relation $t^5=-1$ in the binary icosahedral group using cosets of ${\mathcal T}$
\end{center}

   \begin{center}
\includegraphics[width=0.3\paperwidth]{t2the5thB.pdf}

Figure 70. Compiling the twists  for the relation $t^5=-1$ in the cosets of $A=(1,a,a^2,-1,-a,-a^2)$
\end{center}

   \begin{center}
\includegraphics[width=0.5\paperwidth]{twelvestringPart1.pdf}

FIgure 71. The relation $(at)^2=-1$ in the binary icosahedral group using cosets of ${\mathcal T}$
\end{center}

   \begin{center}
\includegraphics[width=0.6\paperwidth]{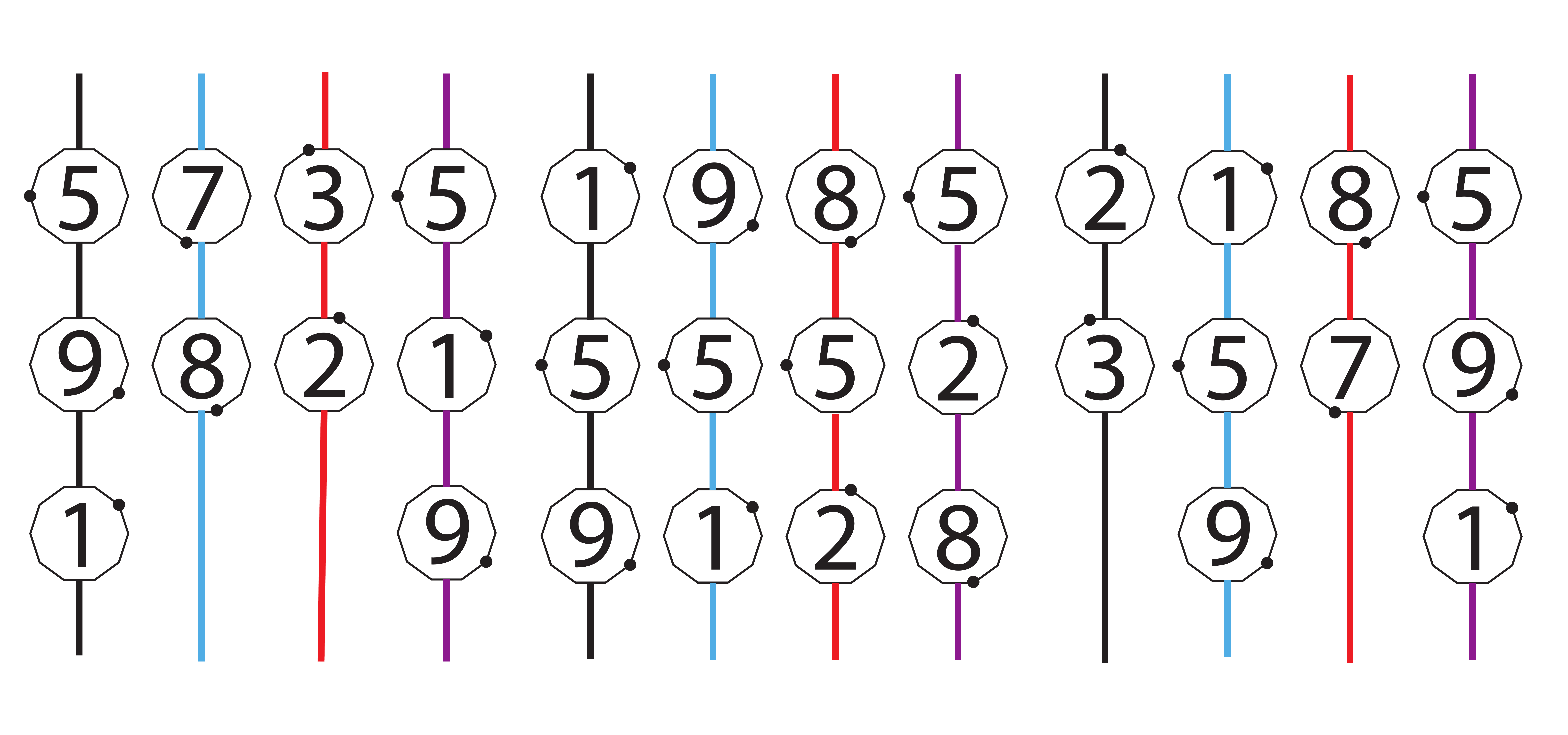}

FIgure 72. The beads on each string of $(at)^2=-1$
\end{center}

\section{Epilogue}

Semi-direct products of finite groups can be expressed via string diagrams that represent the subgroups and quotient groups as permutations acting upon themselves or upon cosets of particular subgroups. That statement is the context of Theorem~\ref{two}. In this way, string and ribbon diagrams can be obtained for the finite subgroups of $\SU(2)$ as well as many other group extensions. In this way novel diagrams for the quaternions and the binary polyhedral groups are obtained. These diagrams are analogues of braid diagrams. 
In many cases a complete list of such representatives has been presented. For other groups, we have expressed generators for the groups in permutation representations. These representations also can be expressed in matrix form in which each row and column has exactly one non-zero entry that is found as an element in an initially chosen subgroup. 

Our approach has been to be as meticulous as possible. However, more work can be done along these lines. In particular, since the binary icosahedral group is isomorphic to $\SL_2(\Z/5)$, its actions upon linear subspaces of $\Z/5 \times \Z/5$ might provide additional or more detailed permutation representations. Furthermore, quandle structures associated to these groups, subgroups, and conjugacy classes are also known to be interesting \cite{Inoue1,Inoue2}. Since conjugation in permutation groups is easily computed using cycle structures, it would be worthwhile to give algorithms to compute conjugation in these and other semi-direct products. 

The diagrams presented can easily be used to compute $2$-dimensional cocycles of the quotient groups that are represented by ignoring the beads upon the strings.

Many aspects of this paper are routine and tedious calculations that usually are left unwritten, or at best given to a student as exercises. Yet we feel strongly that a compilation of these calculations are both useful and, more importantly, entertaining.

\section*{Acknowledgements}
We gratefully acknowledge a grant in The Basic Science Research Program through the National Research Foundation of Korea(NRF) funded by the Ministry of Education(2016R1D1A3B01007669).
JSC was supported in by Simons Foundation Grant  381381 that initiated our conversations. He also enjoyed support from Kyungpook National University in 2016 and 2017. He visited Japan under a Japanese Society for the Promotion of Science grant numbered L18511. B.~Kim was supported by National Research Foundation of Korea(NRF) Grant No. 2019R1Q3B2067839.
 We would like to give personal thanks to Professor Sang Youl Lee, and Seiichi Kamada for several interesting conversations.

\bigskip
\footnotesize
\bibliographystyle{alpha}
\bibliography{biblio}
\end{document}